\documentclass[11pt]{amsart}
\usepackage{times}
\usepackage{eepic}
\usepackage[dvips]{graphicx,color}
\usepackage{graphics}
\usepackage{amssymb,amsfonts,amsmath,amsthm}
\usepackage{float}

\newtheorem{theorem}{Theorem}[section]
\newtheorem{proposition}[theorem]{Proposition}
\newtheorem{lemma}[theorem]{Lemma}

\newtheorem{remark}[theorem]{Remark}
\newtheorem{corollary}[theorem]{Corollary}
\newtheorem{example}[theorem]{Example}

\def\reff#1{{\rm(\ref{#1})}}

\newcommand\EE {\mathbb E}

\newcommand\RR {\mathbb R}

\newcommand\PP {\mathbb P}

\begin{document}

\title[Conservation Laws and FBSDEs]{Singular FBSDEs and Scalar Conservation Laws Driven by Diffusion Processes}

\author{Ren\'e Carmona}
\address{ORFE, Bendheim Center for Finance, Princeton University,
Princeton, NJ  08544, USA.}
\email{rcarmona@princeton.edu}
\thanks{Partially supported  by NSF: DMS-0806591}

\author{Fran\c{c}ois Delarue}
\address{Laboratoire J.-A. Dieudonn\'e, 
Universit\'e de Nice Sophia-Antipolis, 
Parc Valrose, 
06108 Cedex 02, Nice, FRANCE.}
\email{delarue@unice.fr}

\subjclass[2000]{Primary }


\date{August 25, 2011}

\begin{abstract}
Motivated by earlier work on the use of fully-coupled
Forward-Backward Stochastic Differential Equations (henceforth FBSDEs) in the analysis of mathematical models for the CO${}_2$ emissions markets, the present study is concerned with the analysis of these equations 
when the generator of the forward equation has a conservative degenerate structure and 
the terminal condition of the backward equation is a non-smooth function 
of the terminal value of the forward component. We show that a general form of existence and uniqueness result still holds. When the function giving the terminal condition is binary, we also show that the flow property of the forward component of the solution can fail at the terminal time. In particular, we prove that a Dirac point mass appears in its distribution, exactly at the location of the jump of the binary function giving the terminal condition. We provide a detailed analysis of the breakdown of the Markovian representation of the solution at the terminal time.
\end{abstract}

\maketitle

\section{\textbf{Introduction}}
\label{se:intro}

This paper is a contribution to the theory of fully-coupled forward-backward stochastic differential equations (FBSDEs), which goes back to the seminal paper of Pardoux and Peng \cite{PardouxPeng1992} on decoupled FBSDEs. Fundamental existence and uniqueness results for fully-coupled FBSDEs can be found in the subsequent works of Ma, Protter and Yong \cite{MaProtterYong},
Peng and Wu \cite{PengWu} and Delarue \cite{Delarue02}, or in the book of Ma and Yong \cite{MaYong.book}. All these results remain rather technical in nature. As 
far as we know, the terminal value $Y_T$ of the backward component is always specified in the Markovian set-up as a continuous function, say $Y_T=\phi(X_T)$, of the terminal value $X_T$ of the forward component. In this setting, a crucial role is played by the
representation $Y_t=v(t,X_t)$ of the backward component as a function of the forward component, the function $v$ being the natural candidate for 
the solution of a nonlinear Partial Differential Equation (PDE). As we shall see in this paper, the properties of this PDE play a crucial role in the analysis of Markovian FBSDEs.
See for example \cite{Delarue02, DelarueGuatteri}.
\vskip 2pt
Motivated by the analysis of  mathematical models of the CO${}_2$ emissions markets,  see for example \cite{CarmonaFehrHinz,Seifert,CarmonaFehrHinzPorchet,ChesneyTaschini,cdet}, we are interested in the case of forward processes $(X_t=(P_t,E_t))_{0\le t\le T}$ having a one-dimensional component $(E_t)_{0 \leq t \leq T}$ with bounded variations, and a backward component $(Y_t)_{0\le t\le T}$ having a terminal value given by a monotone function $\phi$ of $E_T$, and especially when $\phi$ is an indicator function of the form $\phi={\mathbf 1}_{[\Lambda,+\infty)}$. In \cite{cdet}, we proposed an unrealistic toy example for which we showed that, while the terminal condition 
could not be enforced, it was still possible to prove existence and uniqueness of a solution provided parts of the terminal condition are allowed to become part of the solution.

In the present paper, we propose general models for which we can prove a similar unique solvability result while, at the same time, we carefully investigate the pathological behavior of the solution near and at the threshold $\Lambda$. When  the forward diffusion process driving the model is non-degenerate, we show that there always exists a set of strictly positive probability of scenarii for which the degenerate component ends exactly at the threshold. In other words, the marginal law of the bounded variations component $(E_t)_{0 \leq t \leq T}$ has a Dirac point mass at the terminal time $T$. We also show that, conditionally on this event, the support of the distribution of $Y_T$ covers the entire interval  $[0,1]$.
This demonstrates a breakdown of the expected Markovian structure of the solution according to which $Y_T$ is expected to be a function of $E_T$.  
Restated in terms of the information $\sigma$-fields generated by the random variables $Y_T$ and $E_T$, this can be rewritten as
$\sigma(Y_T) \not \subset \sigma(E_T)$. This fact has dramatic consequences for the emissions market models. Indeed, while a price for the allowance certificates exists and is unique in such a model, 
its terminal value cannot be prescribed as the model would require! Finally, we investigate the formation of the Dirac mass at the threshold $\Lambda$ in some specific models: we study examples in which the noise degenerate at different rates. We also show that the toy model analyzed in \cite{cdet} is critical as we prove that the noise propagates near the threshold with a dramatically small variance.

\vskip 2pt
Our analysis is strongly connected to the standard theory of hyperbolic PDEs:  our FBSDE model appears as a stochastic perturbation of a  first-order conservation law. This connection plays a key role in the analysis: the relaxed notion of solvability and the pathological behavior of the solution at the terminal time are consequences of the shock observed in the corresponding conservation law. In particular, we spend a significant amount of time in determining how the system feels the deterministic shock. Specifically, we compare the intensity of the noise plugged into the system with the typical energy required to escape from the trap resulting from the shock: because of the degeneracy of the forward equation, it turns out that the noise plugged into the system is not strong enough to avoid the collateral effect of the shock. Put differently, the non-standard notion of solution follows from the concomitancy of the degeneracy of the forward equation and of the singularity of the terminal condition.

We are confident that the paper gives the right picture for the pathological behavior of the system at the terminal time. However, several questions about the dynamics remain open. In particular, it seems rather difficult to describe the propagation of the noise before terminal time. As we already pointed out, the specific models investigated in this paper highlight various forms of propagation, but all of them suggest that the way the noise spreads out is very sensitive to the properties of the coefficients, making it extremely difficult to establish general criteria for noise propagation.

\vskip 2pt
We now give more details about the FBSDEs we consider in this paper. We assume that the multivariate forward process $(X_t)_{0\le t\le T}$ has non-degenerate components which we denote $P_t$ and a degenerate component $E_t$. To be specific, our FBSDE has the form:
\begin{equation}
\label{eq:4:11:1}
\left\{
\begin{split}
&\displaystyle dP_t = b(P_t) dt + \sigma(P_t) dW_t,
\\
&\displaystyle dE_t = -f(P_t,Y_t)dt 
\\
&\displaystyle dY_t = \langle Z_t,dW_t\rangle, \quad 0 \leq t \leq T,
\end{split}
\right.
\end{equation}
As usual, the forward components have to satisfy an initial condition, say $(P_0,E_0)=(p,e)$ while the backward component needs to satisfy a terminal condition given by a function of the terminal value $(P_T,E_T)$ of the forward components. In this paper, we consider terminal conditions of the form $Y_T=\phi(E_T)$. Our goal is to study systems for which the function $\phi$ is singular, in contrast to most of the existing literature on FBSDEs in which the function $\phi$ is required to be Lipschitz continuous. In Section \ref{se:generalexistence} we prove a general existence and uniqueness result including the case of indicator functions $\phi$.
\vskip 2pt
It is demonstrated in \cite{cdet} that forward backward systems of the form \reff{eq:4:11:1} appear naturally in the analysis of mathematical models for the emissions markets.
See for example \cite{CarmonaFehrHinz,Seifert,CarmonaFehrHinzPorchet,ChesneyTaschini} for mathematical models of these markets.
As exemplified in \cite{cdet}, the components of the stochastic process $P$ can be viewed as the prices of the goods whose production is the source of the emissions, for example  electricity, $E$ as the cumulative emissions of the producers, and $Y$ as the price of an emission allowance typically representing one ton of CO${}_2$ equivalent. The terminal condition is of the form $\phi=\lambda {\mathbf 1}_{[\Lambda,+\infty)}$
where the real constant $\Lambda>0$ represents the \emph{cap}, emission target set up by the regulator. A penalty of $\lambda$ is applied to each ton of CO${}_2$ exceeding the threshold $\Lambda$ at the end $T$ of the regulation period. Changing \emph{num\'eraire} if necessary, we can assume that $\lambda=1$ without any loss of generality.   A typical example of function $f$ is given by $f(p,y) = \tilde c(p-y e_c)$ where $\tilde c$ is the inverse of the function giving the marginal costs of production of the goods and  $e_c$ the vector (with the same dimension as $p$) giving the rates of emissions associated with the production of the various goods. It is important to emphasize that the terminal condition
$Y_T={\mathbf 1}_{[\Lambda,+\infty)}(E_T)$ is given by a non-smooth deterministic function $ {\mathbf 1}_{[\Lambda,+\infty)}$ of the forward component of the system. A very particular case of \reff{eq:4:11:1} corresponding  to $d=1$, $b\equiv 0$, $\sigma\equiv 1$ and $f(p,y) = p-y$, was partially analyzed in \cite{cdet}, and served as motivation for the detailed analysis presented in this paper.

\vskip 2pt
The paper is organized as follows. Section \ref{se:generalexistence} gives our general existence and uniqueness result including analytic a priori estimates which are crucial for the subsequent analysis. In Section \ref{se:dirac}  we restrict ourselves to a terminal condition of the form $ \phi={\mathbf 1}_{[\Lambda,+\infty)}$ and we investigate the terminal value of the process $(E_t)_{0 \leq t \leq T}$. We show that, under suitable assumptions, the law of $E_T$ has a Dirac mass at the threshold $\Lambda$
for some initial conditions of the process $(P_t,E_t)_{0 \leq t \leq T}$. 
Under the additional assumption that the dynamics of the process $(P_t)_{0 \leq t \leq T}$ are uniformly elliptic, we show that this Dirac mass is present for all the initial conditions and that the support of the conditional law of $Y_T$ given $E_T=\Lambda$ is equal to the entire interval $[0,1]$. 
Finally, Section \ref{se:ac} is concerned with the analysis of the smoothness properties of the distribution of $E_t$ for $t<T$. We show how the paths of the process $E$ coalesce at the threshold $\Lambda$, and how the flow $e\hookrightarrow E^{0,p,e}$ looses its homeomorphic property at time $T$.

The paper ends with an appendix devoted to the proofs of technical estimates which would have distracted from the main thrust of the paper, should have they been included where used.

\section{\textbf{A General Existence and Uniqueness Theorem for Singular FBSDEs}}
\label{se:generalexistence}

We here discuss the existence and uniqueness of solutions to FBSDEs of the form
\reff{eq:4:11:1}, with $Y_T = \phi(E_T)$ as terminal condition at a given terminal time $T>0$ for a non-decreasing bounded function $\phi : \RR \hookrightarrow \RR$. Without any loss of generality we shall assume that
\begin{equation}
\label{eq:25:7:1}
\inf_{x\in\RR} \phi(x)=0,\qquad\text{ and }\qquad \sup_{x\in\RR}\phi(x)=1.
\end{equation}
The process $(P_t)_{0 \leq t \leq T}$ is of dimension $d$, while $(E_t)_{0 \leq t \leq T}$ and $(Y_t)_{0 \leq t \leq T}$ are one-dimensional. The process $W=(W_t)_{0 \leq t \leq T}$ is a $d$-dimensional Brownian motion on some complete probability space $(\Omega,{\mathcal F},\PP)$, and the filtration $({\mathcal F}_t)_{0 \leq t \leq T}$ is the filtration generated by $W$ augmented with $\PP$-null sets. (See also\footnote{\label{foo:25:7:1} When specified below, we will sometimes enlarge the filtration to carry additional randomness.}.) 

\vskip 6pt\noindent
Throughout the paper, the coefficients $b : \RR^d \hookrightarrow \RR^d$, 
$\sigma : \RR^{d} \hookrightarrow \RR^{d \times d}$ and 
$f:\RR^d \times \RR \hookrightarrow \RR$ satisfy 

({\bf A}) There exist three constants $L \geq 1$ and $\ell_1,\ell_2>0$, $1/L \leq \ell_1 \leq \ell_2 \leq L$, such that

(A.1) $b$ and $\sigma$ are at most of $L$-linear growth in the sense that
$$
|b(p)|+|\sigma(p)|\le L(1+|p|), \qquad\qquad p\in\RR^d
$$

and  $L$-Lipschitz continuous in the sense that:
$$
|b(p)-b(p')|+|\sigma(p)-\sigma(p')|\le L|p-p'|, \qquad\qquad p,p'\in\RR^d.
$$

(A.2) for any $y \in {\mathbb R}$, the function $\RR^d\ni p \hookrightarrow f(p,y)$ is $L$-Lipschitz continuous and satisfies:
$$
|f(p,y)|\le  L(1+|p|+|y|), \qquad\qquad p \in \RR^d.
$$ 

Moreover, for any $p \in \RR^d$, the function $y \mapsto f(p,y)$ is strictly increasing and satisfies
$$
\ell_1 |y-y'|^2 \leq (y-y')[f(p,y)-f(p,y')] \leq \ell_2 |y-y'|^2, \qquad\qquad y,y' \in \RR.
$$

Since $\ell_2 \leq L$, $f$ is also $L$-Lipschitz continuous in $y$.

\begin{remark}
The strict monotonicity property of $f$ must be understood as a strict convexity property of the anti-derivative of $f$, as typically assumed in the theory of scalar conservation laws.
\end{remark}

\vskip 2pt\noindent
The main result of this section is stated in terms of the left continuous and right continuous versions $\phi_-$ and $\phi_+$ of the function $\phi$ giving the terminal condition. They are defined as:
\begin{equation}
\label{fo:phi-+}
\phi_-(x)=\sup_{x'<x}\phi(x'),\qquad\text{ and }\qquad \phi_+(x)=\inf_{x'>x}\phi(x').
\end{equation}
\begin{theorem}
\label{thm:4:11:1}
Given any initial condition $(p,e) \in \RR^d \times \RR$, there exists a unique progressively-measurable 4-tuple
$(P_t,E_t,Y_t,Z_t)_{0 \leq t \leq T}$ satisfying \reff{eq:4:11:1} together with the initial conditions $P_0=p$ and $E_0=e$ and the terminal condition
\begin{equation}
\label{eq:4:11:2}
{\mathbb P}\bigl\{ \phi_-(E_T) \leq Y_T \leq \phi_+(E_T) \bigr\}=1.
\end{equation}
Moreover, there exists a constant $C$, depending on $L$ and $T$ only, such that almost surely
$|Z_t| \leq C$ for $t \in [0,T]$.
\end{theorem}

\begin{remark}
We first emphasize that Theorem \ref{thm:4:11:1} is still valid when Eq. \eqref{eq:4:11:1} is set on an interval 
$[t_0,T]$, $0 \leq t_0 <T$, with $(P_{t_0},E_{t_0})=(p,e)$ as initial condition. The solution is then denoted
by $(P_t^{t_0,p},E_t^{t_0,p,e},Y_t^{t_0,p,e},Z_t^{t_0,p,e})_{t_0 \leq t \leq T}$.

We now discuss the meaning of
Theorem \ref{thm:4:11:1}: generally speaking, it must be understood as an existence and uniqueness theorem with a relaxed terminal condition. As we shall see below, there is no way to construct a solution to \eqref{eq:4:11:1} satisfying the terminal condition $Y_T = \phi(E_T)$ exactly: relaxing the terminal condition is thus necessary to obtain an existence property. In this framework, \eqref{eq:4:11:2} must be seen as the right relaxed terminal condition since it guarantees both existence and uniqueness of a solution to \eqref{eq:4:11:1}. As explained in Introduction, the need for a relaxed terminal condition follows from the simultaneity of the degeneracy of the forward equation and of the singularity of the terminal condition. (We will go back to this point later.)

Eq. \eqref{eq:4:11:2} raises a natural question: is there any chance that $\phi_-(E_T)$ and 
$\phi_+(E_T)$ really differ with a non-zero probability? Put differently, does the process $(E_t)_{0 \leq t \leq T}$ really feel the discontinuity points of the terminal condition? (Otherwise, there is no need for a relaxed terminal condition.) The answer is given in Section \ref{se:dirac} in the case when $\phi$ is the Heaviside function: under some additional conditions, it is proven that the random variable $E_T$ has a Dirac mass at the singular point of the terminal condition; that is \eqref{eq:4:11:2} is relevant.
\end{remark}

\subsection{\textbf{Existence via mollifying}}

The analysis relies on a mollifying argument. We first handle the case when the function $\phi$ giving the terminal condition is a non-decreasing smooth function
(in which case $\phi_-=\phi_+=\phi$) with values in $[0,1]$: we then prove that Eq. \eqref{eq:4:11:1}
admits a unique solution $(P_t^{t_0,p},E_t^{\phi,t_0,p,e},Y_t^{\phi,t_0,p,e},Z_t^{\phi,t_0,p,e})_{t_0 \leq t \leq T}$ for any initial condition $(P_{t_0}^{t_0,p},E_{t_0}^{\phi,t_0,p,e})=(p,e)$. This permits to define the value function $v^{\phi} : [0,T] \times \RR^d \times \RR \ni (t_0,p,e) \hookrightarrow Y_{t_0}^{\phi,t_0,p,e}$ when $\phi$ is smooth. In a second step, we approximate the \emph{true} terminal condition in Theorem \ref{thm:4:11:1}
by a sequence $(\phi_n)_{n \geq 1}$ of smooth terminal conditions: we prove that the value functions
$(v^{\phi_n})_{n \geq 1}$ are uniformly continuous on every compact subset of $[0,T) \times \RR^d \times \RR$. By a compactness argument, we then complete the proof of existence in Theorem \ref{thm:4:11:1}.

\subsubsection{Existence and Uniqueness for a Smooth Terminal Condition}
We here assume that the terminal condition $\phi$ is a non-decreasing smooth function with values in $[0,1]$. (In particular, it is assumed to be Lipschitz.) In such a case, existence and uniqueness are known to hold in small time, see for example Delarue
\cite{Delarue02}. A standard way to establish existence and uniqueness over an interval of arbitrary length consists in forcing the system by non-degenerate noise. Below, we thus consider the case when the dynamics of $P$ and $E$ include additional noise terms of small variance
$\varepsilon^2 \in (0,1)$, $\varepsilon >0$.  To be more specific, we call \emph{mollified equation} the system
\begin{equation}
\label{eq:4:11:5}
\begin{cases}
&dP_t^{\varepsilon} = b(P_t^{\varepsilon}) dt + \sigma(P_t^{\varepsilon}) dW_t + \varepsilon dW_t',
\\
&dE_t^{\varepsilon} = -f(P_t^{\varepsilon},Y_t^{\varepsilon})dt + \varepsilon dB_t,
\\
&dY_t^{\varepsilon} = \langle Z_t^{\varepsilon}, dW_t \rangle +  \langle Z_t^{\prime,\varepsilon}, dW_t'\rangle +  \Upsilon_t^{\varepsilon} dB_t, \quad 0 \leq t \leq T,
\end{cases}
\end{equation}
with $Y_T^{\varepsilon} = \phi(E_T^{\varepsilon})$ as terminal condition, $(W_t)_{0 \leq t \leq T}$, $(W_t')_{0 \leq t \leq T}$ and $(B_t)_{0 \leq t \leq T}$
standing for three independent Brownian motions, of dimension $d$, $d$ and $1$ respectively, on the same probability space as above, the filtration being enlarged to accommodate them. Here, $(Z_t^{\varepsilon},Z_t^{\prime,\varepsilon},\Upsilon_t^{\varepsilon})_{0 \leq t \leq T}$ stands for the integrand of the martingale representation of $(Y_t^{\varepsilon})_{0 \leq t \leq T}$ in the complete filtration 
generated by $(W_t,W_t',B_t)_{0 \leq t \leq T}$.

When the coefficients $b$ and $\sigma$ are bounded and the coefficient $f$ is bounded w.r.t. $p$, the mollified system \eqref{eq:4:11:5} admits a unique solution for any initial condition, see \cite{Delarue02}. For any initial condition $(t_0,p,e) \in [0,T] \times \RR^d \times \RR$, we then denote by $(Y_{t}^{\varepsilon,t_0,p,e} )_{t_0\le t\le T}$ the unique solution of the mollified equation \reff{eq:4:11:5}.
Then, the value function 
\begin{equation}
\label{eq:25:7:5}
v^{\varepsilon} : (t_0,p,e) \hookrightarrow 
Y_{t_0}^{\varepsilon,t_0,p,e} \in [0,1]
\end{equation}
is of class ${\mathcal C}^{1,2}$ on $[0,T] \times \RR^d \times \RR$, with bounded and continuous derivatives, and satisfies the PDE:
\begin{equation}
\label{eq:7:11:10}
 \bigl[ \partial_t v^{\varepsilon} + {\mathcal L}_p v^{\varepsilon}
+ \frac{\varepsilon^2}{2} \partial^2_{pp} v^{\varepsilon}
+ \frac{\varepsilon^2}{2} \partial^2_{ee} v^{\varepsilon} \bigr](t,p,e)
- f\bigl(p,v^{\varepsilon}(t,p,e)\bigr) \partial_e v^{\varepsilon}(t,p,e) = 0,
\end{equation} 
with $v^{\varepsilon}(T,p,e)=\phi(e)$ as terminal condition, ${\mathcal L}_p$  standing for the second order partial differential linear operator
\begin{equation}
\label{eq:3:8:2}
{\mathcal L}_p = \langle b(\cdot), \partial_p \rangle
+ \frac{1}{2} {\rm Trace} \bigl[a(\cdot) \partial^2_{pp} \bigr], 
\end{equation}
where we use the notation $a$ for the symmetric non-negative definite matrix $a = \sigma \sigma^{\top}$.
The solution to \reff{eq:4:11:5} then satisfies: 
for $0 \leq t \leq T$, $Y_t^{\varepsilon} = 
v^{\varepsilon}(t,P_t^{\varepsilon},E_t^{\varepsilon})$, and
$$
Z_t^{\varepsilon}=  \sigma^{\top}(P_t^{\varepsilon})
\partial_p v^{\varepsilon}(t,P_t^{\varepsilon},E_t^{\varepsilon}), \quad Z_t^{\prime,\varepsilon}=\varepsilon
\partial_p v^{\varepsilon}(t,P_t^{\varepsilon},E_t^{\varepsilon})\quad \text{and}\quad \Upsilon_t^{\varepsilon} = \varepsilon
\partial_e v^{\varepsilon}(t,P_t^{\varepsilon},E_t^{\varepsilon}).
$$

The strategy for the proof of existence and uniqueness to Eq. \eqref{eq:4:11:1} when driven by the smooth terminal condition $\phi$ is as follows. In Propositions \ref{prop:4:11:1} and \ref{prop:20:07:2} below, we establish several crucial \emph{a-priori} estimates on $v^{\varepsilon}$. In particular, we prove that the gradient of $v^{\varepsilon}$ can be bounded on the whole $[0,T] \times \RR^d \times \RR$ by a constant only depending on the Lipschitz constant of $\phi$ and on $L$ in (A.1-2). 
Referring to the induction scheme in \cite{Delarue02} and using a standard truncation argument of the coefficients, this implies that Eq. \eqref{eq:4:11:1} is uniquely solvable without any boundedness assumption on 
$b$, $\sigma$ and $f$ and without any additional viscosity. (See Corollary \ref{cor:25:7:1}.) 

The first \emph{a-priori} estimate is similar to Proposition 3 in \cite{cdet}:

\begin{proposition}
\label{prop:4:11:1}
Assume that the coefficients $b$, $\sigma$ and $f$ are bounded in $p$.
Then, for the mollified equation \reff{eq:4:11:5}, we have:
\begin{equation}
\label{eq:9:11:1}
\forall (t,p,e) \in [0,T) \times \RR^d \times \RR, \quad 0 \leq \partial_e v^{\varepsilon}(t,p,e) \leq \frac{1}{\ell_1(T-t)}.
\end{equation}
Moreover, the $L^{\infty}$-norm  
of $\partial_e v^{\varepsilon}$ on the whole $[0,T] \times \RR^d \times \RR$ can be bounded
in terms of $L$ and the Lipschitz norm of $\phi$ only.
\end{proposition}

\begin{remark}
Pay attention that the bounds for $\partial_e v^{\varepsilon}$ are independent of the bounds of the coefficients $b$, $\sigma$ and $f$ w.r.t. the variable $p$. This point is crucial in the sequel.
\end{remark}

\noindent {\bf Proof.}
To prove \reff{eq:9:11:1}, without any loss of generality, we can assume that the coefficients 
are infinitely differentiable with bounded derivatives of any order. Indeed, if \reff{eq:9:11:1} holds in the
infinitely differentiable setting, it is shown to hold in the initial 
framework as well by a mollifying argument whose details may be found in Section 2 in \cite{Delarue02}. In that case, 
$L^{-1} \leq \partial_y f \leq L$. The point then consists in 
differentiating the processes $E^{\varepsilon}$ and $Y^{\varepsilon}$ with respect to the initial condition. 

If we fix an initial condition $(t_0,p,e) \in [0,T] \times \RR^d \times \RR$, then,  $(E^{\varepsilon,t_0,p,e}_t)_{t_0 \leq t \leq T}$ satisfies:
$$
E_t^{\varepsilon,t_0,p,e} = e- \int_{t_0}^t f(P_s^{\varepsilon,t_0,p},v^{\varepsilon}(s,E_s^{\varepsilon,t_0,p,e})) ds
$$ 
so that we can consider the derivative process $(\partial_e E_t^{\varepsilon,t_0,p,e})_{t_0 \leq t \leq T}$. By Pardoux and Peng \cite{PardouxPeng1992}, we can also consider
the four derivative processes $(\partial_e Y_t^{\varepsilon,t_0,p,e})_{t_0 \leq t \leq T}$,
$(\partial_e Z_t^{\varepsilon,t_0,p,e})_{t_0 \leq t \leq T}$, $(\partial_e Z_t^{\prime,\varepsilon,t_0,p,e})_{t_0 \leq t \leq T}$ and 
$(\partial_e \Upsilon_t^{\varepsilon,t_0,p,e})_{t_0 \leq t \leq T}$. 
Of course, $\partial_e Y_t^{\varepsilon,t_0,p,e} = \partial_e v^{\varepsilon}(t,P_t^{\varepsilon,t_0,p},E_t^{\varepsilon,t_0,p,e})
\partial_e E_t^{\varepsilon,t_0,p,e}$.
It is plain to see that (below, we do not specify the index $(t_0,p,e)$ to simplify the notations)
\begin{equation*}
d  \bigl[\partial_e E_t^{\varepsilon}\bigr]  = -\partial_y f\bigl(P_t^{\varepsilon},Y_t^{\varepsilon}\bigr) \partial_e Y_t^{\varepsilon} dt 
= - \partial_y f\bigl(P_t^{\varepsilon},Y_t^{\varepsilon}\bigr) \partial_e v^{\varepsilon}(t,P_t^{\varepsilon},E_t^{\varepsilon}) \partial_e E_t^{\varepsilon} dt, \quad t_0 \leq t \leq T,
\end{equation*}
so that
\begin{equation}
\label{eq:5:3:2}
\partial_e E_t^{\varepsilon}= \exp \biggl( - \int_{t_0}^t \partial_y f(P_s^{\varepsilon},Y_s^{\varepsilon}) \partial_e v^{\varepsilon}(s,P_s^{\varepsilon},E_s^{\varepsilon}) ds \biggr),
\quad t_0 \leq t \leq T,
\end{equation}
which is bounded from above and from below by positive constants, uniformly in time and randomness.
Now, we can compute
\begin{equation*}
d \bigl[ \partial_e Y_t^{\varepsilon} \bigr]
=  \langle \partial_e Z_t^{\varepsilon}, dW_t\rangle + \langle \partial_e Z_t^{\prime,\varepsilon}, dW_t'\rangle +
\partial_e \Upsilon_t^{\varepsilon} dB_t, \quad t_0 \leq t \leq T.
\end{equation*}
Taking expectations on both sides (note that 
$(\partial_e Z_t^{\varepsilon})_{t_0 \leq t \leq T}$, $(\partial_e Z_t^{\prime,\varepsilon})_{t_0 \leq t \leq T}$ and $(\partial_e \Upsilon_t^{\varepsilon})_{t_0 \leq t \leq T}$
belong to $L^2([t_0,T] \times \Omega,d{\mathbb P} \otimes dt)$),
we deduce that
$\partial_e Y_{t_0}^{\varepsilon} = {\mathbb E}[\partial_e Y_T^{\varepsilon}] = {\mathbb E}[\partial_e \phi(E_T^{\varepsilon}) \partial_e E_T^{\varepsilon}] \geq 0$, so that $\partial_e v^{\varepsilon}(t_0,p,e) \geq 0$.
To get the upper bound, we compute
\begin{equation*}
d \bigl[ \partial_e E_t^{\varepsilon} \bigr]^{-1}  = \partial_y f(P_t^{\varepsilon},Y_t^{\varepsilon}) \partial_e Y_t^{\varepsilon} \bigl[\partial_e E_t^{\varepsilon}]^{-2} dt, \quad t_0 \leq t \leq T,
\end{equation*}
so that
\begin{equation}
\label{eq:25:7:2}
\begin{split}
d \bigl[\partial_e Y_t^{\varepsilon} (\partial_e E_t^{\varepsilon})^{-1} \bigr] 
&= (\partial_e E_t^{\varepsilon})^{-1} \bigl[ \langle \partial_e Z_t^{\varepsilon}, dW_t \rangle + \langle \partial_e Z_t^{\prime,\varepsilon}, dW_t'\rangle +
\partial_e \Upsilon_t^{\varepsilon} dB_t \bigr]
\\
&\hspace{15pt} 
+ \partial_y f(P_t^{\varepsilon},Y_t^{\varepsilon}) \bigl[\partial_e Y_t^{\varepsilon}\bigr]^2 \bigl[\partial_e E_t^{\varepsilon}]^{-2} dt.
\end{split}
\end{equation}
Taking expectations on both sides 
and using the lower bound for $\partial_y f$,
we deduce that
\begin{equation}
\label{eq:25:7:3}
d \bigl( {\mathbb E} \bigl[\partial_e Y_t^{\varepsilon} (\partial_e E_t^{\varepsilon})^{-1} \bigr] \bigr)
\geq  \ell_1 {\mathbb E} \bigl( \bigl[\partial_e Y_t^{\varepsilon}\bigr]^2 \bigl[\partial_e E_t^{\varepsilon}]^{-2} \bigr) dt.
\end{equation}
By Cauchy-Schwarz inequality, we obtain
\begin{equation*}
d \bigl( {\mathbb E} \bigl[\partial_e Y_t^{\varepsilon} (\partial_e E_t^{\varepsilon})^{-1} \bigr] \bigr)
\geq  \ell_1 \bigl[ {\mathbb E} \bigl( \bigl[\partial_e Y_t^{\varepsilon}\bigr] \bigl[\partial_e E_t^{\varepsilon}]^{-1} \bigr) \bigr]^2 dt, \quad t_0 \leq t \leq T.
\end{equation*}
Without loss of generality, we can assume that $\partial_e v^{\varepsilon}(t_0,p,e) \not =0$, as otherwise, the upper bound for the derivative is obvious. Therefore, 
$\partial_e Y_{t_0}^{\varepsilon} (\partial_e E_{t_0}^{\varepsilon})^{-1} = \partial_e v^{\varepsilon}(t_0,p,e) \not =0$.
We then consider the first time $\tau$ at which the continuous deterministic   function
${\mathbb E} \bigl[\partial_e Y_t^{\varepsilon} (\partial_e E_t^{\varepsilon})^{-1} \bigr]$ hits zero. 
For $t \in [t_0,\tau \wedge T)$, we have
\begin{equation*}
\frac{d \bigl( {\mathbb E} \bigl[\partial_e Y_t^{\varepsilon} (\partial_e E_t^{\varepsilon})^{-1} \bigr] \bigr)}{
\bigl[ {\mathbb E} \bigl( \bigl[\partial_e Y_t^{\varepsilon}\bigr] \bigl[\partial_e E_t^{\varepsilon}]^{-1} \bigr) \bigr]^2}
\geq  \ell_1 dt,
\end{equation*}
i.e.
$\displaystyle \bigl[\partial_e Y_{t_0}^{\varepsilon} (\partial_e E_{t_0}^{\varepsilon})^{-1} \bigr]^{-1}
- \bigl({\mathbb E} \bigl[\partial_e Y_t^{\varepsilon} (\partial_e E_t^{\varepsilon})^{-1} \bigr] \bigr)^{-1}
\geq \ell_1 (t-t_0)$.

We then notice that the function $t \in [t_0,T] \mapsto
{\mathbb E} [\partial_e Y_t^{\varepsilon} (\partial_e E_t^{\varepsilon})^{-1}]$ cannot vanish, as otherwise, the left-hand side would be $-\infty$ since explosion can only occur in $+\infty$. Therefore, we can let $t$ tend to $T$ to deduce that
\begin{equation}
\label{eq:7:3:10}
\bigl[\partial_e Y_{t_0}^{\varepsilon} (\partial_e E_{t_0}^{\varepsilon})^{-1} \bigr]^{-1} \geq \ell_1 (T-t_0),
\end{equation}
which completes the first part of the proof. 

The second part of the proof is a straightforward consequence of Delarue \cite{Delarue02}: in small time, i.e. for $t$ close to $T$, the Lipschitz constant
of $v^{\varepsilon}(t,\cdot)$ can be bounded in terms of the Lipschitz constants of the coefficients and of the terminal condition. For $t$ away from $T$, the result follows from the first part of the statement directly. \qed

\vskip 2pt\noindent

We now estimate the derivative of the value function in the direction $p$, and then derive the regularity in time:
\begin{proposition}
\label{prop:20:07:2}
Assume that the coefficients $b$, $\sigma$ and $f$ 
in the mollified equation \reff{eq:4:11:5} are bounded in $p$. Then, there exists a constant 
$C$, depending on $L$ and $T$ only, such that, for any $(t,p,e) \in [0,T) \times \RR^d \times \RR$, we have:
\begin{equation}
\label{eq:20:07:3}
\bigl| \partial_p v^{\varepsilon}(t,p,e) \bigr| \leq C.
\end{equation}
As a consequence, for any $\delta \in (0,T)$ and any compact set $K \subset \RR^d$, the $1/2$-H\"older norm of the function $[0,T-\delta] \ni t\hookrightarrow v^{\varepsilon}(t,p,e)$, $p \in K$ and $e \in \RR$,
is bounded in terms of $\delta$, $K$, $L$ and $T$ only; the $1/2$-H\"older norm of the function $[0,T] \ni t\hookrightarrow v^{\varepsilon}(t,p,e)$ (that is the same function but on the whole $[0,T]$), $p \in K$ and $e \in \RR$, is bounded in terms of $K$, $L$, $T$ and the Lipschitz norm of $\phi$ only.
\end{proposition}

{\bf Proof.} Again, using a mollifying argument if necessary, we can assume without any loss of generality that the coefficients are infinitely differentiable with bounded derivatives of any order. By Section 3 in Crisan and Delarue \cite{CrisanDelarue}, $v^{\varepsilon}$ is infinitely differentiable w.r.t. $p$ and $e$, with bounded and continuous derivatives of any order on $[0,T] \times \RR^d \times \RR$. The idea then consists in differentiating the equation satisfied by $v^{\varepsilon}$ with respect to $p$. Writing $p=(p_i)_{1\le i\le d}$, we see that  
$(\partial_{p_i} v^{\varepsilon})_{1 \leq i \leq d}$ satisfies the system of PDEs:
\begin{equation}
\label{eq:13:03:1}
\begin{split}
& \bigl[ \partial_t \bigl( \partial_{p_i} v^{\varepsilon} \bigr) + 
{\mathcal L}_p \bigl( \partial_{p_i} v^{\varepsilon} \bigr)
+ \frac{\varepsilon^2}{2} \Delta_{pp} \bigl( \partial_{p_i} v^{\varepsilon} \bigr)
+ \frac{\varepsilon^2}{2} \partial^2_{ee} \bigl( \partial_{p_i} v^{\varepsilon} \bigr) \bigr](t,p,e)+ \langle \partial_{p_i} b(p),\partial_p v^{\varepsilon}(t,p,e) \rangle
\\
&\hspace{25pt} + \frac{1}{2} {\rm Trace} \bigl[ \partial_{p_i} a(p) \partial_{pp}^2 v^{\varepsilon}(t,p,e) \bigr] 
 - f \bigl(p,v^{\varepsilon}(t,p,e) \bigr) \partial_{e} \bigl[ \partial_{p_i} v^{\varepsilon} \bigr](t,p,e) 
\\
&\hspace{25pt} - \bigl[ \partial_{p_i} f\bigl(p,v^{\varepsilon}(t,p,e) \bigr)
+ \partial_y f \bigl(p,v^{\varepsilon}(t,p,e) \bigr) 
\partial_{p_i} v^{\varepsilon}(t,p,e) \bigr] \partial_{e} v^{\varepsilon}(t,p,e) = 0,
\end{split}
\end{equation}
for $(t,p,e) \in [0,T) \times \RR^d \times \RR$ and $1 \leq i \leq d$, with the terminal condition $\partial_{p_i} v^{\varepsilon}(T,p,e) = 0$.
For a given initial condition $(t_0,p,e) \in [0,T) \times \RR^d \times \RR$,
we set 
$$(U_t^{\varepsilon},V_t^{\varepsilon}) = (\partial_p v^{\varepsilon}(t,P_t^{\varepsilon,t_0,p},E_t^{\varepsilon,t_0,p,e}),\partial_{pp}^2 v^{\varepsilon}(t,P_t^{\varepsilon,t_0,p},E_t^{\varepsilon,t_0,p,e})), \quad t_0 \leq t \leq T.
$$ 
By \eqref{eq:13:03:1} and It\^o's formula, we can write
\begin{equation}
\label{eq:13:03:10}
\begin{split}
&dU_t^{\varepsilon} = - [\partial_p b(P_t^{\varepsilon})]^{\top} U_t^{\varepsilon} dt
- \frac{1}{2} {\rm Trace} \bigl[ \partial_{p} a(P_t^{\varepsilon}) V_t^{\varepsilon} \bigr] dt + \partial_{p} f\bigl(P_t^{\varepsilon},Y_t^{\varepsilon} \bigr)
\partial_{e} v^{\varepsilon}(t,P_t^{\varepsilon},E_t^{\varepsilon}) dt
\\
&\hspace{5pt} 
+ \partial_y f \bigl(P_t^{\varepsilon},Y_t^{\varepsilon}\bigr) 
\partial_{e} v^{\varepsilon}(t,P_t^{\varepsilon},E_t^{\varepsilon})
U_t^{\varepsilon} dt + V_t^{\varepsilon} \sigma(P_t^{\varepsilon}) dW_t
+ \varepsilon V_t^{\varepsilon} dW_t' + \varepsilon \partial_{pe}^2 v^{\varepsilon}(t,P_t^{\varepsilon},E_t^{\varepsilon}) dB_t, 
\end{split}
\end{equation}
$t_0 \leq t \leq T$, ${\rm Trace}[\partial_p a(P_t^{\varepsilon}) V_t^{\varepsilon}]$ standing for the vector
$({\rm Trace}[\partial_{p_i} a(P_t^{\varepsilon}) V_t^{\varepsilon}])_{1 \leq i \leq d}$. Introducing the exponential weight
\begin{equation*}
{\mathcal E}_t^{\varepsilon} = \exp \biggl( - \int_{t_0}^t 
\partial_y f \bigl(P_s^{\varepsilon},Y_s^{\varepsilon}\bigr) 
\partial_{e} v^{\varepsilon}(s,P_s^{\varepsilon},E_s^{\varepsilon})
ds\biggr), \quad t_0 \leq t \leq T,
\end{equation*}
and setting $(\bar{U}_t^{\varepsilon},\bar{V}_t^{\varepsilon})=
{\mathcal E}_t^{\varepsilon} (U_t^{\varepsilon},V_t^{\varepsilon})$, we obtain
\begin{equation}
\label{eq:13:03:2}
\begin{split}
d \bar{U}_t^{\varepsilon}  &= - [ \partial_p b(P_t^{\varepsilon})]^{\top} \bar{U}_t^{\varepsilon} dt
- \frac{1}{2} {\rm Trace} \bigl[ \partial_{p} a(P_t^{\varepsilon}) \bar{V}_t^{\varepsilon} \bigr] dt
 + \partial_{p} f\bigl(P_t^{\varepsilon},Y_t^{\varepsilon} \bigr)
\partial_{e} v^{\varepsilon}(t,P_t^{\varepsilon},E_t^{\varepsilon}) {\mathcal E}_t^{\varepsilon} dt
\\
&\hspace{15pt} + \bar{V}_t^{\varepsilon} \sigma(P_t^{\varepsilon}) dW_t
+ \varepsilon \bar{V}_t^{\varepsilon} dW_t' + \varepsilon {\mathcal E}_t^{\varepsilon} \partial_{pe}^2 v^{\varepsilon}(t,P_t^{\varepsilon},E_t^{\varepsilon}) dB_t, \quad t_0 \leq t \leq T,
\end{split}
\end{equation}
with $\bar{U}_T^{\varepsilon} = 0$ as terminal condition.
Noting that $|{\rm Trace}[\partial_{p} a(P_t^{\varepsilon})
\bar{V}_t^{\varepsilon}]|\leq L |\bar{V}_t^{\varepsilon} \sigma(P_t^{\varepsilon})|$
and applying It\^o's formula to $(|\bar{U}_t^{\varepsilon}|^2)_{t_0 \leq t \leq T}$, 
we obtain
\begin{equation*}
\begin{split}
&{\mathbb E} \bigl[ |\bar{U}_t^{\varepsilon}|^2 \bigr]
+ {\mathbb E} \int_t^T \bigl[ |\bar{V}_s^{\varepsilon}
\sigma(P_s^{\varepsilon})|^2 + \varepsilon^2  
|\bar{V}_s^{\varepsilon}|^2
+ \varepsilon^2 | {\mathcal E}_s^{\varepsilon} \partial_{pe}^2 v^{\varepsilon}(s,P_s^{\varepsilon},E_s^{\varepsilon})|^2 \bigr] ds
\\
&\leq C {\mathbb E} \int_t^T \bigl[ |\bar{U}_s^{\varepsilon}|^2 + 
|\bar{U}_s^{\varepsilon}| |\bar{V}_s^{\varepsilon} \sigma(P_t^{\varepsilon})|
\bigr] ds + C {\mathbb E} \biggl[ \sup_{t_0 \leq s \leq T}
|\bar{U}_s^{\varepsilon}| \int_t^T 
\partial_{e} v^{\varepsilon}(s,P_s^{\varepsilon},E_s^{\varepsilon}) {\mathcal E}_s^{\varepsilon} ds
\biggr],
\end{split}
\end{equation*}
for some constant $C$ depending on $L$ only, possibly varying from line to line. By (A.2), the integral 
$\int_t^T 
\partial_{e} v^{\varepsilon}(s,P_s^{\varepsilon},E_s^{\varepsilon}) {\mathcal E}_s^{\varepsilon} ds$ can be bounded by $L$. Using the convexity inequality $2xy \leq ax^2+a^{-1}y^2$, $a>0$, we deduce 
\begin{equation*}
\begin{split}
&{\mathbb E} \bigl[ |\bar{U}_t^{\varepsilon}|^2 \bigr]
+ \frac{1}{2} {\mathbb E} \int_t^T \bigl[ |\bar{V}_s^{\varepsilon}
\sigma(P_s^{\varepsilon})|^2 + \varepsilon^2  
|\bar{V}_s^{\varepsilon}|^2
+ \varepsilon^2 | {\mathcal E}_s^{\varepsilon} \partial_{pe}^2 v(s,P_s^{\varepsilon},E_s^{\varepsilon})|^2 \bigr] ds
\\
&\leq C {\mathbb E} \int_t^T |\bar{U}_s^{\varepsilon}|^2 ds + C {\mathbb E} \bigl[ \sup_{t_0 \leq s \leq T}
|\bar{U}_s^{\varepsilon}| \bigr].
\end{split}
\end{equation*}
By \eqref{eq:13:03:2} and the B\"urkholder-Davies-Gundy inequality, we can bound 
${\mathbb E} [ \sup_{t_0 \leq s \leq T}
|\bar{U}_s^{\varepsilon}| ]$ directly. We obtain
\begin{equation*}
\begin{split}
&{\mathbb E} \bigl[ |\bar{U}_t^{\varepsilon}|^2 \bigr]
+ \frac{1}{2} {\mathbb E} \int_t^T \bigl[ |\bar{V}_s^{\varepsilon}
\sigma(P_s^{\varepsilon})|^2 + \varepsilon^2  
|\bar{V}_s^{\varepsilon}|^2
+ \varepsilon^2 | {\mathcal E}_s^{\varepsilon} \partial_{pe}^2 v(s,P_s^{\varepsilon},E_s^{\varepsilon})|^2 \bigr] ds
\\
&\leq C \biggl[ 1 + {\mathbb E} \int_t^T \bigl[ |\bar{U}_s^{\varepsilon}|^2 +  |\bar{U}_s^{\varepsilon}| + |\bar{V}_s^{\varepsilon}\sigma(P_s^{\varepsilon})| \bigr] ds 
\\
&\hspace{15pt}
+ \biggl(  {\mathbb E} \int_t^T \bigl[ |\bar{V}_s^{\varepsilon}
\sigma(P_s^{\varepsilon})|^2 + \varepsilon^2  
|\bar{V}_s^{\varepsilon}|^2 
+ \varepsilon^2 | {\mathcal E}_s^{\varepsilon} \partial_{pe}^2 v(s,P_s^{\varepsilon},E_s^{\varepsilon})|^2 
\bigr] ds
 \biggr)^{1/2}
\biggr].
\end{split}
\end{equation*}
Using the convexity inequality $2xy \leq ax^2+a^{-1}y^2$ again (with $a>0$), together with Gronwall's Lemma, we deduce the announced bound for ${\mathbb E}[|U_{t_0}^{\varepsilon}|^2]$.

The end of the proof is quite standard. For $(t_0,t,p,e) \in [0,T) \times [0,T) \times \RR^d \times \RR$, $t_0 \leq t$, we deduce from the martingale property 
of $(v^{\varepsilon}(s,P_s^{\varepsilon,t_0,p},E_s^{\varepsilon,t_0,p,e}))_{t_0 \leq s \leq T}$:
\begin{equation*}
|v^{\varepsilon}(t_0,p,e) - v^{\varepsilon}(t,p,e)|  \leq
{\mathbb E} \bigl[ |v^{\varepsilon}(t,P_t^{\varepsilon,t_0,p},E_t^{\varepsilon,t_0,p,e}) - v(t,p,e)| \bigr].
\end{equation*}
By \eqref{eq:9:11:1} and \eqref{eq:20:07:3}, there exists a constant $C$, depending on $L$ and $T$ only, such that
\begin{equation*}
|v^{\varepsilon}(t_0,p,e) - v^{\varepsilon}(t,p,e)| 
\leq \frac{C}{T-t}  {\mathbb E} \bigl[ |P_t^{\varepsilon,t_0,p}-p| + 
|E_t^{\varepsilon,t_0,p,e}-e| \bigr]
\leq \frac{C}{T-t} (1+|p|) (t-t_0)^{1/2}.
\end{equation*}
Alternatively, by the second part in Proposition \ref{prop:4:11:1} and by \eqref{eq:20:07:3},
there exists a constant $C'$, depending on $L$, $T$ and the Lipschitz constant of $\phi$ only, such that
\begin{equation*}
|v^{\varepsilon}(t_0,p,e) - v^{\varepsilon}(t,p,e)|
\leq C'  {\mathbb E} \bigl[ |P_t^{\varepsilon,t_0,p}-p| + 
|E_t^{\varepsilon,t_0,p,e}-e| \bigr] \leq C'  (1+|p|) (t-t_0)^{1/2}.
\end{equation*} 
This completes the proof. \qed

\vskip 2pt\noindent

As announced, we deduce 
\begin{corollary}
 \label{cor:25:7:1}
Assume that Assumptions (A.1) and (A.2) are in force and that $\phi$ is a non-decreasing Lipschitz smooth function with values in $[0,1]$. Then, for any $(t_0,p,e) \in [0,T] \times \RR^d \times \RR$ and $\varepsilon \in (0,1)$, the mollified Eq. \eqref{eq:4:11:5} is uniquely solvable under the initial condition $(P^{\varepsilon}_{t_0},E^{\varepsilon}_{t_0})=(p,e)$ (even if the coefficients are not bounded). In particular, the value function $v^{\varepsilon}$ in \eqref{eq:25:7:5} still makes sense: it is a ${\mathcal C}^{1,2}$ solution 
of the PDE \eqref{eq:7:11:10} on the whole $[0,T] \times \RR^d \times \RR$ and it satisfies \eqref{eq:9:11:1} and 
\eqref{eq:20:07:3}.

Similarly, original Eq. \eqref{eq:4:11:1} is uniquely solvable under the initial condition $(P_{t_0},E_{t_0})=(p,e)$ and the terminal condition $Y_T=\phi(E_T)$. 
Denoting by 
$(P_t^{t_0,p},E_t^{\phi,t_0,p,e},Y_t^{\phi,t_0,p,e},Z_t^{\phi,t_0,p,e})_{t_0 \leq t \leq T}$ the solution with $(t_0,p,e)$ as initial condition, the value function $v^{\phi}: [0,T] \times \RR^d \times \RR \ni (t_0,p,e) \hookrightarrow 
Y_{t_0}^{\phi,t_0,p,e}$ is the limit of $v^{\varepsilon}$ as $\varepsilon$ tends to 0, the convergence being uniform on compact subsets of $[0,T] \times \RR^d \times \RR$.
\end{corollary}

{\bf Proof.} The proof of unique solvability is the same as in Delarue \cite{Delarue02}, both for \eqref{eq:4:11:5} and for \eqref{eq:4:11:1}. In both cases, the coefficients are Lipschitz continuous: existence and uniqueness hold in small time; this permits to define the value functions $v^{\varepsilon}$ and $v^{\phi}$ on some interval $[T-\delta_0,T]$, for a small positive real $\delta_0$ depending on $L$, $T$ and the Lipschitz norm of $\phi$ only. By Theorem 1.3 in \cite{Delarue02}, the functions $(v^{\varepsilon})_{0 < \varepsilon <1}$ converge towards $v^{\phi}$ as $\varepsilon$ tends $0$, uniformly on compact subsets of $[T-\delta,T] \times \RR^d \times \RR$.

Following Section 2 in \cite{Delarue02}, we can approximate Eq. \eqref{eq:4:11:5} by equations of the same type, but driven by bounded smooth coefficients. In particular, we can apply Propositions \ref{prop:4:11:1} and \ref{prop:20:07:2} to the approximated equations: passing to the limit along the approximation, this shows that $v^{\varepsilon}(T-\delta_0,\cdot,\cdot)$ is $C$-Lipschitz continuous w.r.t. the variables $p$ and $e$, for some constant $C$ depending on $\delta_0$, $L$ and $T$ only. Letting $\varepsilon$ tend to $0$, we deduce that $v^{\phi}(T-\delta_0,\cdot,\cdot)$ is $C$-Lipschitz as well. We then follow the induction scheme in \cite{Delarue02}: we can solve Eqs. \eqref{eq:4:11:5} and \eqref{eq:4:11:1} on some interval $[T- (\delta_0 + \delta_1),T-\delta_0]$, with $v^{\varepsilon}(T-\delta_0,\cdot,\cdot)$ and $v^{\phi}(T-\delta_0,\cdot,\cdot)$ as respective terminal conditions. Here, $\delta_1$ depends on $C$, $L$ and $T$ only, i.e. on $L$, $T$ and $\delta_0$ only. This permits to extend the value functions  $v^{\varepsilon}$ and $v^{\phi}$ to 
$[T-(\delta_0 + \delta_1),T]$. Iterating the process, we can define the value functions on the whole $[0,T] \times \RR^d \times \RR$: they are Lipschitz continuous w.r.t. the variables $p$ and $e$; and the functions $(v^{\varepsilon})_{0 <\varepsilon < 1}$ converge towards $v^{\phi}$ as $\varepsilon$ tends to $0$, uniformly on compact subsets of
$[0,T] \times \RR^d \times \RR$. 
The value functions being now constructed, existence and uniqueness are proven as in Theorem 2.6 in \cite{Delarue02}. 

It then remains to check that $v^{\varepsilon}$ is a
${\mathcal C}^{1,2}$ solution 
of the PDE \eqref{eq:7:11:10} on the whole $[0,T] \times \RR^d \times \RR$. Basically, this follows from Schauder's estimates, see Chapter 8 in Krylov \cite{Krylov:book}. Indeed, when the coefficients are bounded, PDE \eqref{eq:7:11:10} may be seen as a uniformly elliptic PDE with locally H\"older continuous coefficients, the local H\"older norms of the coefficients being independent of the $L^{\infty}$-bounds of the coefficients: this is a consequence of the last assertion in Proposition \ref{prop:20:07:2}. By Schauder's estimates, the H\"older norms of the first-order derivatives in time and space and second-order derivatives in space are bounded on compact subsets of $[0,T] \times \RR^d \times \RR$, independently of the $L^{\infty}$-bounds of the coefficients. 
Approximating the coefficients in \eqref{eq:4:11:5} by bounded coefficients as in the previous paragraph and passing to the limit along the approximation, we deduce that $v^{\varepsilon}$ is a classical solution of the PDE 
\eqref{eq:7:11:10}. Similarly, we deduce that \eqref{eq:9:11:1} and \eqref{eq:20:07:3} hold true at the limit as well.
\qed

\subsubsection{Existence for a Singular Terminal Condition}
To pass to the limit along the mollification of the terminal condition, we need first to understand the boundary behavior of the solution of the mollified equation \eqref{eq:4:11:5}. We thus claim (the results provided in the following proposition will be refined in the next section when we restrict ourselves to binary terminal conditions $\phi$):

\begin{proposition}
\label{prop:03:11:1}
Consider the mollified equation \reff{eq:4:11:5} with a non-decreasing Lipschitz smooth terminal condition $\phi$ satisfying \eqref{eq:25:7:1}. Then, for any $\rho >0$ and $q \geq 1$, there exists a constant $C(\rho,q) >0$, only depending on $\rho$, $L$ and $T$, such that for any $t \in [0,T)$, $e, \ \Lambda \in\RR$, and $|p| \leq \rho$, we have:
\begin{equation}
\label{eq:7:11:20}
\begin{split}
&e >  \Lambda  \Rightarrow
v^{\varepsilon}(t,p,e) \geq 
\phi(\Lambda) -  C(\rho,q) \bigl( \frac{e-\Lambda}{L(T-t)} \bigr)^{-q},
\\
&e < \Lambda  \Rightarrow
v^{\varepsilon}(t,p,e) \leq \phi(\Lambda)+ C(\rho,q) \bigl( \frac{\Lambda-e}{L(T-t)} \bigr)^{-q}.
\end{split}
\end{equation}
In particular, for any $t<T$ and $p \in \RR^d$
\begin{equation}
\label{eq:26:7:11}
\lim_{e \rightarrow + \infty} v^{\phi}(t,p,e) = 1, \quad
\lim_{e \rightarrow - \infty} v^{\phi}(t,p,e) = 0,
\end{equation}
$v^{\phi}$ being given by Corollary \ref{cor:25:7:1}.
\end{proposition}

{\bf Proof.} In \eqref{eq:7:11:20}, we will prove the lower bound only as the upper bound can be proven in a similar fashion. For a given starting point $(t_0,p,e) \in [0,T) \times \RR^d \times \RR$, 
we consider the solution $(P_t^{\varepsilon,t_0,p},E_t^{\varepsilon,t_0,p,e},Y_t^{\varepsilon,t_0,p,e})_{t_0 \leq t \leq T}$ of the system \reff{eq:4:11:5}. We ignore the superscrit $(t_0,p,e)$ for the sake of convenience. We have:
\begin{equation*}
v^{\varepsilon}(t_0,p,e) =\EE[\phi(E^{\varepsilon}_T)]
\geq \phi(\Lambda) - \PP \{E^{\varepsilon}_T < \Lambda\}.
\end{equation*}
Since $f$ is increasing w.r.t. $y$, for $e>\Lambda$,
\begin{equation*}
\begin{split}
\PP \{E^{\varepsilon}_T < \Lambda\}
&= \PP \biggl\{ \int_{t_0}^T f(P_s^{\varepsilon},1) ds + \varepsilon (B_T - B_{t_0}) > e-\Lambda \biggr\}
\\
&\leq {\mathbb P} \biggl\{  L \bigl( 1 + \sup_{t_0 \leq s \leq T} |P_s^{\varepsilon}| \bigr) 
+ \varepsilon (B_T - B_{t_0}) > \frac{e-\Lambda}{T-t_0} \biggr\} \leq C(p,q)\left(\frac{e-\Lambda}{L(T-t_0)}\right)^{-q}.
\end{split}
\end{equation*}
Indeed, under the linear growth assumption (A.2),
$\sup_{t_0 \leq s \leq T} |P_s^{\varepsilon}|$ has finite $q$-moments, for any $q \geq 1$. (Keep in mind that $\varepsilon \in (0,1)$, see \eqref{eq:4:11:5}.)
The proof of \eqref{eq:7:11:20} is easily completed. 

Eq. \eqref{eq:26:7:11} easily follows. By Corollary \ref{cor:25:7:1}, we can let $\varepsilon$ tend to $0$ in \eqref{eq:7:11:20}. For $e>0$, we obtain 
$v^{\phi}(t,p,e) \geq \phi(e/2) - C(\rho,1)[e/(2LT)]^{-1}$. By \eqref{eq:25:7:1}, we deduce that $\lim_{e \rightarrow + \infty} v^{\phi}(t,p,e) = 1$. The limit in $-\infty$ is proven in the same way.
\qed

\vskip 6pt

Below, we turn to the case when 
the function $\phi$ giving the terminal condition $Y_T=\phi(E_T)$ in \eqref{eq:4:11:1}
is possibly discontinuous. As announced, we go back to the smooth setting by mollification:

\begin{example}
\label{ex:mollifying}
Consider a non-decreasing function $\phi$ as in \eqref{eq:25:7:1} and $\phi_-$ and $\phi_+$ as in \eqref{fo:phi-+}. Notice that $\phi_+$ is a cumulative distribution function as a non-decreasing right-continuous function matching $0$ at $-\infty$ and $1$ at $+ \infty$. Notice also that $\phi_-$ is the left-continuous version of $\phi_+$.

We now construct two specific mollifying sequences for $\phi$.
Let $j$ be the density of a positive random variable, $j$ 
being $C^\infty$ with a compact support. Let $\xi$ and $\vartheta$ be independent random variables, $\xi$ with $\phi_+$ as cumulative distribution function and $\vartheta$ with $j$ as density. For each integer $n \geq 1$, denote by $\phi_+^n$ and $\phi_-^n$ the cumulative distribution functions of the random variables $\xi- n^{-1}\vartheta$ and $\xi+n^{-1}\vartheta$ respectively. Then, the functions $\phi_+^n$ and $\phi_-^n$ are non-decreasing with values in $[0,1]$. They are $C^\infty$, with bounded derivatives of any order. Moreover, 
$\phi_-^n\le \phi$ and $\phi_+^n\ge \phi$ and the sequences
 $(\phi_+^n)_{n \geq 1}$ and  $(\phi_-^n)_{n \geq 1}$ converge pointwise towards $\phi_+$ and $\phi_-$ respectively as $n$ tends $+\infty$. Finally,
\begin{equation*}
\int_{\RR} |\phi_+^n(e) - \phi_+(e)| de \leq \int_{\RR \times \RR^+} {\mathbb P}\{ e  \leq \xi \leq e + t/n\} j(t) dt de = \frac{1}{n} \int_{\RR^+} t j(t) dt\rightarrow 0
\end{equation*}
as $n$ tends to $+ \infty$,
so that the convergence of
 $(\phi_+^n)_{n \geq 1}$ towards $\phi_+$
holds in $L^1(\RR)$ as well. (Obviously, the same holds for
$(\phi_-^n)_{n \geq 1}$ and $\phi_-$.)
\end{example}

Existence of a solution in Theorem \ref{thm:4:11:1} follows from
\begin{proposition}
\label{prop:7:11:10}
There exists a continuous function $v : [0,T) \times \RR^d \times \RR \rightarrow [0,1]$ satisfying
\begin{enumerate}
\item $v(t,\,\cdot\,,\,\cdot\,)$ is $1/[\ell_1(T-t)]$-Lipschitz continuous with respect to $e$ for any $t \in [0,T)$,
\item  $v(t,\,\cdot\,,\,\cdot\,)$ is $C$-Lipschitz continuous with respect to $p$ for any 
$t \in [0,T)$, $C$ as in \eqref{eq:20:07:3},
\item for each $e\in\RR$, $\phi_-(e)\le \lim_{t \rightarrow T} v(t,p,e) \le \phi_+(e)$ uniformly in $p$ in compact subsets of $\RR^d$.
\end{enumerate}
Moreover, for any $(t_0,p,e) \in [0,T) \times \RR^d \times \RR$, the strong solution
$(E_t^{t_0,p,e})_{t_0 \leq t <T}$ of 
\begin{equation}
\label{eq:25:7:10}
E_t = e - \int_{t_0}^t f\bigl(P_s^{t_0,p},v(s,P_s^{t_0,p},E_s)\bigr) ds, \quad t_0 \leq t <T,
\end{equation}
is such that $(v(t,P_t^{t_0,p},E_t^{t_0,p,e}))_{t_0 \leq t <T}$ is a $[0,1]$-valued martingale with respect to
the complete filtration generated by $W$, the integrand in the martingale representation of 
$(v(t,P_t^{t_0,p},E_t^{t_0,p,e}))_{t_0 \leq t <T}$ being bounded by a constant depending on $L$ and $T$ only. Moreover, $\PP$-almost surely,
\begin{equation}
\label{fo:squeeze}
\phi_-(E_T^{t_0,p,e}) \leq \lim_{t \nearrow T}
v(t,P_t^{t_0,p},E_t^{t_0,p,e}) \leq \phi_+(E_T^{t_0,p,e}).
\end{equation}
\end{proposition}
\noindent
Here $(P_t^{t_0,p})_{t_0 \leq t \leq T}$ is the solution of the forward equation in \reff{eq:4:11:1} starting from $p$ at time $t_0$. We emphasize that the limit $\lim_{t \nearrow T}
v(t,P_t^{t_0,p},E_t^{t_0,p,e})$ exists as the a.s. limit of a non-negative martingale. 

\vskip 6pt
\textbf{Proof.}
Let $(\phi^n)_{n\ge 1}$ be a mollified approximation of $\phi_+$ as constructed 
in Example \ref{ex:mollifying}. By Corollary \ref{cor:25:7:1}, we can consider the value functions $(v^n = v^{\phi^n})_{n \geq 1}$ associated with the terminal conditions 
$(\phi^n)_{n \geq 1}$. Following the proof of Corollary \ref{cor:25:7:1}, we deduce from Propositions \ref{prop:4:11:1} and \ref{prop:20:07:2} that, for any $T'<T$ and any compact $K \subset \RR^d$, the functions $(v^{n})_{n \geq 1}$ are equicontinuous on $[0,T']\times K\times \RR$. Let us denote by $v$ a function constructed on $[0,T)\times\RR^d\times \RR$ as the limit
of a subsequence $(v^{\varphi(n)})_{n \geq 1}$ of $(v^{n})_{n \geq 1}$ that converges locally uniformly. 
By \eqref{eq:9:11:1} and \eqref{eq:20:07:3}, $v$ is Lipschitz continuous w.r.t. 
$(p,e)$, uniformly in $t$ in compact subsets of $[0,T)$: thus, 
\eqref{eq:25:7:10} has a strong unique solution on $[t_0,T)$ for any initial condition
$(t_0,p,e)$. By Cauchy criterion, the limit of $E_t^{t_0,p,e}$ exists a.s. as $t \nearrow T$: it is denoted by $E_T^{t_0,p,e}$. 

For any initial condition $(t_0,p,e) \in [0,T) \times \RR^d \times \RR$ and any integer $n \geq 1$, we also denote by $(P_t^{t_0,p},E_t^{n,t_0,p,e},Y_t^{n,t_0,p,e},Z_t^{n,t_0,p,e})_{t_0 \leq t \leq T}$ the solution to \eqref{eq:4:11:1} with $\phi^n$ as initial condition. Since
$(v^{\varphi(n)})_{n \geq 1}$ converges towards $v$ uniformly on compact subsets of $[0,T)
\times \RR^d \times \RR$, we deduce that, a.s., $E_t^{\varphi(n),t_0,p,e} \rightarrow
E_t^{t_0,p,e}$ uniformly in time $t$ in compact subsets of $[t_0,T)$. Since the process
$(v^{\varphi(n)}(t,P_t^{t_0,p},E_t^{\varphi(n),t_0,p,e})=Y_t^{\varphi(n),t_0,p,e})_{t_0 \leq t <T}$ is a $[0,1]$-valued martingale
for any $n \geq 1$, $(v(t,P_t^{t_0,p},E_t^{t_0,p,e}))_{t_0 \leq t <T}$ is also a  martingale as the a.s. limit of bounded martingales. The bound for the integrand in the martingale representation follows from \eqref{eq:20:07:3}.

Applying \reff{eq:7:11:20} to each $\phi^{\varphi(n)}$, $n \geq 1$, and letting $\varepsilon$ tend to 0 therein, we understand that each $v^{\varphi(n)}$, $n \geq 1$, satisfies 
\eqref{eq:7:11:20} as well. Letting $n$ tend to $+\infty$, we deduce that 
$v$ also satisfies \eqref{eq:7:11:20}. As a consequence, for any family $(p_t,e_t)_{0 \leq t \leq T}$ converging towards $(p,e)$ as $t \nearrow T$, we have:
$$
\phi_-(e)\le\liminf_{t\nearrow T}v(t,p_t,e_t) \le\limsup_{t\nearrow T}v(t,p_t,e_t)\le \phi_+(e),
$$
which implies \reff{fo:squeeze}. \qed

\subsection{Uniqueness via Conservation Law}

\subsubsection{Key Lemmas}
As in \cite{cdet}, the proof of uniqueness relies on two key ingredients: a comparison lemma for the solutions to \eqref{eq:4:11:1}, see Lemma \ref{lem:7:11:2} below, and a conservation lemma for the value functions of the mollified equation
\eqref{eq:4:11:5}, see Lemma \ref{lem:7:11:1}. 

\begin{lemma}
\label{lem:7:11:2}
Let $\phi'$ be another non-decreasing terminal condition satisfying condition \eqref{eq:25:7:1} and let $(P_t,E_t',Y_t')_{t_0 \leq t \leq T}$, $t_0 \in [0,T)$, be a solution of equation \reff{eq:4:11:1} satisfying
\begin{equation*}
\PP\{\phi_-'(E_T') \leq Y_T' \leq \phi_+'(E_T') \}=1,
\end{equation*}
in lieu of terminal condition. 

Let $(w^{\varepsilon,n})_{0 < \varepsilon <1,n \geq 1}$ be a family of classical solutions of \eqref{eq:7:11:10} associated with a non-increasing sequence of mollified non-decreasing 
$[0,1]$-valued terminal conditions $(\chi^{n})_{n \geq 1}$ satisfying
$\chi^n \searrow \phi_+$. 
(In particular, $\chi^n \geq \phi_+$.)
If $\phi \geq \phi'$, then
\begin{equation}
\label{eq:26:7:10}
\PP \bigl\{
\liminf_{n \rightarrow + \infty} \liminf_{\varepsilon \rightarrow 0} w^{\varepsilon,n}(t,P_t,E_t') \geq Y_t'
\bigr\}=1, \quad t \in [t_0,T).
\end{equation}
Similarly, if $(\chi^{n})_{n \geq 1}$ is a non-decreasing sequence of mollified non-decreasing $[0,1]$-valued terminal conditions satisfying 
$\chi^n \nearrow \phi_-$ and $\phi \leq \phi'$, then 
$\limsup_{n \rightarrow + \infty} \limsup_{\varepsilon \rightarrow 0} w^{\varepsilon,n}(t,P_t,E_t') \leq Y_t'$, $t \in [t_0,T)$.
\end{lemma}

{\bf Proof.} We prove the first part only, the proof of the second part being similar. We apply 
It\^o's formula to $(Y_t^{\varepsilon,n} =w^{\varepsilon,n}(t,P_t^{\varepsilon},E_t'+\varepsilon B_t))_{t_0 \leq t \leq T}$
where $dP_t^{\varepsilon}=b(P_t^{\varepsilon})dt + \sigma(P_t^{\varepsilon})dW_t + \varepsilon dW_t'$, with $P_{t_0}^{\varepsilon}=P_{t_0}$. Using the PDE \reff{eq:7:11:10} satisfied by $w^{\varepsilon,n}$, we obtain
\begin{equation}
\label{eq:22:12:1}
\begin{split}
dY_t^{\varepsilon,n} 
&= \bigl[ 
f\bigl(P_t^{\varepsilon},w^{\varepsilon,n}(t,P_t^{\varepsilon},E_t'+ \varepsilon B_t) \bigr)
-f\bigl(P_t,Y_t'\bigr)
\bigr] \partial_e w^{\varepsilon,n}(t,P_t^{\varepsilon},E_t'+\varepsilon B_t)
dt
\\
&\hspace{15pt}
+ \langle \partial_p w^{\varepsilon,n}(t,P_t^{\varepsilon},E_t'+\varepsilon B_t), \sigma(P_t^{\varepsilon}) dW_t + \varepsilon dW_t' \rangle
+
\varepsilon \partial_e w^{\varepsilon,n}(t,P_t^{\varepsilon},E_t'+\varepsilon B_t)  dB_t
\\
&= \bigl[ 
f\bigl(P_t^{\varepsilon},Y_t^{\varepsilon,n} \bigr)
- f\bigl(P_t,Y_t'\bigr)
\bigr]\partial_e w^{\varepsilon,n}(t,P_t^{\varepsilon},E_t'+\varepsilon B_t)
dt + dM_t^{\varepsilon,n},
\end{split}
\end{equation}
where $(M_t^{\varepsilon,n})_{t_0 \leq t \leq T}$ is a square-integrable martingale with respect to the filtration generated by $(W_t)_{t_0 \leq t \leq T}$, $(W_t')_{t_0 \leq t \leq T}$ and $(B_t)_{t_0 \leq t \leq T}$. Modifying $(M_t^{\varepsilon,n})_{t_0 \leq t \leq T}$, we see that:
\begin{equation*}
d\bigl[Y_t^{\varepsilon,n} - Y_t' \bigr] = \bigl[ 
f\bigl(P_t^{\varepsilon},Y_t^{\varepsilon,n} \bigr)
- f\bigl(P_t,Y_t'\bigr)
\bigr] \partial_e w^{\varepsilon,n}(t,P_t^{\varepsilon},E_t'+\varepsilon B_t)
dt + dM_t^{\varepsilon,n}.
\end{equation*}
We now apply It\^o's formula with the function $y \mapsto (y^-)^2=(\min(y,0))^2$. Recall that $f$ is increasing in its second argument, that $\partial_e w^{\varepsilon,n}$ is non-negative and that both $Y^{\varepsilon,n}$ and $Y'$ are $[0,1]$-valued. We obtain
\begin{equation*}
d \bigl[(Y_t^{\varepsilon,n} - Y_t')^-]^2 \geq - 2  L \sup_{\varepsilon' \in (0,1)} \|\partial_e w^{\varepsilon',n}\|_{\infty} |P_t^{\varepsilon}-P_t| dt + dM_t^{\varepsilon,n},
\quad t_0 \leq t \leq T,
\end{equation*}
for a new choice of  $(M_t^{\varepsilon,n})_{t_0 \leq t \leq T}$. Letting $\varepsilon$ tend to $0$ and applying the second part in Proposition \ref{prop:4:11:1}, we deduce that, for any $t \in [t_0,T]$, a.s., 
\begin{equation*}
\limsup_{\varepsilon \rightarrow 0} \bigl[(Y_t^{\varepsilon,n} - Y_t')^- \bigr]^2 
\leq \limsup_{\varepsilon \rightarrow 0}
{\mathbb E} \bigl\{ \bigl[(Y_T^{\varepsilon,n} - Y_T')^- \bigr]^2 | {\mathcal F}_t \bigr\} 
= {\mathbb E} \bigl\{ \bigl[(\chi^n(E_T') - Y_T')^- \bigr]^2 | {\mathcal F}_t \bigr\} = 0.
\end{equation*}
The result easily follows. \qed

\begin{remark}
\label{rem:26:7:1}
When $\phi$ and $\phi'$ are Lipschitz smooth functions, Lemma \ref{lem:7:11:2} reads as a comparison principle for the value functions $v^{\phi}$ and $v^{\phi'}$ given by Corollary \ref{cor:25:7:1}. Indeed, choosing $t=t_0$ in \eqref{eq:26:7:10}, we obtain $\liminf_{n \rightarrow + \infty} \lim_{\varepsilon \rightarrow 0}
w^{\varepsilon,n}(t_0,p,e) \geq v^{\phi'}(t_0,p,e)$ provided $\phi \geq \phi'$. By Corollary \ref{cor:25:7:1}, we deduce  $\lim_{n \rightarrow + \infty} v^{\chi^n}(t_0,p,e) \geq v^{\phi'}(t_0,p,e)$. Following the proof of Proposition \ref{prop:7:11:10}, we know that the sequence $(v^{\chi^n})_{n \geq 1}$ is uniformly continuous on compact subset of $[0,T] \times \RR^d \times \RR$ and that any possible limit $v$ provides a solution to \eqref{eq:4:11:1} with $\phi$ as terminal condition. By the uniqueness part in Corollary \ref{cor:25:7:1}, $v$ always matches $v^{\phi}$, so that $v^{\phi}(t_0,p,e) \geq v^{\phi'}(t_0,p,e)$ when $\phi \geq \phi'$. 

In particular, when $(\chi^n)_{n \geq 1}$ is a non-increasing (resp. non-decreasing) sequence of mollified terminal conditions as in the statement of Lemma \ref{lem:7:11:2}, the sequence of functions $(v^{\chi^n})_{n \geq 1}$ is non-increasing (resp. non-decreasing) as well. In particular, it is convergent. As already explained, when $\phi$ is a Lipschitz smooth function, the limit is $v^{\phi}$ exactly.
When $\phi$ is possibly discontinuous, \eqref{eq:26:7:10} yields: $\lim_{n \rightarrow + \infty}
v^{\chi^n}(t,P_t,E_t') \geq Y_t'$, $t \in [t_0,T)$. 
\end{remark}

\begin{lemma}
\label{lem:7:11:1}
Let us consider two sequences of mollified terminal conditions
$(\chi^{n}_i)_{n \geq 1}$, $i=1,2$, satisfying \eqref{eq:25:7:1}
such that $(\chi^n_i-\phi)_{n \geq 1}$ converges towards $0$
in $L^1(\RR)$ for $i=1,2$.
Let us also consider the associated families $(w^{\varepsilon,n}_i)_{\varepsilon \in (0,1),n \geq 1}$, $i=1,2$, of solutions to the PDE \reff{eq:7:11:10}, $\varepsilon$ standing for the viscosity parameter. Then, for any $(t,p) \in [0,T) \times \RR$,
\begin{equation*}
\lim_{n \rightarrow + \infty} 
\lim_{m \rightarrow +\infty} 
\lim_{\varepsilon \rightarrow 0} 
\int_{-m}^m \bigl[  w^{\varepsilon,n}_2(t,p,e)
- w^{\varepsilon,n}_1(t,p,e) \bigr] de  = 0.
\end{equation*}
\end{lemma}

{\bf Proof.} For an integer $m \geq 1$, we set
$\displaystyle W_i^{\varepsilon,m,n}(t,p) = \int_{-m}^m w_i^{\varepsilon,n}(t,p,e) de$.

By integration of the PDE \reff{eq:7:11:10},
it satisfies the PDE
\begin{equation*}
\begin{split}
&\partial_t W_i^{\varepsilon,m,n}(t,p)
+ \bigl( {\mathcal L}_p + \frac{\varepsilon^2}{2} \partial^2_{pp}
\bigr) \bigl[  W_i^{\varepsilon,m,n} \bigr](t,p)
+ \frac{\varepsilon^2}{2} \bigl[ \partial_e w_i^{\varepsilon,n}(t,p,m) -   \partial_e w_i^{\varepsilon,n}(t,p,-m) \bigr]
\\
&\hspace{15pt} - F\bigl(p,w_i^{\varepsilon,n}(t,p,m)\bigr)
+  F\bigl(p,w_i^{\varepsilon,n}(t,p,-m)\bigr) = 0,
\end{split}
\end{equation*}
where $\displaystyle F(p,v) = \int_0^v f(p,r) dr$.

In particular, we can give a probabilistic representation of $W_i(t,p,e)$ in terms of the process
$(P^{\varepsilon}_s)_{t \leq s \leq T}$, solution of $dP_s^{\varepsilon} = b(P_s^{\varepsilon})ds + 
\sigma(P_s^{\varepsilon})dW_s + \varepsilon dW_s'$, $P_{t}^{\varepsilon}=p$:
\begin{equation*}
\begin{split}
W_i^{\varepsilon,m,n}(t,p) &= {\mathbb E} \bigl[ W_i^{\varepsilon,m,n}(T,P_T^{\varepsilon}) \bigr]
\\
&\hspace{5pt} + {\mathbb E} \int_t^T 
\biggl[
\frac{\varepsilon^2}{2} \bigl[ \partial_e w_i^{\varepsilon,n}(s,P_s^{\varepsilon},m) -   \partial_e w_i^{\varepsilon,n}(s,P_s^{\varepsilon},-m) \bigr]
\\
&\hspace{50pt} - \bigl[ F\bigl(P_s^{\varepsilon},w_i^{\varepsilon,n}(s,P_s^{\varepsilon},m)\bigr)
-  F\bigl(P_s^{\varepsilon},w_i^{\varepsilon,n}(s,P_s^{\varepsilon},-m)\bigr) \bigr] \biggr] ds
\\
&= T_i^{\varepsilon,m,n}(1) + T_i^{\varepsilon,m,n}(2).
\end{split}
\end{equation*}
Let us examine how $T_i^{\varepsilon,m,n}(2)$ behaves as $\varepsilon$ tends to $0$, keeping in mind that, for every mollifying parameter $n \geq 1$,
$\partial_e w_i^{\varepsilon,n}$ is bounded on the whole domain independently of $\varepsilon$ (see Proposition 
\ref{prop:4:11:1})
and $(w_i^{\varepsilon,n})_{0 < \varepsilon <1}$ converges towards $v^{\chi^n_i}$, uniformly on compact sets, as $\varepsilon$ tends to $0$ (see Corollary \ref{cor:25:7:1}).
By uniform integrability of the family $(\sup_{t_0 \leq t \leq T} |P_t^{\varepsilon}|)_{0 < \varepsilon < 1}$,
we deduce that 
\begin{equation*}
\lim_{\varepsilon \rightarrow 0} T_i^{\varepsilon,m,n}(2)
= {\mathbb E} \int_t^T  \bigl[ F\bigl(P_s,v^{\chi^n_i}(s,P_s,-m)
\bigr)
-  F\bigl(P_s,v^{\chi^n_i}(s,P_s,m)\bigr) \bigr]  ds.
\end{equation*}
By \eqref{eq:26:7:11} and by dominated convergence, 
\begin{equation*}
\lim_{m \rightarrow + \infty} \lim_{\varepsilon \rightarrow 0} T_i^{\varepsilon,m,n}(2)
= {\mathbb E} \int_t^T  \bigl[ F(P_s,0)
-  F(P_s,1) \bigr]  ds.
\end{equation*}
In particular, the limit of the difference
$T_2^{\varepsilon,m,n}(2)-T_1^{\varepsilon,m,n}(2)$ is zero (as $\varepsilon$ tends to 0 first and then as $m$ tends to $+ \infty$). 
Let us now consider $T_i^{\varepsilon,m,n}(1)$. Clearly, 
$W_i^{\varepsilon,m,n}(T,p) = \int_{-m}^m \chi_i^n(e) de$
is independent of $p$, so that the limit of the difference 
$T_2^{\varepsilon,m,n}(1)-T_1^{\varepsilon,m,n}(1)$
is 
\begin{equation*}
\lim_{m \rightarrow + \infty} \lim_{\varepsilon \rightarrow 0} \bigl[ T_2^{\varepsilon,m,n}(1) - T_1^{\varepsilon,m,n}(1) \bigr] = \int_{\RR} \bigl( \chi_2^n- \chi_1^n \bigr)(e) de.
\end{equation*}
Note that $\chi_2^n - \chi_1^n = (\chi_2^n - \phi) - (\chi_1^n - \phi)$
converges towards 0 in $L^1(\RR)$ as $n$ tends to $+\infty$. 
\qed

\subsubsection{End of the Proof of Uniqueness in Theorem \ref{thm:4:11:1}.}
By Example \ref{ex:mollifying}, we know that 
$(\phi^n_- - \phi_-)_{n \geq 1}$ and $(\phi^n_+ - \phi_+)_{n \geq 1}$ converge towards
$0$ in $L^1(\RR)$; since $\phi=\phi_+=\phi_-$ almost everywhere for the Lebesgue
measure on $\RR$, $(\phi^n_- - \phi)_{n \geq 1}$ and $(\phi^n_+ - \phi)_{n \geq 1}$ converge towards $0$ in $L^1(\RR)$. 
Therefore, we can apply Lemma \ref{lem:7:11:1} with $(\chi^n_1 =\phi^n_-)_{n \geq 1}$ and $(\chi_2^n = \phi^n_+)_{n \geq 1}$. Denoting by 
$(w_i^{\varepsilon,n})_{\varepsilon \in (0,1),n \geq 1}$ the associated 
 solutions with $(\chi_i^n = \phi^n_+)_{n \geq 1}$, $i=1,2$, we obtain
\begin{equation*}
\lim_{n \rightarrow + \infty} 
 \lim_{m \rightarrow +\infty} 
\lim_{\varepsilon \rightarrow 0}
\int_{-m}^m \bigl[  w^{\varepsilon,n}_2(t,p,e)
- w^{\varepsilon,n}_1(t,p,e) \bigr] de  = 0.
\end{equation*}
By Corollary \ref{cor:25:7:1},
\begin{equation*}
\lim_{n \rightarrow + \infty} 
 \lim_{m \rightarrow +\infty} 
\int_{-m}^m \bigl[  v^{\phi_+^n}(t,p,e)
- v^{\phi^n_-}(t,p,e) \bigr] de  = 0.
\end{equation*}
By Remark \ref{rem:26:7:1}, we have $v^{\phi_+^n} \geq  v^{\phi^n_-}$, so that  
the integrand above is non-negative, that is, for any $m \geq 1$,
\begin{equation*}
\lim_{n \rightarrow + \infty} 
\int_{-m}^m \bigl|  v^{\phi_+^n}(t,p,e)
- v^{\phi^n_-}(t,p,e) \bigr| de  = 0.
\end{equation*}
By Remark \ref{rem:26:7:1}, we know that
the sequences
$(v^{\phi^n_-})_{n \geq 1}$ and $(v^{\phi^n_+})_{n \geq 1}$ are pointwise convergent on $[0,T) \times \RR^d \times \RR$. 
By Propositions \ref{prop:4:11:1} and \ref{prop:20:07:2}, convergence is uniform on compact subsets
of $[0,T) \times \RR^d \times \RR$. Therefore,
\begin{equation*}
\lim_{n \rightarrow + \infty}  v^{\phi_+^n}(t,p,e)
= \lim_{n \rightarrow + \infty} v^{\phi^n_-}(t,p,e), \quad (t,p,e) \in [0,T) \times \RR^d \times \RR.
\end{equation*}
By construction, the limit matches $v(t,p,e)$ in Proposition \ref{prop:7:11:10}. 

Now, for any solution $(P_t,E'_t,Y_t')_{t_0 \leq t \leq T}$ of equation 
\reff{eq:4:11:1}, we can apply Lemma \ref{lem:7:11:2} with $(\chi^n=\phi^{\pm}_n)_{n \geq 1}$. Passing to the limit in \eqref{eq:26:7:10}, we obtain $v(t,P_t,E_t') \leq Y_t' \leq v(t,P_t,E_t')$, $t_0 \leq t < T$,
so that $(E_t')_{t_0 \leq t < T}$ satisfies \eqref{eq:25:7:10}. Uniqueness easily follows. \qed

\vskip 6pt\noindent
The following corollary will be useful in the sequel:

\begin{corollary}
\label{cor:27:7:1}
Consider a sequence of non-decreasing Lipschitz smooth terminal conditions $(\phi^n)_{n \geq 1}$ converging pointiwse towards $\phi$ as $n$ tend to $+\infty$. 
Then, the functions $(v^{\phi^n})_{n \geq 1}$ converge towards 
$v$ in Proposition \ref{prop:7:11:10} as $n$ tends to $+\infty$, uniformly on compact subsets of $[0,T) \times \RR^{d} \times \RR$.
\end{corollary}

{\bf Proof.} By Propositions \ref{prop:4:11:1} and \ref{prop:20:07:2}, the functions $(v^{\phi^n})_{n \geq 1}$ are uniformly continuous on 
every compact subset of $[0,T) \times \RR^d \times \RR$. Passing to the limit in \eqref{eq:7:11:20}, the limit 
$w$ of any converging subsequence of $(v^{\phi^n})_{n \geq 1}$ satisfies Proposition 
\ref{prop:7:11:10} with the right relaxed terminal condition, that is $w=v$ by uniqueness. \qed

\section{Dirac Mass at the Terminal Time}
\label{se:dirac}
\noindent
From now on we restrict ourselves to the case $\phi={\mathbf 1}_{[\Lambda,+\infty)}$, $\Lambda \in \RR$, and we assume:

(A.3) For any $p \in \RR^d$, the function $y \hookrightarrow f(p,y)$ is differentiable with respect to $y$  and there exists  $\alpha \in (0,1]$ such that, for any $(p,p',y,y')
\in \RR^d \times \RR^d \times \RR \times \RR$, we have
$$
|\partial_y f(p,y) - \partial_y f(p',y')| \leq L \bigl( |p'-p|^{\alpha} + |y-y'|^{\alpha} \bigr).
$$

(A.4) The drift $b$ and the matrix $\sigma$ are bounded by $L$.
\vskip 6pt

The main result of this section may be summarized as follows: 
\emph{(i)} Under (A.3) and (A.4), there is a cone of initial conditions $(t_0,p,e)$ for which the distribution of the random variable $E_T^{t_0,p,e}$ has a Dirac mass
at the singular point $\Lambda$. Put differently, there is a non-zero event of scenarii for which the terminal conditions $\phi_-(E_T^{t_0,p,e})$ and $\phi_+(E_T^{t_0,p,e})$ in the terminal condition \eqref{eq:4:11:2} differ: this makes the relaxation of the terminal condition meaningful. The complete statement is given in Proposition \ref{prop:30:12:1}. 
\emph{(ii)} When the diffusion matrix $\RR^d \ni p \hookrightarrow [\sigma \sigma^{\top}](p)$ is uniformly elliptic and the gradient
$\RR^d \ni p \hookrightarrow \partial_p f(p,0)$ is uniformly continuous and uniformly away from zero, the Dirac mass exists for any initial condition $(t_0,p,e) \in [0,T) \times \RR^d \times \RR$. Moreover, the topological support of the conditional law of $Y_T^{t_0,p,e}$ given $E_T^{t_0,p,e} = \Lambda$ is the entire
$[0,1]$, that is, conditionally on the non-zero event $E_T^{t_0,p,e} = \Lambda$, all the values between $\phi_-(\Lambda)=0$ and $\phi_+(\Lambda)=1$ may be observed in the relaxed terminal condition \eqref{eq:4:11:2}. In particular, the $\sigma$-algebra $\sigma(Y_T^{t_0,p,e})$
is not included into the $\sigma$-algebra $\sigma(E_T^{t_0,p,e})$: because of the degeneracy of the forward equation and of the singularity of the terminal condition, the standard Markovian structure breaks down at terminal time. We refer to Proposition \ref{prop:4:1:1} for the complete statement.

The strategy of the proof consists in a careful analysis of the trajectories of the process $(E_t)_{0 \leq t \leq T}$. Precisely, we compare the trajectories of 
$(E_t)_{0 \leq t \leq T}$ with the characteristics of the non-viscous version of PDE \eqref{eq:7:11:10}, i.e. of the first-order PDE $\partial_t u(t,p,e) -f(p,u(t,p,e))\partial_e u(t,p,e) = 0$ with $u(T,p,e)=\phi(e)$ as boundary condition. Because of the singularity of $\phi$, the characteristics of the PDE merge at $\Lambda$ at time $T$. This phenomenon is called a \emph{shock} in the PDE literature. Here the shock acts as a trap for the trajectories of $(E_t)_{0 \leq t \leq T}$ enclosing a non-zero mass of the process into a cone narrowing towards $\Lambda$: because of the degeneracy of the forward process in \eqref{eq:4:11:1}, the noise may not be sufficient enough to let the process $(E_t)_{0 \leq t \leq T}$
escape from the trap. Here is the collateral effect of the simultaneity of the singularity of $\phi$ and the degeneracy of $(E_t)_{0 \leq t \leq T}$.

\subsection{Change of Variable}
For the sake of convenience, we switch from the degenerate component $E$ of the forward process to a process $\bar{E}$ which has the same terminal value, hence leaving the terminal condition of the backward process unchanged, and which will be easier to manipulate. Generalizing the linear transform of the example used in \cite{cdet},
we here introduce the modified process
\begin{equation}
\label{eq:7:1:3}
\bar{E}_t = E_t - {\mathbb E} \biggl[ \int_t^T f(P_s,0) ds | {\mathcal F}_t \biggr].
\end{equation}
In some sense, $\bar{E}_t$ gives an approximation of $E_T$ given ${\mathcal F}_t$: the dependence of $f$ upon $Y$ from time $t$ onward is frozen at $0$, that is $E_T$ is approximated by $E_t - \int_t^T f(P_s,0) ds$ and the conditional expectation provides the best least squares approximation of the resulting \emph{frozen} version of $E_T$ at time $t$. (See Footnote\footnote{One might think that using
$(f(P_s,Y_t))_{t \leq s < T}$ would provide a better approximation of $(f(P_s,Y_s))_{t\leq s < T}$. We prefer the simpler version $(f(P_s,0))_{t \leq s < T}$ because it does not depend upon $(E_s)_{t_0 \leq s \leq t}$.})  In particular, $\bar{E}_T= E_T$.  The coefficients of the dynamics of $P$ being Lipschitz continuous, the conditional expectation appearing in \reff{eq:7:1:3} is given by the deterministic function 
$w: [0,T] \times \RR^d \hookrightarrow \RR$ defined as the expectation
\begin{equation}
\label{eq:19:12:1}
w(t,p) = - {\mathbb E} \biggl[ \int_t^T f(P_s^{t,p},0) ds \biggr] = - {\mathbb E} \biggl[ \int_0^{T-t} f(P_s^{0,p},0) ds \biggr]
\end{equation}
over the solution for the dynamics of $P$ starting from $p$ at time $t$ (or from $p$ at time $0$ by time homogeneity). 

When the coefficients $b$, $\sigma$ and $f$ are smooth (with bounded derivatives of any order), the function $w$ is a classical solution of the PDE:
\begin{equation}
\label{eq:09:03:10}
\partial_t w(t,p) + \frac{1}{2} {\rm Trace} \bigl[a(p) \partial^2_{pp} 
w(t,p) \bigr] + \langle b(p),\partial_p w(t,p) \rangle - f(p,0)=0,
\end{equation}
with $0$ as terminal condition. 
(By \eqref{eq:19:12:1}, $w$ is once continuously differentiable in time; by differentiating the flow associated
with the process $P$ w.r.t. the variable $p$, it is also twice continuously differentiable in space. Moreover, by the standard dynamic programming principle, $w$ is a viscosity solution to the PDE \eqref{eq:09:03:10}. Therefore, it is a classical solution.) Consequently, for a given $0 \leq t_0 < T$, $(\bar{E}_t)_{t_0 \leq t \leq T}$ is an It\^o process with 
\begin{equation}
\label{eq:26:2:5}
\begin{split}
d \bar{E}_t &= d E_t + d\bigl[ w(t,P_t) \bigr]
\\
&= - \bigl[ f(P_t,Y_t) - f(P_t,0) \bigr] dt + \langle \sigma^{\top}(P_t) \partial_p w(t,P_t),dW_t \rangle, \quad t_0 \leq t < T,
\end{split}
\end{equation}
as dynamics. 

In any case, the process $(M^{t_0}_t)_{t_0 \leq t \leq T}$ defined by
\begin{equation*}
M^{t_0}_t=w(t,P_t) - \int_{t_0}^t f(P_s,0)ds = - {\mathbb E} \biggl[ \int_{t_0}^T f(P_s,0) ds | {\mathcal F}_t \biggr]
\end{equation*}
is a square integrable martingale on $[t_0,T]$. By the martingale representation theorem, it can be written as $\int_{t_0}^t \langle \theta_s,dW_s \rangle$ for some $\RR^d$-valued square integrable adapted process $\theta=(\theta_t)_{t_0\le t\le T}$. Eq. \reff{eq:26:2:5} shows that $\theta_t=\sigma^{\top}(P_t) \partial_p w(t,P_t)$ when $b$, $\sigma$ and $f$ are smooth (with bounded derivatives of any order).
In the general case, we will still use the same notation $\sigma^{\top}(P_t) \partial_p w(t,P_t)$ for the integrand appearing in the martingale representation of $M^{t_0}$ as a stochastic integral with respect to $W$, even if the gradient doesn't exist as a true function.

In particular, when 
the coefficients $b,\sigma,f$ satisfy (A.1), (A.2) and (A.3) only,
It\^o's formula \eqref{eq:26:2:5} holds as well as the expansion of $d[E_t+ \int_{t_0}^t f(P_s,0) ds + M_t^{t_0}]$.
As already explained, we always write $(\sigma^{\top}(P_t) \partial_p w(t,P_t))_{t_0 \leq t \leq T}$ for the 
integrand of the martingale part. Notice that, in any case, this integrand is bounded:
\begin{lemma}
\label{lem:22:12:3}
Under (A.1) and (A.2) only, there exists a constant $C$, depending on $L$ and $T$ only, such that 
\begin{equation}
\label{fo:Lw}
\forall (t,p,p') \in [0,T) \times \RR^d \times \RR^d, \quad 
|w(t,p')-w(t,p)| \leq C (T-t) |p'-p|. 
\end{equation}
In particular, when it exists, the function 
$\partial_p w(t,\cdot)$ is uniformly bounded from above by $C(T-t)$. And, in any case, the representation term $(\sigma^{\top}(P_t) \partial_p w(t,P_t))_{t_0 \leq t < T}$ is bounded by $C L(T-t)$ provided $\sigma$ is bounded by $L$. 
\end{lemma}

{\bf Proof.} The Lipschitz property \reff{fo:Lw} follows from the definition \reff{eq:19:12:1} of $w$, of the Lipschitz property of $f$ and of the $L^q(\Omega)$-Lipschitz property of the flow associated with $(P_s)_{t_0 \leq t \leq T}$, $q\geq 1$. 
Hence, when the coefficients are smooth, the integrand in the martingale representation is bounded, the bound depending on $L$ and $T$ only. By mollification, the bound remains true in the general case. \qed

\vskip 6pt\noindent
{\bf Notation.} From now on, 
$v$ denotes the value function 
in Proposition \ref{prop:7:11:10}. Moreover, we adopt the following convention. For 
$(t,p,e) \in [0,T) \times \RR^d \times \RR$, the notation $\bar{e}$ denotes
$\bar{e} = e +w(t,p)$. In particular, given $(t,p,\bar{e}) \in [0,T) \times \RR^d \times \RR$, $e$ is understood as $e=\bar{e} - w(t,p)$: quite often, we are given $\bar{e}$
first so that the value of $e$ follows. 

\subsection{Affine Feedback}
We first consider the case  
\begin{equation}
\label{fo:affine}
f(p,y) = f_0(p) + \ell y,
\end{equation}
for a given $\ell \in [\ell_1,\ell_2]$, $f_0$ being continuously differentiable. The need for the analysis of this particular case comes from the specific choice \reff {eq:7:1:3} of the approximation  $\bar E$ of $E$. Let us set
\begin{equation*}
\psi(e) = e {\mathbf 1}_{[0,1]}(e) + {\mathbf 1}_{(1,+\infty)}(e), \quad
e \in \RR,
\end{equation*}
so that the function $(t,e) \in [0,T) \times \RR \hookrightarrow \psi\big(e/(T-t)\big)$ is the continuous solution of the inviscid
Burgers' equation
\begin{equation*}
\partial_t u(t,e) - u(t,e) \partial_e u(t,e) = 0, \quad (t,e) \in [0,T) \times \RR,
\end{equation*}
with $u(T,\cdot) = {\mathbf 1}_{[0,+\infty)}$ as terminal condition. See \cite{cdet} or Lax \cite{Lax}.
By a change of variable, the function $e \hookrightarrow  \psi\big(\ell^{-1}(e-\Lambda)/(T-t)\big)$ satisfies the inviscid Burgers' equation
\begin{equation*}
\partial_t u(t,e) - \ell u(t,e) \partial_e u(t,e) = 0, \quad (t,e) \in [0,T) \times \RR,
\end{equation*}
with $u(T,\cdot) = {\mathbf 1}_{[\Lambda,+\infty)}$ as terminal condition. 
The specific affine form \reff{fo:affine} of the feedback function $f$ implies that \reff{eq:26:2:5} 
becomes:
\begin{equation}
\label{fo:debar}
d \bar{E}_t = -\ell Y_t dt +
\langle \sigma^{\top}(P_t)  \partial_p w(t,P_t),dW_t\rangle,
\quad t \in [0,T).
\end{equation}

We have:
\begin{lemma}
\label{lem:22:12:1}
There exists a constant $C$, depending on $L$ and $T$ only, such that
\begin{equation*}
\bigl|v(t_0,p,e) - \psi \left( \frac{\bar{e}-\Lambda}{\ell(T-t_0)} \right) \bigr|
\leq C (T-t_0)^{1/4}, \quad (t_0,p,e) \in [0,T) \times \RR^d \times \RR.
\end{equation*}
\end{lemma}
Recall that $\bar{e}$ stands for $\bar{e} =e + w(t_0,p)= e - {\mathbb E} \int_{t_0}^T f_0(P_s^{t_0,p}) ds$.

\vskip 2pt\noindent
{\bf Proof.} Since a similar bound was given in the proof of Proposition 4 in \cite{cdet}, we
only give a sketch of the proof. 

\textit{First Step.}
We prove that there exists a constant $c>0$, depending on $L$ and $T$ only, such that
\begin{equation}
\label{eq:26:2:6}
\begin{split}
&v(t_0,p,e) \geq 1 - \exp \bigl( - c \frac{\delta^2}{(T-t_0)^3} \bigr), \quad 
\bar{e} \geq \Lambda + \ell(T-t_0) + \delta,
\\
&v(t_0,p,e) \leq \exp \bigl( - c \frac{\delta^2}{(T-t_0)^3} \bigr), \quad 
\bar{e} \leq \Lambda - \delta.
\end{split}
\end{equation}
For the proof of the first inequality in \eqref{eq:26:2:6} we notice that \reff{fo:debar}, with the obvious initial condition 
$(P_{t_0}^{t_0,p},E_{t_0}^{t_0,p,e})=(p,e)$, implies
\begin{equation*}
\bar{E}_T \geq \bar{e} - \ell (T-t_0) + \int_{t_0}^T \langle
\sigma^{\top}(P_s) \partial_p w(s,P_s),dW_s \rangle.
\end{equation*}
So, when $\bar{e} \geq \Lambda + \ell(T-t_0) + \delta$, 
\begin{equation*}
1 - v(t_0,p,e) \leq \PP\{E_T \leq \Lambda\}
\leq {\mathbb P}
\biggl\{
\int_{t_0}^T \langle
\sigma^{\top}(P_s) \partial_p w(s,P_s),dW_s \rangle \leq -\delta \biggr\}.
\end{equation*}
By the the bound we have for the integrand  $(\sigma^{\top}(P_s) \partial_p w(s,P_s))_{t_0 \leq s <T}$ in the martingale representation in Lemma \ref{lem:22:12:3}, the bracket of the stochastic integral above is less than $C(T-t_0)^3$. By the exponential inequality
for continuous martingales, we complete the proof of the first inequality in \eqref{eq:26:2:6}. A similar argument gives the second inequality.

\vspace{5pt}\noindent
\textit{Second Step.} Next we prove:
\begin{equation*}
\begin{split}
&v(t_0,p,e) \geq \frac{\bar{e}-\Lambda}{\ell(T-t_0)} - \exp\bigl( - \frac{c \ell^2}{(T-t_0)^{1/2}}
\bigr) - (T-t_0)^{1/4}, \ \bar{e} < \Lambda + \ell (T-t_0) + \ell (T-t_0)^{5/4},
\\
&v(t_0,p,e) \leq \frac{\bar{e}-\Lambda}{\ell(T-t_0)} +  \exp\bigl( - \frac{c \ell^2}{(T-t_0)^{1/2}}
\bigr) + (T-t_0)^{1/4}, \ \bar{e} >\Lambda - \ell (T-t_0)^{5/4}.
\end{split}
\end{equation*}
Again, we prove the first inequality only. 
Choosing $\ell_1=\ell_2=\ell$ in the statement of Proposition \ref{prop:4:11:1}, we deduce that $v(t_0,\cdot,\cdot)$ is $1/[\ell(T-t_0)]$-Lipschitz w.r.t. $e$, so that
\begin{equation*}
v\bigl(t_0,p,\Lambda + \ell (T-t_0) + \ell (T-t_0)^{5/4} \bigr) - v(t_0,p,\bar{e})
\leq \frac{\Lambda - \bar{e}}{\ell(T-t_0)} + 1 + (T-t_0)^{1/4}.
\end{equation*}
Using the first step to bound from below 
$v(t_0,p,\Lambda + \ell (T-t_0) + \ell (T-t_0)^{5/4})$ in a similar way, we complete the proof of
the second step.

\vspace{5pt}\noindent
\textit{Third Step.} The proof is easily completed by using the second step when $\Lambda - \ell (T-t_0)^{5/4} < \bar{e} < \Lambda + \ell (T-t_0) + \ell (T-t_0)^{5/4}$ and the first step when
$\bar{e} \leq \Lambda - \ell (T-t_0)^{5/4}$ or
$\bar{e} \geq \Lambda + \ell (T-t_0) + \ell (T-t_0)^{5/4}$.
\qed

\subsection{Comparison with Burgers' Equation: General Case.}
We now  generalize Lemma \ref{lem:22:12:1} to the case of feedback functions $f$ of general form:
\begin{proposition}
\label{prop:26:2:1}
There exists a constant $C$ and an exponent $\beta \in (0,1)$, depending on $\alpha$, $L$ and $T$ only, such that
\begin{equation*}
\forall (t_0,p,e) \in [0,T) \times \RR^d \times \RR,
\quad \biggl|v(t_0,p,e) - \psi\biggl( 
\frac{\bar{e}-\Lambda}{\ell(t_0,p,e)[T-t_0]} \biggr) \biggr| \leq C (T-t_0)^{\beta},
\end{equation*}
where $\displaystyle \ell(t_0,p,e) =  \int_0^1 \frac{\partial f}{\partial y}\bigl(p,\lambda 
v(t_0,p,e) \bigr) d \lambda$.
\end{proposition}
\noindent
Note that by definition we have
\begin{equation}
\label{fo:vl}
v(t_0,p,e)\ell(t_0,p,e)=f(p,v(t_0,p,e))- f(p,0).
\end{equation}
{\bf Proof.} The proof is based on a comparison argument allowing us to \emph{piggy-back} on the affine case studied  above. Given an initial condition $(t_0,p,e)$, we set
$\ell = \ell(t_0,p,e) \in [\ell_1,\ell_2]$.
For a small $\varepsilon >0$, we consider the regularized systems
\begin{equation}
\label{eq:26:2:10}
\begin{split}
&dP_t^{\varepsilon} = b(P_t^{\varepsilon}) dt + \sigma(P_t^{\varepsilon}) dW_t + \varepsilon dW_t'
\\
&dE_t^{\varepsilon} =  - f(P_t^{\varepsilon},Y_t^{\varepsilon}) dt + \varepsilon dB_t,
\\
&dY_t^{\varepsilon} = \langle Z_t^{\varepsilon}, dW_t \rangle + \langle Z_t^{\varepsilon,\prime}, dW_t' \rangle + \Upsilon_t^{\varepsilon}dB_t,
\quad t_0 \leq t \leq T,
\end{split}
\end{equation}
and
\begin{equation}
\label{eq:26:2:11}
\begin{split}
&dP_t^{\varepsilon} = b(P_t^{\varepsilon}) dt + \sigma(P_t^{\varepsilon}) dW_t + \varepsilon dW_t'
\\
&dE_t^{\varepsilon,\ell} =  - \ell {Y}_t^{\varepsilon,\ell} dt - f(P_t^{\varepsilon},0) dt + \varepsilon dB_t,
\\
&dY_t^{\varepsilon,\ell} = \langle Z_t^{\varepsilon,\ell}, dW_t \rangle + \langle Z_t^{\varepsilon,\ell,\prime}, dW_t' \rangle + \Upsilon_t^{\varepsilon,\ell}dB_t,
\quad t_0 \leq t \leq T,
\end{split}
\end{equation}
with
$Y_T^{\varepsilon} = \phi(E_T^{\varepsilon})$ and
 $Y_T^{\varepsilon,\ell} = \phi(E_T^{\varepsilon,\ell})$ as terminal conditions, $\phi$ standing for a smooth non-decreasing function with values in $[0,1]$ (understood as an approximation
of the Heaviside funtion ${\mathbf 1}_{[\Lambda,+\infty)}$), and as before, $W'$ and $B$ being independent Brownian motions also independent of $W$. The associated value functions are denoted by $v^{\varepsilon}$ and $v^{\varepsilon,\ell}$. 
(See \eqref{eq:7:11:10}.)
The function $v^{\varepsilon,\ell}$ satisfies the PDE
\begin{equation*}
\bigl[ 
\partial_t v^{\varepsilon,\ell}+ {\mathcal L}_p v^{\varepsilon,\ell} + \frac{\varepsilon^2}{2} \Delta_{pp}
v^{\varepsilon,\ell} + \frac{\varepsilon^2}{2} \partial^2_{ee}
v^{\varepsilon,\ell} \bigr](t,p,e)
 - \bigl[ f(p,0) + \ell v^{\varepsilon,\ell}(t,p,e)
\bigr] \partial_e v^{\varepsilon,\ell}(t,p,e) = 0,
\end{equation*}
with $\phi$ as terminal condition. Compute now
$d [ v^{\varepsilon,\ell}(t,P_t^{\varepsilon},E_t^{\varepsilon})]$.  Following \reff{eq:22:12:1}, we obtain
\begin{equation*}
d \bigl[v^{\varepsilon,\ell}(t,P_t^{\varepsilon},E_t^{\varepsilon})\bigr]
= \bigl[ f(P_t^{\varepsilon},0) -f(P_t^{\varepsilon},Y_t^{\varepsilon})  + \ell v^{\varepsilon,\ell}(t,P_t^{\varepsilon},E_t^{\varepsilon}) \bigr]
\partial_e v^{\varepsilon,\ell}(t,P_t^{\varepsilon},E_t^{\varepsilon}) dt + dm_t,
\end{equation*}
where we use the notation $(m_t)_{0 \leq t \leq T}$ for a generic martingale which can change from one formula to the next.
Up to a modification of $(m_t)_{0 \leq t \leq T}$, we deduce
\begin{equation*}
d \bigl[v^{\varepsilon,\ell}(t,P_t^{\varepsilon},E_t^{\varepsilon})-Y_t^{\varepsilon}\bigr]
= \bigl[ f(P_t^{\varepsilon},0) -f(P_t^{\varepsilon},Y_t^{\varepsilon})  + \ell v^{\varepsilon,\ell}(t,P_t^{\varepsilon},E_t^{\varepsilon}) \bigr]
\partial_e v^{\varepsilon,\ell}(t,P_t^{\varepsilon},E_t^{\varepsilon}) dt + dm_t,
\end{equation*}
with $0$ as terminal condition at time $T$. 
This may also be written as
\begin{equation*}
\begin{split}
d \bigl[v^{\varepsilon,\ell}(t,P_t^{\varepsilon},E_t^{\varepsilon})-Y_t^{\varepsilon}\bigr]
&=   \ell \bigl[ v^{\varepsilon,\ell}(t,P_t^{\varepsilon},E_t^{\varepsilon}) - Y_t^{\varepsilon} \bigr]\partial_e v^{\varepsilon,\ell}(t,P_t^{\varepsilon},E_t^{\varepsilon}) dt
\\
&\hspace{15pt} + 
\biggl[ \ell - \int_0^1 \partial_y f(P_t^{\varepsilon},\lambda Y_t^{\varepsilon}) d\lambda
 \biggr] Y_t^{\varepsilon} \partial_e v^{\varepsilon,\ell}(t,P_t^{\varepsilon},E_t^{\varepsilon}) dt + dm_t, 
 \end{split}
\end{equation*}
$0 \leq t \leq T$.
Clearly,
\begin{equation*}
\begin{split}
&v^{\varepsilon,\ell}(t_0,p,e) - Y_{t_0}^{\varepsilon}
\\
&= {\mathbb E}
\biggl[ \int_{t_0}^T 
\biggl( \ell - \int_0^1 \partial_y f(P_t^{\varepsilon},\lambda Y_t^{\varepsilon}) d\lambda
 \biggr) Y_t^{\varepsilon} \partial_e v^{\varepsilon,\ell}(t,P_t^{\varepsilon},E_t^{\varepsilon}) 
\exp \biggl( - \ell \int_{t_0}^t 
\partial_e v^{\varepsilon,\ell}(s,P_s^{\varepsilon},E_s^{\varepsilon}) ds \biggr)
dt \biggr].
\end{split}
\end{equation*}
Therefore,
\begin{equation*}
\begin{split}
&\bigl|v^{\varepsilon,\ell}(t_0,p,e) - Y_{t_0}^{\varepsilon} \bigr| 
\\
&\leq {\mathbb E} \biggl[ \sup_{t_0 \leq t \leq T}
\biggl(  \int_0^1 \bigl| \ell - \partial_y f(P_t^{\varepsilon},\lambda Y_t^{\varepsilon}) \bigr| d\lambda \biggr)
\\
&\hspace{15pt} \times
\int_{t_0}^T \partial_e v^{\varepsilon,\ell}(t,P_t^{\varepsilon},E_t^{\varepsilon}) 
\exp \biggl( - \ell \int_{t_0}^t 
\partial_e v^{\varepsilon,\ell}(s,P_s^{\varepsilon},E_s^{\varepsilon}) ds \biggr)
dt \biggr]
\\
&= \ell^{-1} {\mathbb E}
\biggl\{ \sup_{t_0 \leq t \leq T}  \biggl|\ell - \int_0^1  \partial_y f(P_t^{\varepsilon},\lambda Y_t^{\varepsilon}) d\lambda  \biggr|
\times
\biggl[1- \exp
\biggl( - \ell
\int_{t_0}^T 
\partial_e v^{\varepsilon,\ell}(s,P_s^{\varepsilon},E_s^{\varepsilon}) ds \biggr)
 \biggr] \biggr\},
\end{split}
\end{equation*}
and finally,
\begin{equation*}
\bigl|(v^{\varepsilon,\ell} - v^{\varepsilon})(t_0,p,e) \bigr| 
\leq \ell^{-1} {\mathbb E} \biggl\{ \sup_{t_0 \leq t \leq T}
\biggl[ \int_0^1 \bigl| 
\partial_y f(P_t^{\varepsilon},\lambda Y_t^{\varepsilon})
- \partial_y f\bigl(p,\lambda v(t_0,p,e)\bigr) \bigr| d\lambda \biggr] \biggr\}.
\end{equation*}
Using the H\"older continuity of $\partial_y f$, we get
\begin{equation*}
\bigl|(v^{\varepsilon,\ell} - v^{\varepsilon})(t_0,p,e) \bigr| 
\leq C {\mathbb E} \bigl\{ \sup_{t_0 \leq t \leq T} \bigl[
|P_t^{\varepsilon} - p|^{\alpha} + |Y_t^{\varepsilon} - Y_{t_0}^{\varepsilon}|^{\alpha}\bigr] \bigr\} 
 + C|(v^{\varepsilon}- v)(t_0,p,e)|^{\alpha}.
\end{equation*}
Following the proof of Lemma \ref{lem:7:11:2}, we let $\varepsilon$ tend first to $0$. By Corollary \ref{cor:25:7:1}, 
$v^{\varepsilon}$ converges towards $v^{\phi}$
and $v^{\varepsilon,\ell}$ converges towards $v^{\phi,\ell}$, convergences being uniform on compact subsets of $[0,T] \times \RR^d \times \RR$ and
$v^{\phi}$ and $v^{\phi,\ell}$ standing for the value functions associated with 
\eqref{eq:26:2:10} and \eqref{eq:26:2:11} when $\varepsilon=0$ therein. By the gradient bound \eqref{eq:20:07:3},
the integrand in the martingale representation of $(Y_t^{\phi})_{t_0 \leq t \leq T}$
is bounded, independently of $\phi$, so that the increments of 
$(Y_t^{\phi})_{t_0 \leq t \leq T}$ are well-controlled. We deduce that
\begin{equation}
\label{eq:23:2:1}
\bigl|(v^{\phi,\ell}  - v^{\phi})(t_0,p,e) \bigr| 
\leq C(T-t_0)^{\alpha/2} + C \bigl|(v^{\phi}  - v)(t_0,p,e) \bigr|^{\alpha},
\end{equation}
for a constant $C$ depending on $\alpha$, $L$ and $T$ only.  
As $\phi$ converges towards the Heaviside function ${\mathbf 1}_{[\Lambda,+\infty)}$ as in Section \ref{se:generalexistence}, we know from Corollary \ref{cor:27:7:1} that 
$v^{\phi}$ converges towards $v$ and $v^{\phi,\ell}$ towards 
$v^{\ell}$, where $v^{\ell}$ is the value function associated with \eqref{eq:26:2:11} when $\varepsilon=0$, but with 
${\mathbf 1}_{[\Lambda,+\infty)}$ as terminal condition.
Applying Lemma \ref{lem:22:12:1} (with $f_0(p) = f(p,0)$)
to estimate $v^{\ell}$, we complete the proof.
\qed

\subsection{Proof of the Existence of a Dirac Mass}
We claim

\begin{proposition}
\label{prop:30:12:1}
There exists a constant $c \in (0,1)$, depending on $\alpha$ and $L$ only,
such that, if $T-t_0 \leq c$, $p \in \RR^d$ and $(\bar{e}-\Lambda)/(T-t_0) \in [\ell_1/4,3\ell_1/4]$, then:
\begin{equation}
\label{fo:ptmass}
\PP\{E_T^{t_0,p,e}=\Lambda\} \geq c. 
\end{equation}
\end{proposition}

\begin{remark}
\label{remark:1}
We emphasize that, in general, Proposition \ref{prop:30:12:1} cannot be true for all starting point $(t_0,p,e)$. Indeed, in the non-viscous case, i.e. when $E$ doesn't depend upon $P$ (or, equivalently, when $f$ depends on $Y$ only),
the dynamics of 
$E$ coincide with the dynamics of the characteristics of the associated inviscid equation of conservation law. The typical example is the Burgers' equation: the characteristics satisfy the equation
\begin{equation*}
dE_t = - \psi \bigl( \frac{E_t-\Lambda}{T-t} \bigr) dt,
\end{equation*}
and consequently,
\begin{equation*}
E_t = 
\left\{\begin{array}{ll}
\displaystyle e - (t-t_0) \quad &{\rm when} \ e - (T-t_0) > \Lambda,
\vspace{5pt}
\\
\displaystyle e \quad &{\rm when} \ e < \Lambda,
\\
\displaystyle e - \frac{e-\Lambda}{T-t_0}(t-t_0) \quad &{\rm when} \
\Lambda \leq e \leq \Lambda + (T-t_0),
\end{array}
\right.
\end{equation*}
where $(t_0,e)$ stands for the initial condition of the process $(E_t)_{t_0 \leq t \leq T}$, i.e. $E_{t_0}=e$. 
This corresponds to the following picture:

\unitlength .7cm
\begin{center}
\begin{picture}(10.5,10.5)(-0.5,-0.5)%
\put(8.75,2.4){\makebox(0,0){\textcolor{black}{$T$}}}
\color{black}
\path(0.5026,6.0053)(8.7566,6.0053)
\path(8.7434,9.2328)(8.7434,2.7646)
\path(8.0026,9.246)(8.7434,8.5053)
\path(6.4947,9.246)(8.7434,6.9841)
\path(8.7434,6.0053)(5.4894,9.246)
\path(6.0053,9.246)(8.7434,6.4947)
\path(6.9974,9.246)(8.7434,7.5)
\path(7.5,9.246)(8.7434,8.0026)
\path(5,9.246)(8.7434,6.0053)
\path(4.4974,9.246)(8.7434,6.0053)
\path(3.9947,9.246)(8.7434,6.0053)
\path(3.5053,9.246)(8.7434,6.0053)
\path(3.0026,9.246)(8.7434,6.0053)
\path(2.5,9.246)(8.7434,6.0053)
\path(1.9974,9.246)(8.7434,6.0053)
\path(1.4947,9.246)(8.7434,6.0053)
\path(1.0053,9.246)(8.7434,6.0053)
\path(0.5026,9.246)(8.7434,6.0053)
\path(0.5026,8.9947)(8.7434,6.0053)
\path(0.5026,8.7566)(8.7434,6.0053)
\path(0.5026,8.5053)(8.7434,6.0053)
\path(0.5026,8.254)(8.7434,6.0053)
\path(0.5026,8.0026)(8.7434,6.0053)(0.5026,7.7513)
\path(0.5026,7.5)(8.7434,6.0053)(0.5026,7.2487)
\path(0.5026,6.9974)(8.7434,6.0053)(0.5026,6.746)
\path(0.5026,6.4947)(8.7434,6.0053)(0.5026,6.2434)
\put(9.1,6.1053){\makebox(0,0){\textcolor{black}{$\Lambda$}}}
\color{black}
\path(0.5026,5.5026)(8.7434,5.5026)
\path(0.5026,5)(8.7434,5)
\path(0.5026,4.4974)(8.7434,4.4974)
\path(0.5026,3.9947)(8.7434,3.9947)
\path(8.7434,3.5053)(0.5026,3.5053)
\path(0.5026,3.0026)(8.7434,3.0026)
\end{picture}%
\vspace{-30pt}

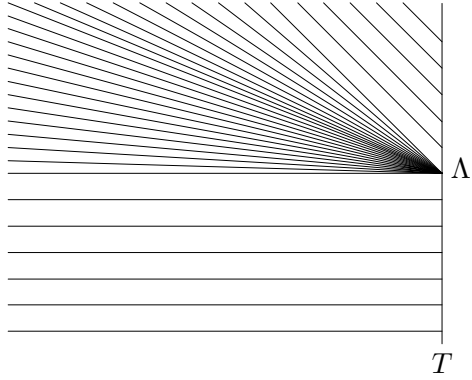
\textit{Fig 1. Characteristics of the inviscid Burgers' equation.}
\vspace{5pt}
\end{center}
Clearly, the singular point $\Lambda$ is hit when  $\Lambda \leq e \leq \Lambda + (T-t_0)$ only. 
In order for $\Lambda$ to be hit starting from $e$ outside the cone shown in the figure, noise must be plugged into the system, i.e. noise must be transmitted from the first to the second equation. This point is investigated in the next subsection.
\end{remark}

\noindent
{\bf Proof of Proposition \ref{prop:30:12:1}.}
Given an initial condition $(t_0,p,e) \in [0,T) \times \RR^d \times \RR$ for the process $(P,E)$, we consider the stochastic differential equations
\begin{equation}
\label{eq:01:03:1}
\begin{split}
d \bar{E}_t^{\pm} &= \left( -\ell(t,P_t,E_t) \psi \bigl[ \ell^{-1}(t,P_t,E_t)
\frac{E_t^{\pm} - \Lambda}{T-t} \bigr] \pm C'(T-t)^{\beta} \right) dt
\\
&\hspace{15pt} + \langle \sigma^{\top}(P_t) \partial_p w(t,P_t), dW_t \rangle,
\end{split}
\end{equation}
with $\bar{E}_{t_0}^{\pm}=\bar{e}$ as initial conditions, the constant $C'$ being chosen later on. Notice that 
the process appearing in $\ell$ and $\ell^{-1}$ above is $E$ and not $\bar{E}^{\pm}$.
From \eqref{eq:26:2:5} and \reff{fo:vl} it follows that
\begin{equation*}
d \bar{E}_t = - \ell(t,P_t,E_t) v(t,P_t,E_t) dt 
+ \langle \sigma^{\top}(P_t) \partial_p w(t,P_t),dW_t \rangle, \quad t \in [t_0,T),
\end{equation*}
with
\begin{equation*}
\bigl|\ell(t,P_t,E_t) v(t,P_t,E_t) - \ell(t,P_t,E_t) \psi
\bigl( \ell^{-1}(t,P_t,E_t) \frac{\bar{E}_t - \Lambda}{T-t} \bigr)
\bigr| \leq L C (T-t)^{\beta}, \quad t \in [t_0,T),
\end{equation*}
where $C$ is given by Proposition \ref{prop:26:2:1}.
We now choose $C'=LC$. By the comparison theorem for one-dimensional SDEs, we deduce
\begin{equation*}
\bar{E}_t^{-} \leq \bar{E}_t \leq \bar{E}_t^+, \qquad t \in [t_0,T]. 
\end{equation*}
Next, we introduce the bridge equations
\begin{equation}
\label{eq:01:03:2}
d \bar{Z}_t^{\pm} = \left( - \frac{\bar{Z}_t^{\pm} - \Lambda}{T-t} \pm C'(T-t)^{\beta}
\right) dt
+ \langle \sigma^{\top}(P_t) \partial_p w(t,P_t), dW_t \rangle, \quad
\bar{Z}_{t_0}^{\pm} = \bar{e}.
\end{equation}
The solution is given by
\begin{equation}
\label{eq:14:2:10}
\bar{Z}_t^{\pm} = \Lambda + (T-t) \biggl[ \frac{\bar{e}-\Lambda}{T-t_0}
\pm C' \int_{t_0}^t (T-s)^{\beta-1} ds + \int_{t_0}^t
(T-s)^{-1}
\langle \sigma^{\top}(P_s) \partial_p w(s,P_s),dW_s \rangle \biggr],
\end{equation}
so that $\bar{Z}_t^{\pm} \rightarrow \Lambda$ as $t \rightarrow T$. 
(The stochastic integral is well-defined up to time $T$
by Lemma \ref{lem:22:12:3}.)

Now, we choose $\bar{e}$ such that $(\bar{e}-\Lambda)/(T-t_0) \in [\ell_1/4,3\ell_1/4]$ and $t_0$ such that
$C' \int_{t_0}^T (T-s)^{\beta-1}ds \in [0,\ell_1/16]$, and we introduce the stopping time
\begin{equation}
\label{eq:08:03:2}
\tau = \inf \biggl\{t \geq t_0 :  \biggl| \int_{t_0}^t
(T-s)^{-1}
\langle \sigma^{\top}(P_s) \partial_p w(s,P_s), dW_s \rangle 
\biggr| \geq \frac{\ell_1}{16} \biggr\} \wedge T.
\end{equation}
We obtain
\begin{equation*}
\frac{\ell_1}{8}\leq \frac{\bar{Z}^{\pm}_t - \Lambda}{T-t} \leq \frac{7\ell_1}{8},
\end{equation*}
for any $t \in [t_0,\tau)$,
so that 
\begin{equation}
\label{eq:08:03:1}
\frac{\bar{Z}^{\pm}_t - \Lambda}{T-t}
= \ell(t,P_t,E_t) \psi \biggl[ 
\ell^{-1}(t,P_t,E_t) \frac{\bar{Z}_t^{\pm} - \Lambda}{T-t} \biggr],
\quad t_0 \leq t < \tau,
\end{equation}
in other words, $(\bar{Z}^{\pm}_t)_{t_0 \leq t < \tau}$ and
$(\bar{E}_t^{\pm})_{t_0 \leq t <\tau}$ coincide. (Compare \eqref{eq:01:03:1} and \eqref{eq:01:03:2}.) We deduce that, on the event $F = \{\tau=T\}$,
\begin{equation*}
\bar{Z}^{\pm}_t = \bar{E}^{\pm}_t, \quad t \in [t_0,T].
\end{equation*}
Finally, by Markov inequality and Lemma \ref{lem:22:12:3}, the probability of the event $F$ is strictly greater than zero for $T-t_0$ small enough.
This completes the proof. \qed

\begin{remark}
\label{rem:4:9:1}
We emphasize that the boundedness of $b$ in Assumption (A.4) plays a minor role in the proof of Proposition \ref{prop:26:2:1}. Basically, it is used in \eqref{eq:23:2:1} only to bound the increments of the process $P$. When $b$ is not bounded but at most of linear growth, the constant $C$ in \eqref{eq:23:2:1} may depend on $p$: in the end, the right-hand side in Proposition \ref{prop:26:2:1} has the form $C(1+|p|)(T-t_0)^{\beta}$, $C$ being independent of $p$. This affects the proof of Proposition \ref{prop:26:2:1} in the following way: to adapt the proof to the case when $b$ is at most of linear growth, the constant $C'$ in \eqref{eq:01:03:1} and \eqref{eq:01:03:2} must be changed into $C'(1+|P_t|)$. As a consequence, the term $C' \int_{t_0}^t (1+|P_s|)(T-s)^{\beta-1} ds$ 
in \eqref{eq:14:2:10} 
is small with large probability provided $T-t_0 \leq c$, $c$ being uniform w.r.t. the initial condition of the process $P$, namely $P_{t_0}=p$, in compact subsets of ${\mathbb R}^d$. Therefore, Proposition \ref{prop:26:2:1} still holds when $b$ is at most of linear growth provided the constant $c$ therein is assumed to be uniform w.r.t. $p$ in compact subsets only. 
\end{remark}

\subsection{Dirac Mass in the Non-Degenerate Regime}
We now discuss the attainability of the criticial area described in the statement of 
Proposition \ref{prop:30:12:1} when the initial point $\bar{e}$ does not belong to it. As noticed in Remark \ref{remark:1}, additional  assumptions of non-degeneracy type are necessary to let the critical area be attainable with a non-zero probability. 

\begin{proposition}
\label{prop:4:1:1}
In addition to (A.1), (A.2), (A.3) and (A.4), let us assume that
the noisy component of the forward process is uniformly elliptic in the sense that (up to a modification of $L$)
\begin{equation}
\label{eq:4:8:15}
\sigma^{\top}(p)\sigma(p)  \geq L^{-1} ,\qquad\qquad p \in\RR^d.
\end{equation}
Furthermore, let us also assume that for any $p \in {\mathbb R^d}$,
$|\partial_p f(p,0)| \geq L^{-1}$, and that 
$p \hookrightarrow \partial_p f(p,0)$ is uniformly continuous. 
Then, for any starting point $(t_0,p,e) \in [0,T) \times \RR^d \times \RR$,
$$
\PP\{E_T^{t_0,p,e}=\Lambda \}>0
$$
and the topological support of the conditional law of 
$Y_T^{t_0,p,e}$ given $E_T^{t_0,p,e}= \Lambda$ is $[0,1]$. 
\end{proposition}

{\bf Proof.} 
\textit{First Step. Positivity of $\PP\{E_T=\Lambda\}$.}
Since $\bar{E}_T=E_T$, it is enough to prove that, for any starting point 
$(t_0,p,\bar{e}) \in [0,T) \times \RR^d \times \RR$ of the process
$(t,P_t,\bar{E}_t)_{t_0 \leq t \leq T}$, 
$\PP\{\bar{E}_T= \Lambda\} >0$. (As usual, we omit below to specify the superscript
$(t_0,p,\bar{e})$ in $(t,P_t,\bar{E}_t)_{t_0 \leq t \leq T}$.)

By Proposition \ref{prop:30:12:1}, it is enough to prove that there exists $t$ close to $T$ such that $\bar{E}_t \in [\Lambda+\ell_1(T-t)/4,\Lambda+3\ell_1(T-t)/4]$
with a non-zero probability. Since the pair $(P,E)$ is a Markov process, we can assume $t_0$ itself to be close to $T$: we then aim at proving that, with a non-zero probability, 
the path $(\bar{E}_{t})_{t_0 \leq t \leq (T+t_0)/2}$ hits
the interval $[\Lambda+\ell_1(T-t)/4,\Lambda+3\ell_1(T-t)/4]$. 
It is sufficient to prove that, 
with a non-zero probability, $(\bar{E}_{t})_{t_0 \leq t \leq (T+t_0)/2}$ falls at least once into
the interval $[\Lambda+\ell_1(T-t_0)/4,\Lambda+3\ell_1(T-t_0)/8]$.

We first prove that the diffusion coefficient in \eqref{eq:26:2:5}
is away from zero on $[t_0,(T+t_0)/2]$.
By uniform ellipticity of $\sigma \sigma^{\top}$, $w$ is a classical solution to the PDE \eqref{eq:09:03:10}, so that $\partial_p w$ exists as a true function. (See Chapter 8 in \cite{Krylov:book}.) Moreover, $|\partial_p w(t,p)|^2$ is away from zero, uniformly in $(t,p) \in [t_0,(T+t_0)/2] \times \RR^d$, when $t_0$ is close enough to $T$. Indeed, going back
to \eqref{eq:19:12:1}, when the coefficients are smooth, 
$\partial_p w(t,p)$ is given by
\begin{equation}
\label{eq:09:03:5}
\partial_{p_i} w(t,p) = - {\mathbb E} \int_t^T \langle \partial_p f(P_s^{t,p},0),
\partial_{p_i} P_s^{t,p} \rangle ds, \quad i \in \{1,\dots,d\},
\end{equation}
and when $T-t$ is small, $\partial_p P_s^{t,p}$
is close to the identity (uniformly in $p$ since the coefficients $b$ and $\sigma$ are Lipschitz continuous) and $\partial_p f(P_s^{t,p},0)$ is close to 
$\partial_p f(p,0)$ (uniformly in $p$ since the coefficients $b$ and $\sigma$
are bounded and $\partial_p f(\cdot,0)$ is uniformly continuous). Therefore, $\partial_{p_i} w(t,p)$ is close to 
$\partial_{p_i} f(p,0)$, uniformly in $p$. Since the norm of
$\partial_p f(p,0)$ is away from zero, uniformly in $p$, we deduce that the same holds for $\partial_p w(t,p)$ when $T-t$ is small, uniformly in $p$. By a mollification argument, the result remains true under the assumption of Proposition \ref{prop:4:1:1}.

Therefore, $(\bar{E}_t)_{t_0 \leq t \leq (T+t_0)/2}$ is an It\^o process with a bounded drift and a (uniformly) non-zero diffusion coefficient: by Girsanov theorem and by change of time, the problem is equivalent to proving that a Brownian motion falls, with a non-zero probability, into a given interval of non-zero length before a given positive time, which is obviously true.
\vspace{5pt}

\textit{Second Step. Support of the Conditional Law.}
We investigate the support of the conditional law of 
$Y_T$ given $E_T= \Lambda$ or equivalently of $Y_T$ given 
$\bar{E}_T = \Lambda$.
The desired result is a consequence of the following two lemmas

\begin{lemma}
\label{lem:26:2:1}
Assume that
\begin{equation}
\label{eq:5:3:1}
f(p,0) - f(p,y) + \frac{\bar{e}-\Lambda}{T-t} = 0,
\end{equation}
for some $(t,p,\bar{e}) \in [0,T) \times \RR^d \times \RR$ and
$y \in [\varepsilon,1-\varepsilon]$, $\varepsilon \in (0,1)$. 
Then, there exists $\delta_1(\varepsilon)$, independent of $(t,p,\bar{e},y)$, such that
$T-t < \delta_1(\varepsilon)$ implies $v(t,p,e) \in (y-\varepsilon,y+\varepsilon)$, with 
$e=\bar{e}-w(t,p)$.
\end{lemma}

\begin{lemma}
\label{lem:26:2:2}
For any $\varepsilon >0$, there exists $\delta_2(\varepsilon)>0$, such that,
for any $y \in [0,1]$ and any $(t_0,p,e) \in [0,T) \times \RR^d \times \RR$
satisfying 
$|v(t_0,p,e)-y| \leq \varepsilon$ and $T-t_0 \leq \delta_2(\varepsilon)$, 
it holds 
$$
{\mathbb P}\{|Y_T^{t_0,p,e}- y| < 2 \varepsilon\} \geq 1/2.
$$
\end{lemma}

Here is the application of Lemmas \ref{lem:26:2:1} and \ref{lem:26:2:2} 
to the proof of Proposition \ref{prop:4:1:1}. Given $y,\varepsilon \in (0,1)$ such that
$(y-\varepsilon,y+\varepsilon) \subset (0,1)$, we are now proving that
$
{\mathbb P}\{Y_T^{t_0,p,e} \in (y-\varepsilon,y+\varepsilon)\} >0
$
for any initial condition $(t_0,p,e)$. (Below, we do not specify the superscript $(t_0,p,e)$.)
By the Markov property, we can assume $T-t_0 \leq \delta_1(\varepsilon/2) \wedge \delta_2(\varepsilon/2)$.
It is then sufficient to prove that, with a non-zero probability, 
the stopping time 
\begin{equation*}
\tau = \inf\bigl\{t \in [t_0,T] : f(P_t,0) - f(P_t,y) + \frac{\bar{E}_t-\Lambda}{T-t}
=0 \bigr\} \wedge T,
\end{equation*}
is in $[t,T)$. Indeed, by Lemma \ref{lem:26:2:1}, 
$\tau < T$ implies $|Y_{\tau} - y| < \varepsilon/2$; by Lemma \ref{lem:26:2:2} and the strong Markov property, this implies
${\mathbb P}\{|Y_T - y| < \varepsilon|{\mathcal F}_{\tau}\} \geq 1/2$, so that ${\mathbb P}\{\tau < T\} >0$ implies ${\mathbb P}\{
|Y_T - y| < \varepsilon\}>0$.

To prove that $\tau <T$ with a non-zero probability, we apply the following simple inequality
\begin{equation*}
0 \leq f(P_t,y) - f(P_t,0) \leq \ell_2, \quad t \in [t_0,T].
\end{equation*}
Assume indeed that we can find two times $t_1 < t_2 \in [t_0,T)$ such that
\begin{equation}
\label{eq:26:2:3}
\bar{E}_{t_1} > \Lambda + \ell_2 T, \qquad \text{and}, \qquad \bar{E}_{t_2} < \Lambda.
\end{equation}
Then,
\begin{equation*}
\frac{\bar{E}_{t_1} - \Lambda}{T-t_1} > f(P_{t_1},y) - f(P_{t_1},0),
\qquad \text{and}, \qquad \frac{\bar{E}_{t_2} - \Lambda}{T-t_2} < f(P_{t_2},y) - f(P_{t_2},0).
\end{equation*}
By continuity, there exists some $t \in (t_1,t_2)$ at which 
$(\bar{E}_t - \Lambda)/(T-t) = f(P_t,y) - f(P_t,0)$.

To prove \eqref{eq:26:2:3}, we follow the same strategy as in the first step.
The process $(\bar{E}_t)_{t_0 \leq t \leq (T+t_0)/2}$ is an It\^o process with a bounded drift and a uniformly non-zero diffusion coefficient: by Girsanov theorem and by a change of time, it is sufficient to prove that a Brownian motion starting from
$\bar{e}$ at time $t_0$ satisfies \eqref{eq:26:2:3} with a non-zero probability on an interval of small length, which is obviously true. \qed

\vspace{5pt}

{\bf Proof of Lemma \ref{lem:26:2:2}.}
For any $t \in [t_0,T]$,
\begin{equation*}
Y_t^{t_0,p,e} = v(t_0,p,e) + \int_{t_0}^t Z_s^{t_0,p,e} dB_s.
\end{equation*}
By Theorem \ref{thm:4:11:1}, $|Z_s^{t_0,p,e}| \leq C$, $C$ depending on $L$ and $T$ only, so that
${\mathbb E} [ |Y_t^{t_0,p,e} - v(t_0,p,e)|^2 ] \leq C(T-t_0)$.
In particular,
\begin{equation*}
{\mathbb P}\bigl\{ |Y_T^{t_0,p,e} - y| > 2 \varepsilon \bigr\}
\leq {\mathbb P}\bigl\{ |Y_T^{t_0,p,e} - v(t_0,p,e)| >  \varepsilon \bigr\}
\leq \frac{C}{\varepsilon^2} (T-t_0).
\end{equation*}
Choosing $T-t_0$ appropriately, we complete the proof.
\qed
\vspace{5pt}
\\

{\bf Proof of Lemma \ref{lem:26:2:1}.} In the proof, we denote $v(t,p,e)$ by $v$ and
$\ell(t,p,e)$ by $\ell(v)$. By Proposition \ref{prop:26:2:1}, we know that
\begin{equation}
\label{eq:26:2:2}
\bigl| v - \psi \bigl( \frac{\bar{e}-\Lambda}{\ell(v)(T-t)}  \bigr)
\bigr| \leq C(T-t)^{\beta},
\end{equation}
that is
\begin{equation*}
\bigl| v - \psi \bigl( \frac{F(y)}{\ell(v)} \bigr) \bigr| \leq C (T-t)^{\beta},
\end{equation*}
with $F(z) = f(p,z) - f(p,0)$, by \eqref{eq:5:3:1}. Multiplying by $\ell(v)$ and observing that 
$F(v) = \ell(v)v$, we deduce, for a new value of $C$,
\begin{equation}
\label{eq:26:2:1}
\bigl| F(v) - \ell(v) \psi \bigl( \frac{F(y)}{\ell(v)} \bigr) \bigr| \leq C (T-t)^{\beta}.
\end{equation}

\textit{First Case: $v \geq y$.} Since $F$ is increasing and $F(y),\ell(v) \geq 0$, we notice that
\begin{equation*}
F(v) - \ell(v) \psi \bigl( \frac{F(y)}{\ell(v)} \bigr)
\geq F(v) - F(y) \geq 0,
\end{equation*}
so that 
$F(v) - F(y) \leq C(T-t)^{\beta}$.
Since $\partial_z F(z) \geq \ell_1$, this proves that 
$v-y \leq C\ell_1^{-1}(T-t)^{\beta}$. Choosing $T-t$ small enough, we deduce that 
$0 \leq v-y < \varepsilon$.
\vspace{5pt}

\textit{Second Case: $v < y$.}
From the inequality
\begin{equation*}
\frac{\bar{e}- \Lambda}{\ell(v)(T-t)} = \frac{F(y)}{\ell(v)} \geq (\ell_1/\ell_2) y \geq (\ell_1/\ell_2) \varepsilon,
\end{equation*}
we deduce that (w.l.o.g. we can assume
$(\ell_1/\ell_2) \varepsilon < 1$)
\begin{equation*}
\psi \bigl( \frac{\bar{e}-\Lambda}{\ell(v)(T-t)} \bigr) - C(T-t)^{\beta}
\geq (\ell_1/\ell_2) \varepsilon  - C(T-t)^{\beta}.
\end{equation*}
Choosing $T-t$ small enough, we deduce that 
$\psi [( \bar{e}-\Lambda)/(\ell(v)(T-t))] >0$. 
Similarly, from \eqref{eq:26:2:2},
\begin{equation*}
1-\varepsilon \geq y \geq v \geq \psi \bigl( \frac{\bar{e}-\Lambda}{\ell(v)(T-t)} \bigr)
- C(T-t)^{\beta}.
\end{equation*}
Again, for $T-t$ small enough, we deduce that
\begin{equation*}
0 < \psi \bigl( \frac{\bar{e}-\Lambda}{\ell(v)(T-t)} \bigr) < 1,
\end{equation*}
so that
\begin{equation*}
\psi \bigl( \frac{F(y)}{\ell(v)} \bigr)=
\psi \bigl( \frac{\bar{e}-\Lambda}{\ell(v)(T-t)} \bigr)
= \frac{\bar{e}-\Lambda}{\ell(v)(T-t)} = \frac{F(y)}{\ell(v)}.
\end{equation*}
By \eqref{eq:26:2:1}, we obtain
$0 \leq F(y) - F(v) \leq C(T-t)^{\beta}$. We complete the proof as in the first case.
\qed

\section{Absolute Continuity before Terminal Time $T$}
\label{se:ac}
We have just shown that the paths of the process $E$ coalesce at the singular point
with a non-zero probability when \eqref{eq:4:11:1} is driven by a binary terminal condition. We here investigate the dynamics before terminal time $T$, a first objective being to put in evidence the existence of a transition in the system at time $T$, a second and more refined one being to establish the smoothness properties of the marginal distribution of the process $E$ before time $T$ in the cone limited by the characteristics of the inviscid regime. A quite simple evidence of the transition at the terminal time $T$ is the fact that the flow $e \hookrightarrow E^{0,p,e}$ loses its homeomorphic property at that time $T$.

\begin{proposition}
\label{prop:6:1:1}
Assume that (A.1--4) are in force and that $\phi = {\mathbf 1}_{[\Lambda,+\infty)}$
as in Section \ref{se:dirac}. Then, at any time $t<T$ and for any $p \in \RR^d$, the mapping $\RR\ni e \hookrightarrow E^{0,p,e}_t$ is an homeomorphism with probability $1$, and with non-zero probability, it is not a homeomorphism at time $t=T$. 
\end{proposition}

We chose $0$ as initial time for notation convenience only as it can clearly be shifted.
\vspace{5pt}

{\bf Proof.} The proof is a straightforward consequence of formula 
\eqref{eq:5:3:2} for mollified coefficients. In the limit, it says that 
\begin{equation}
\label{eq:2:8:2}
e-e'\geq E_t^{t_0,p,e} - E_t^{t_0,p,e'} \geq (e-e') \exp \biggl( - \frac{\ell_2}{\ell_1} \int_{t_0}^t (T-s)^{-1} ds \biggr)
 = \bigl( \frac{T-t}{T-t_0} \bigr)^{\ell_2/\ell_1} (e-e'), 
 \end{equation}
for  $e>e'$.
In particular, the mapping $ \RR\ni e \hookrightarrow E_t^{t_0,p,e}$ is continuous and increasing. \qed

\subsection{H\"ormander Property}
Proposition \ref{prop:6:1:1} describes the transition regime in rather simplistic terms. We now refine the analysis by studying how the process $E$ feels the noise coming from the diffusion process. It was proven in \cite{cdet} that the marginal laws of $E$ have densities before time $T$ when the diffusion process $P$ is a Brownian motion and the transmission function $f$ is linear in $p$ and $y$. The proof used Bouleau-Hirsch's criterion for the Malliavin derivative of $E$. There, the crucial step was to prove that
\begin{equation}
\label{eq:6:1:11}
- \partial_p \bigl[ f\bigl(p,u(t,p,e) \bigr) \bigr] >0,
\end{equation}
in the specific case $d=1$ and $-f(p,y)=p-y$. Equation \reff{eq:6:1:11} may be viewed as a non-degeneracy condition of H\"ormander type: if some noise is plugged into the process $P$ and if $P$ feels the noise itself, then the noise is transmitted to the process $E$, and $E$ oscillates whenever $P$ does. We refer to this condition as ``a first-order H\"ormander structure''.

Given \reff{eq:6:1:11}, we may ask two questions: (i) what happens in the general case? (ii) what is the typical size of the noise transmitted from $P$ to $E$ close to time $T$? As we are about to show, addressing these questions is far from trivial and the answers we give below are only partial.

In Proposition \ref{prop:20:07:1} we prove property \reff{eq:6:1:11} for some important cases,
 and in Proposition \ref{prop:08:03:1} we show that  \reff{eq:6:1:11} can fail in other cases. When it does, the process $E$ has pathological points inside the critical cone where either the noise is transmitted from $P$ to  $E$ in a singular way, or its marginal distribution is not absolutely continuous with respect to Lebesgue's measure. Here, ``singular way'' refers to a possible transmission of the noise from $P$ to $E$ through the second (or higher) order derivatives of $f(p,u(t,p,e))$, in other words, H\"ormander's condition is satisfied but because of Lie brackets of lengths greater than $2$. 

We also show in Proposition \ref{prop:7:1:1} that the linear case is critical: in equation \reff{eq:6:1:11}, the transmission coefficient is non-zero, but so small when $t$ approaches the terminal time $T$ that the noise propagation cannot be observed numerically: the typical fluctuations  of $E$ produced at the end of a time interval $[t,t+\varepsilon]$ are of order $c(t) \varepsilon^{3/2}$ with $c(t)$ decaying exponentially fast as $t$ tends to $T$. This suggests the presence of a \emph{phase transition} in the linear case: below this critical case, the noise is either propagated on a scale smaller than $\varepsilon^{3/2}$ (i.e. $\varepsilon^{\beta}$,  $\beta > 3/2$) or not propagated at all. We hope to be able to come back to this question later.

\subsection{Criticality of the Constant Coefficient Case}
In this subsection, we prove that the constant coefficient case is critical in the sense that the noise transmitted from the diffusion process $P$ to the absolutely continuous process $E$ is exponentially small as $t$ approaches $T$. We show that this transmission is so small that the first-order H\"ormander structure becomes unstable and we can find a perturbation of the constant coefficients case for which \eqref{eq:6:1:11} fails.

\begin{proposition}
\label{prop:7:1:1}
Assume that the coefficients $b$ and $\sigma$ are constant, ${\rm det}(\sigma)$ being non-zero, and that $f$ has the form
$-f(p,y)= \langle \alpha, p \rangle - \gamma y$, for some $\alpha \in \RR^d \setminus \{0\}$
and $\gamma >0$. Then, for any initial condition $(t_0,p,e) \in [0,T) \times \RR^d \times \RR$ and $t<T$, the pair $(P_t^{t_0,p},E_t^{t_0,p,e})$ has an infinitely differentiable density. 

Moreover, there exists $c' \geq 1$, depending on known parameters only,
such  that, for $T-t_0 \leq 1/c'$ and $0 \leq t-t_0 \leq (T-t_0)/2$,
\begin{equation}
\label{eq:14:2:18}
(c')^{-1} \exp \bigl( - \frac{c'}{(T-t_0)} \bigr) (t-t_0)^3
\leq {\text var} \bigl( E_t^{t_0,p,e} \bigr)
\leq c' \exp \bigl( - \frac{1}{c'(T-t_0)} \bigr) (t-t_0)^3,
\end{equation}
when $(t_0,\bar{e}) \in \{ (t,e') \in [0,T) \times \RR : 
(e'-\Lambda)/(T-t) \in [4\gamma/9,5\gamma/9]\}$, with 
$\bar{e} = e + w(t_0,p)$. 
\end{proposition}

Equation \eqref{eq:14:2:18} gives the typical size of the conditional fluctuations of
the process $E$ inside the critical cone: the power $3$ in $(t-t_0)^3$ is natural since $E$ behaves as the integral of a diffusion process, but the coefficient in front of $(t-t_0)^3$ is dramatically small when $t_0$ is close to time $T$.
In some sense, this says that the regime is nearly degenerate: inside the cone and close to time $T$, the trap formed by the characteristics of the first-order conservation law tames most of the randomness inherited from the diffusion process. In particular, the order of the variance at time $(T-t_0)/2$ is not $(T-t_0)^3$ but is exponentially small in $(T-t_0)$.
 
 \vspace{5pt}\noindent
 {\bf Proof.} The proof is divided in several steps.
First, we start with the following lemma
\begin{lemma}
\label{lem:7:3:1}
Under the assumption of Proposition \ref{prop:7:1:1}, the function 
$v$ defined in Proposition \ref{prop:7:11:10}
is infinitely differentiable on $[0,T) \times \RR^d \times \RR$. 
Moreover, for a sequence of smooth terminal conditions $(\phi^n)_{n \geq 1}$ converging towards the Heaviside terminal condition 
${\mathbf 1}_{[\Lambda,+\infty)}$, the sequences $(\partial_p v^{\phi_n})_{n \geq 1}$
and $(\partial_e v^{\phi_n})_{n \geq 1}$ converge towards $\partial_p v$ and
$\partial_e v$ respectively, uniformly on compact subsets of $[0,T) \times \RR^d \times \RR$. 
\end{lemma}

\vskip 1pt\noindent
{\bf Proof .}
Recalling the definitions \eqref{eq:7:1:3}  of $\bar{E}$ and \eqref{eq:19:12:1} of $w$, we see that
$w(t,p) = (T-t) \langle \alpha,(p+b(T-t)/2)\rangle$ and 
$\partial_p w(t,p) = (T-t) \alpha$, so that (compare with 
\eqref{eq:26:2:5})
\begin{equation}
\label{eq:30:7:1}
d \bar{E}_t = - \gamma Y_t dt + (T-t) \langle \sigma^{\top} \alpha,dW_t \rangle, \quad
t \in [t_0,T].
\end{equation}
Here, $(Y_t)_{t_0 \leq t \leq T}$ is a martingale with 
${\mathbf 1}_{(\Lambda,+\infty)}(\bar{E}_T) \leq Y_T \leq 
{\mathbf 1}_{[\Lambda,+\infty)}(\bar{E}_T)$ as terminal value. 
Therefore, the value function $v(t_0,p,e) = Y_{t_0}^{t_0,p,e}$ may be understood as the value function of the BSDE
$dY_t = \langle Z_t, dW_t \rangle$,
with ${\mathbf 1}_{(\Lambda,+\infty)}(\bar{E}_T) \leq Y_T \leq 
{\mathbf 1}_{[\Lambda,+\infty)}(\bar{E}_T)$ as terminal condition. 
That is, we expect $v(t_0,p,e)$ to coincide with $\bar{v}(t_0,\bar{e})$,  the solution $\bar{v}$  of
the PDE
\begin{equation}
\label{eq:7:3:3}
\partial_t \bar{v}(t,\bar{e}) + \frac{1}{2} (T-t)^2 
|\sigma^{\top} \alpha|^2 \partial^2_{\bar{e}\bar{e}} \bar{v}(t,\bar{e})
- \gamma \bar{v}(t,\bar{e}) \partial_{\bar{e}} \bar{v}(t,\bar{e}) = 0,
\quad t \in [0,T), \ \bar{e} \in \RR,
\end{equation}
with $\bar{v}(T,\bar{e}) = {\mathbf 1}_{[\Lambda,+\infty)}(\bar{e})$ as terminal condition, i.e. 
\begin{equation}
\label{eq:7:3:2}
v(t_0,p,e) = \bar{v}(t_0,\bar{e}) = \bar{v}(t_0,e+w(t_0,p))
= \bar{v}\bigl(t_0,e+  (T-t_0) \langle \alpha,(p+b(T-t_0)/2)\rangle\bigr).
\end{equation}
To prove \eqref{eq:7:3:2} rigorously, we follow the strategy of the first section: we first consider the case when the terminal condition is smooth, given by some non-decreasing $[0,1]$-valued smooth function $\phi$ approximating the Heaviside function; in this case, the PDE \eqref{eq:7:3:3} admits a continuous solution on the whole 
$[0,T] \times \RR$ which is ${\mathcal C}^{1,2}$ on 
$[0,T) \times \RR$ (the proof is \emph{mutatis mutandis} the proof given in  \cite{cdet}): we denote it by $\bar{v}^{\phi}$. By
It\^o's formula, $(\bar{v}^{\phi}(t,\bar{E}_t^{\phi,t_0,\bar{e}}))_{t_0 \leq t \leq T}$ is a martingale with 
$\phi(\bar{E}_T^{\phi,t_0,\bar{e}})=
\phi({E}_T^{\phi,t_0,\bar{e}})$ as terminal variable. (Here we use the same notation as in Corollary \ref{cor:25:7:1}.) Therefore, it coincides with 
$(Y_t^{\phi,t_0,p,e})_{t_0 \leq t \leq T}$ when $\bar{e} = e + w(t_0,p)$. This gives
the connection between $v^{\phi}$ and $\bar{v}^{\phi}$, i.e.
 \eqref{eq:7:3:2} for a smooth terminal condition. By Corollary \ref{cor:27:7:1}, $v^{\phi}$ converges towards $v$ when $\phi$ converges towards ${\mathbf 1}_{[\Lambda,+\infty)}$. Following the proof of Proposition 5 in \cite{cdet}, $\bar{v}^{\phi}$, $\partial_{\bar{e}} \bar{v}^{\phi}$ and $\partial_{\bar{e}\bar{e}}^2 \bar{v}^{\phi}$ converge towards $\bar{v}$, $\partial_{\bar{e}} \bar{v}$ and 
 $\partial_{\bar{e} \bar{e}}^2 \bar{v}$ respectively, $\bar{v}$ standing for the classical solution of \eqref{eq:7:3:3} on $[0,T) \times \RR$ satisfying $\bar{v}(t,\bar{e}) \rightarrow {\mathbf 1}_{[\Lambda,+\infty)}(\bar{e})$ as $t \nearrow T$ for $\bar{e} \not = \Lambda$. Passing to the limit along the regularization of the terminal condition, we complete the proof of \eqref{eq:7:3:2}. As a by-product, $\partial_e v^{\phi}$ and $\partial_p v^{\phi}$ converge towards
 $\partial_e v$ and $\partial_p v$ respectively, uniformly on compact subsets of $[0,T) \times \RR^d \times \RR$.

We emphasize that equation \eqref{eq:7:3:3} is uniformly elliptic on every 
$[0,T-\varepsilon] \times \RR$. Moreover, by Proposition \ref{prop:4:11:1},  $\partial_{\bar{e}} \bar{v}$ is uniformly bounded
on $[0,T-\varepsilon] \times \RR$ so that 
the product $\bar{v}(t,\bar{e}) \times \partial_{\bar{e}} \bar{v}(t,\bar{e})$ can be understood as $F(\bar{v}(t,\bar{e}),\partial_{\bar{e}} \bar{v}(t,\bar{e}))$ for a smooth function $F$ with compact support. By Section 3 in Crisan and Delarue in \cite{CrisanDelarue},
$\bar{v}$ is infinitely differentiable on $[0,T) \times \RR$. By \eqref{eq:7:3:2}, $v$ is infinitely differentiable on $[0,T) \times \RR^d \times \RR$.  
\qed

\vspace{5pt}\noindent
The second step of the proof of Proposition \ref{prop:7:1:1} provides two-sided bounds for $\partial_p v$:
\begin{lemma}
\label{lem:2:8:1}
Under the assumption of Proposition \ref{prop:7:1:1}, there exists a constant $C \geq 1$, depending on $L$ and $T$ only, such that, for any $1 \leq i \leq d$ satisfying $\alpha_i \not = 0$,
\begin{equation}
\label{eq:1:8:1}
\begin{split}
&C^{-1} {\mathbb P} \bigl\{ \inf_{(t_0+T)/2 \leq s
\leq T} |\bar{E}_s^{t_0,p,e} - \Lambda| > C (T-t_0) \bigr\}
 \\
&\hspace{15pt} \leq 1 - \gamma \alpha_i^{-1} \partial_{p_i} v(t_0,p,e) \leq 
{\mathbb E} \biggl[ \exp \biggl( - \int_{t_0}^T \gamma \partial_e v(s,P_s^{t_0,p},E_s^{t_0,p,e}) ds \biggr) \biggr].
\end{split}
\end{equation}
\end{lemma}

{\bf Proof.} When Eq. \eqref{eq:4:11:1} is driven by a smooth terminal condition 
$\phi$, $\partial_p v^{\phi}$ satisfies the system of PDEs (compare with \eqref{eq:13:03:1} and pay attention that $\partial_p v^{\phi}$ is a vector)
\begin{equation}
\label{eq:15:2:1}
\begin{split}
\partial_t \bigl[ \partial_p v^{\phi}  \bigr](t,p,e) &+ {\mathcal L}_p \bigl[ 
\partial_p v^{\phi} \bigr] (t,p,e) +
 \bigl[ \langle \alpha,p \rangle - \gamma 
v(t,p,e) \bigr] \partial_e \bigl[\partial_p v^{\phi}\bigr](t,p,e) 
\\
& + \bigl[ \alpha - \gamma 
\partial_p v^{\phi}(t,p,e) \bigr] \partial_e v^{\phi}(t,p,e) = 0,
\end{split}
\end{equation}
with $0$ as terminal condition. Setting $U_t = \partial_p v^{\phi}(t,P_t^{t_0,p},E_t^{t_0,p,e})$, we obtain
by It\^o's formula:
\begin{equation}
\label{eq:15:2:2}
dU_t = - \bigl[ \alpha - \gamma U_t \bigr] \partial_e v^{\phi}(t,P_t^{t_0,p},E_t^{\phi,t_0,p,e}) dt, 
\end{equation}
and the \emph{variation of the constant} formula gives
\begin{equation*}
\begin{split}
U_{t_0} &= \alpha  {\mathbb E} \biggl[ \int_{t_0}^T 
\partial_e v^{\phi}(t,P_t^{t_0,p},E_t^{\phi,t_0,p,e}) \exp \biggl( - \int_{t_0}^t \gamma \partial_e v^{\phi}(s,P_s^{t_0,p},E_s^{\phi,t_0,p,e}) ds \biggr) dt \biggr]
\\
&= \gamma^{-1} \alpha  \biggl\{ 1 - {\mathbb E} \biggl[
\exp \biggl( - \int_{t_0}^T \gamma \partial_e v^{\phi}(s,P_s^{t_0,p},E_s^{\phi,t_0,p,e}) ds \biggr) \biggr]
\biggr\}.
\end{split}
\end{equation*} 
As a consequence, we obtain
\begin{equation}
\label{eq:7:1:5}
\alpha - \gamma \partial_p v^{\phi}(t_0,p,e) 
= \alpha {\mathbb E} \biggl[
\exp \biggl( - \int_{t_0}^T \gamma \partial_e v^{\phi}(s,P_s^{t_0,p},E_s^{\phi,t_0,p,e}) ds \biggr) \biggr].
\end{equation} 
Now the point is to pass to the limit in \eqref{eq:7:1:5} along a mollification of the Heaviside terminal condition by using the same approximation procedure as in 
the proof of Lemma \ref{lem:7:3:1}. Clearly, the bound \eqref{eq:9:11:1} for $\partial_e v$ is not sufficient to apply Dominated Convergence Theorem. At least, (two-sided) Fatou's Lemma yields the upper bound in \eqref{eq:1:8:1}.

To get the lower bound, we apply Proposition \ref{prop:20:07:5} below. Approximating 
the Heaviside function by a non-decreasing sequence of non-decreasing $[0,1]$-valued smooth functions $(\phi^n)_{n \geq 1}$ satisfying $\phi^n(e)=1$ for 
$e \geq \Lambda$, we deduce from \eqref{eq:1:8:2} that $v^n=v^{\phi^n}$ satisfies
$\partial_e v^{n}(t,p,e) \leq C(T-t_0)^2$ for $\bar{e} - \Lambda > C(T-t_0)$,
$t_0 \leq t \leq T$, $p \in \RR^d$, for some constant $C \geq 1$.
Applying \eqref{eq:7:1:5} and Proposition \ref{prop:4:11:1} and modifying $C$ if necessary, 
\begin{equation}
\label{eq:4:8:10}
1 - \gamma \alpha_i^{-1} \partial_{p_i} v^{n}(t_0,p,e)
\geq C^{-1} {\mathbb P} \bigl\{ \inf_{(t_0+T)/2 \leq s
\leq T} \bigl[ \bar{E}_s^{\phi^n,t_0,p,e} - \Lambda \bigr] > C (T-t_0) \}.
\end{equation}
Following the proof of Proposition \ref{prop:7:11:10}, 
$E^{\phi^n,t_0,p,e}_t \rightarrow E_t^{t_0,p,e}$ uniformly in time $t$ in compact subsets of $[t_0,T)$, a.s. as $n \rightarrow + \infty$. Actually, convergence is a.s. uniform on the whole $[t_0,T]$, since 
$|dE^{\phi^n,t_0,p,e}_t/dt| \leq |\alpha| \sup_{t_0 \leq t  \leq T} |P_t^{t_0,p}| + \gamma$.
As $\bar{E}_t^{\phi^n,t_0,p,e}
= E_t^{\phi^n,t_0,p,e} + w(t,P_t^{t_0,p})$, we deduce that, a.s., $(\bar{E}_t^{\phi^n,t_0,p,e})_{t_0 \leq t \leq T}$ converges towards $(\bar{E}_t^{t_0,p,e})_{t_0 \leq t \leq T}$ uniformly on $[t_0,T]$.
 Therefore, we can pass to the limit in the above inequality. 
We obtain \eqref{eq:1:8:1} but without the absolute value in the infimum.
Choosing the approximating sequence $(\phi^n)_{n \geq 1}$  such that $\phi^n(e)=0$ for $e \leq \Lambda$ and repeating the argument, we obtain 
\eqref{eq:4:8:10} with $\bar{E}_s^{\phi^n,t_0,p,e} - \Lambda$ replaced by 
$\Lambda - \bar{E}_s^{\phi^n,t_0,p,e}$. Passing to the limit, we complete the proof.
\qed
\vspace{5pt}

\noindent
The third step of the proof of Proposition \ref{prop:7:1:1} provides a sharp estimate of $\partial_p v(t_0,p,e)$ for $(t,p,e)$ in the critical cone in terms of 
the distance to the terminal time $T$. 

\begin{proposition}
\label{prop:14:2:1}
There exists a constant $c \geq 1$, depending on known parameters only, such that, 
for $T-t_0 \leq 1/c$ and $\bar{e} = e + w(t_0,p) \in [\Lambda +
3(T-t_0) \gamma/8,\Lambda + 5(T-t_0) \gamma/8]$,
\begin{equation*}
c^{-1} \exp\bigl( - \frac{c}{T-t_0} \bigr)
\leq 
\alpha_i^{-1} \bigl[ \alpha - \gamma \partial_p v(t_0,p,e) \bigr]_i
\leq c \exp\bigl( - \frac{1}{c(T-t_0)} \bigr), \quad
1 \leq i \leq d : \alpha_i \not = 0.
\end{equation*}
Actually, the lower bound (i.e. the left-hand side inequality)
holds for any $(t_0,p,e) \in [0,T) \times \RR^d \times \RR$.
\end{proposition}

{\bf Proof.} 
Once more, we use $\bar{E}$ defined in \eqref{eq:30:7:1}. The initial condition of $\bar{E}$ is denoted by $(t_0,\bar{e})$, that is 
$\bar{E}_{t_0} = \bar{e}$. It satisfies $\bar{e}=
e+w(t_0,p)$, with $E_{t_0}=e$. 
\vskip 2pt

\noindent
\emph{Lower Bound.} By Lemma \ref{lem:2:8:1}, it is sufficient to bound the probability
${\mathbb P}\{ \inf_{(t_0+T)/2 \leq t \leq T} |\bar{E}_t - \Lambda| > C(T-t_0)\}$ from below.
Clearly, 
\begin{equation}
\label{eq:14:2:5}
\begin{split}
&{\mathbb P} \bigl\{ \inf_{(t_0+T)/2 \leq t < T} |\bar{E}_t - \Lambda| >  C(T-t_0) \bigr\}
\\
&\geq {\mathbb P} \bigl\{ |\bar{E}_{(t_0+T)/2} - \Lambda| \geq (C+\frac{\gamma}{2}+1)(T-t_0) \bigr\}
\\
&\hspace{5pt} \times
\sup_{p' \in \RR^d,\bar{e}' \in \RR} 
 {\mathbb P} \bigl\{ \sup_{(t_0+T)/2 \leq t \leq T} \bigl| \bar{E}_t - \bar{E}_{(t_0+T)/2} \bigr| \leq (\frac{\gamma}{2}+1)(T-t_0) \big| 
(P,\bar{E})_{\frac{T+t_0}{2}} = (p',\bar{e}')\bigr\}
\\
&= \pi_1 \times \pi_2.
\end{split}
\end{equation}
By \eqref{eq:30:7:1} and by maximal inequality (IV.37.12) in Rogers and Williams \cite{RogersWilliams.book}, the conditional probability $\pi_2$ can be easily estimated:
\begin{equation}
\label{eq:14:2:6}
\pi_2
\geq {\mathbb P} \biggl\{ \sup_{(t_0+T)/2 \leq t \leq T} \biggl| \int_{(t_0+T)/2}^t 
(T-s) \langle \sigma^{\top} \alpha ,dW_s \rangle  \biggr| \leq (T-t_0)
 \biggr\}
\geq 1 - \exp \bigl( - \frac{1}{c(T-t_0)} \bigr),
\end{equation}
for some constant $c \geq 1$ depending on known parameters only.

Turning to $\pi_1$ and using standard Gaussian bounds, if $\bar{e} \geq \Lambda$,
\begin{equation}
\label{eq:09:03:1}
\begin{split}
&{\mathbb P}\bigl\{ \bar{E}_{t_0+T/2} \geq \Lambda + (C+\frac{\gamma}{2}+1)(T-t_0) \bigr\}
\\
&\geq 
{\mathbb P}\biggl\{ - \gamma \frac{T-t_0}{2} + \int_{t_0}^{(t_0+T)/2}
(T-s) \langle \sigma^{\top} \alpha, dW_s \rangle
 \geq  (C+\frac{\gamma}{2}+1)(T-t_0) \biggr\}
\\
&\geq {\mathbb P} \biggl\{ \int_{t_0}^{(t_0+T)/2} 
(T-s)
\langle \sigma^{\top} \alpha, dW_s \rangle \geq (C+ \gamma+1)(T-t_0) \biggr\}
\geq \exp \bigl( - \frac{1}{c'(T-t_0)} \bigr),
\end{split}
\end{equation}
for some constant $c' \geq 1$ depending on known parameters only.
The case when $\bar{e} < \Lambda$ can be handled in a similar way.
Using the left-hand side in \eqref{eq:1:8:1}, we complete the proof of the lower bound.

\vskip 2pt\noindent
\emph{Upper Bound.} Consider again $\bar{E}$ as in \reff{eq:30:7:1} with $\bar{E}_{t_0}= \bar{e}$. Then
\begin{equation*}
d \bigl[ \gamma (T-t) Y_t - \bar{E}_t \bigr]
= \gamma(T-t) \langle Z_t, dW_t \rangle - (T-t) \langle \sigma^{\top} \alpha, dW_t \rangle,
\quad t_0 \leq t < T,
\end{equation*}
i.e. $(\gamma (T-t) Y_t - \bar{E}_t)_{t_0 \leq t < T}$ is a martingale. In particular,
\begin{equation*}
\gamma (T-t_0) v(t_0,p,e) = \gamma (T-t_0) v\bigl(t_0,p,\bar{e}-w(t_0,p)\bigr) = \bar{e} - {\mathbb E} \bigl[ \bar{E}_T^{t_0,p,\bar{e}} \bigr]. 
\end{equation*}
We then aim at bounding from below the derivative of $v$ with respect to $e$ (keep in mind that $v$ is differentiable before $T$). 
By \eqref{eq:2:8:2}, the Radon-Nykodym
of the non-decreasing Lipschitz function $\bar{e} \hookrightarrow \bar{E}_T^{t_0,p,\bar{e}}$ is a.s. less than 1. For $(\bar{e}-\Lambda)/(T-t_0) \in (\gamma/4,(3\gamma)/4)$, the bound 
$\bar{e} \hookrightarrow \bar{E}_T^{t_0,p,\bar{e}}$ cannot be achieved with probability 1. Indeed, by Proposition \ref{prop:30:12:1} (with $\ell_1=\ell_2=\gamma$), for $(\bar{e}-\Lambda)/(T-t_0) \in (\gamma/4,(3\gamma)/4)$
and for $T-t_0$ small enough (independently of $p$ and $\bar{e}$), we have
$\bar{E}_T^{t_0,p,\bar{e}} = \Lambda$ on the set
$F = \{ \sup_{t_0 \leq t \leq T} | \int_{t_0}^t
\langle \sigma^{\top} \alpha, dB_s \rangle | \leq \gamma/16 \}$ (see also \eqref{eq:08:03:2}). In particular, for $(\bar{e}-\Lambda)/(T-t_0) \in (\gamma/4,(3 \gamma)/4)$, the Radon-Nykodym derivative $\partial_{\bar{e}} \bar{E}_T^{t_0,p,\bar{e}}$ is zero on the set $F$. It follows that
\begin{equation*}
\gamma (T-t_0) \partial_{e} v(t_0,p,e) \geq  1 - \PP\bigl(F^{\complement}\bigr) = \PP(F), \quad 
\Lambda + \gamma/4 \leq \bar{e} \leq  \Lambda+ 3\gamma/4. 
\end{equation*}
This proves that
\begin{equation*}
\partial_{e} v(t_0,p,e) \geq \frac{ \PP(F)}{\gamma(T-t_0)}, \quad 
\Lambda + \gamma/4 \leq \bar{e} \leq  \Lambda+ 3\gamma/4.
\end{equation*}
Using maximal inequality once more, there exists a constant $c \geq 1$ depending on known parameters only such that $\PP(F) \geq 1 - \exp ( -  1/[c(T-t_0])$.
We deduce that, for $\bar{e} \in [\Lambda + (\gamma/4)(T-t_0),\Lambda + (3\gamma/4)(T-t_0)]$,
\begin{equation}
\label{eq:7:1:4}
\partial_e v(t_0,p,e) \geq \frac{1 - \exp \bigl[ - [c (T-t_0)]^{-1} \bigr]}{\gamma(T-t_0)}.
\end{equation}
Equation \reff{eq:7:1:4} is crucial: it says that the gradient of $v$ with respect to $e$
is not integrable inside the critical cone. We can plug it into \reff{eq:1:8:1}. We obtain that, for any 
$i \in \{1,\dots,d\}$ such that $\alpha_i \not = 0$,
\begin{equation}
\label{eq:7:1:10}
1 - \gamma \alpha_i^{-1} \partial_{p_i} v(t_0,p,e) \leq \EE
\biggl[ \exp \biggl( - \int_{t_0}^T \gamma \partial_e v(s,P_s,E_s) ds \biggr)
{\mathbf 1}_{\{(t,\bar{E}_t)_{t_0 \leq t < T} \in {\mathcal C}^{\complement}\}}\biggr],
\end{equation}
with $(P_{t_0},E_{t_0})=(p,e)$ and where ${\mathcal C} = \{ (t,e) \in [0,T) \times \RR : 
(e-\Lambda)/(T-t) \in [\gamma/4,3\gamma/4]\}$. Indeed, inside
the cone, the integral inside the exponential explodes so that the exponential vanishes. 

We deduce that, for any coordinate $i \in \{1,\dots,d\}$ for which $\alpha_i \not= 0$,
\begin{equation*}
\alpha_i^{-1} \bigl[\alpha - \gamma \partial_p v(t_0,p,e) \bigr]_i \leq 
\PP \bigl\{ (t,\bar{E}_t)_{t_0 \leq t < T} \in {\mathcal C}^{\complement} \bigr\}.
\end{equation*}
Choose now $\bar{e}$ deep in the middle of the critical cone, say
$\bar{e}$ in the interval $[\Lambda+(3\gamma/8)(T-t_0),$ \ $\Lambda + (5 \gamma/8)(T-t_0)]$. Then, 
\begin{equation*}
\begin{split}
\PP \bigl\{ (t,\bar{E}_t)_{t_0 \leq t < T} \in {\mathcal C}^{\complement} \bigr\}
\leq 
\PP \bigl\{ \exists t \in [t_0,T) : |\bar{E}_t - \Lambda - \frac{\gamma}{2} (T-t) | \geq \frac{\gamma}{4} (T-t) \bigr\}.
\end{split}
\end{equation*}
Using the same notation as in the proof of Proposition \ref{prop:30:12:1}, we claim
\begin{equation}
\label{eq:08:03:4}
\begin{split}
&\PP \bigl\{ \exists t \in [t_0,T) : |\bar{E}_t - \Lambda - \frac{\gamma}{2} (T-t) | \geq \frac{\gamma}{4} (T-t) \bigr\}
\\
&\leq \PP \bigl\{ \exists t \in [t_0,T) : 
\bar{E}_t^+ \geq \Lambda + \frac{3\gamma}{4}(T-t) \bigr\}
+ \PP \bigl\{ \exists t \in [t_0,T) : 
\bar{E}_t^- \leq \Lambda + \frac{\gamma}{4} (T-t) \bigr\}
\\
&= \pi_3 + \pi_4.
\end{split}
\end{equation}
We now use $\bar{Z}^{\pm}$. Until
$\bar{E}^{+}$ reaches $\Lambda + (3\gamma/4)(T-t)$, 
it holds $\bar{E}^+ \leq \bar{Z}^+$. Indeed, the drift of the process $\bar{Z}^+$
is then greater than the drift of the process $\bar{E}$. Therefore, by \eqref{eq:14:2:10} and by the relationship $\partial_p w(s,\cdot)=(T-s) \alpha$,
\begin{equation*}
\begin{split}
\pi_3 &\leq \PP \bigl\{ \exists t \in [t_0,T) : 
\bar{Z}_t^+ \geq \Lambda + \frac{3\gamma}{4}(T-t) \bigr\}
\\
&= \PP \bigl\{ \exists t \in[t_0,T) : 
C' \int_{t_0}^t (T-s)^{\beta-1} ds + \langle \sigma^{\top} \alpha, W_t- W_{t_0} \rangle \geq \frac{3 \gamma}{4} - 
\frac{\bar{e}-\Lambda}{T-t_0} \bigr\}.
\end{split}
\end{equation*}
Since $(\bar{e}-\Lambda)/(T-t_0) \leq 5 \gamma/8$, we get 
\begin{equation*}
\pi_3
\leq \PP \bigl\{ \exists t \in[t_0,T) : 
C' \int_{t_0}^t (T-s)^{\beta-1} ds + 
\langle \sigma^{\top} \alpha, W_t- W_{t_0} \rangle\geq \frac{\gamma}{8} \bigr\}.
\end{equation*}
We deduce that there exists a constant $c' \geq 1$, depending on known parameters only, such that, for $T-t_0 \leq 1/c'$, 
\begin{equation}
\label{eq:3:8:1}
\PP \bigl\{ \exists t \in [t_0,T) : 
\bar{E}_t^+ \geq \Lambda + \frac{3\gamma}{4}(T-t) \bigr\}
\leq \exp \bigl( - \frac{1}{c'(T-t_0)} \bigr).
\end{equation}
Handling the second term in \eqref{eq:08:03:4} in a similar way, we get, for $(t_0,\bar{e}) \in {\mathcal C}$,
\begin{equation*}
\alpha_i^{-1} \bigl[ \alpha_i  - \gamma [ \partial v/\partial p_i] (t_0,p,e) \bigr]_i
\leq 2 \exp \bigl( - \frac{1}{c'(T-t_0)} \bigr), \quad
i =1, \dots,d : \alpha_i \not = 0.
\end{equation*}
This completes the proof. 
\qed
\vspace{5pt}

We are now ready to complete the proof of Proposition \ref{prop:7:1:1}. 
The existence of an infinitely differentiable density for the pair 
$(P_t,E_t)$, $t<T$, follows directly from H\"ormander Theorem applied to 
the parabolic adjoint operator $\partial_t - ({\mathcal L}_p - f(p,v(t,p,e)) \partial_e)^*$, with ${\mathcal L}_p$ as in \eqref{eq:3:8:2}, see Theorem 1.1 in H\"ormander \cite{Hormander} with $X_0 = \partial_t + b \partial_p - f(p,v(t,p,e)) \partial_e$
therein. See also Delarue and Menozzi \cite{DelarueMenozzi2} for a specific application of H\"ormander Theorem to the current setting. 

To estimate the conditional variance of $E_t$, we compute the Malliavin derivative of $E_t$, $t_0 \leq t < T$. (Below, the initial condition is $(P_{t_0},E_{t_0})=(p,e)$.)
For any coordinate $i \in \{1,\dots,d\}$ and any $t_0 \leq s \leq t <T$,
\begin{equation*}
dD_s^i E_t = \bigl[ 
\langle \alpha - \gamma \partial_p v(t,P_t,E_t),D_s^i P_t \rangle - 
\gamma \partial_e v(t,P_t,E_t) D_s^i E_t \bigr] dt,
\end{equation*}
with $D_s^i E_s =0$, by Theorem 2.2.1 in Nualart \cite{Nualart}. We deduce
\begin{equation}
\label{eq:18:2:10}
\begin{split}
D_s^i E_t &= \int_s^t  \langle \alpha - \gamma \partial_p v(r,P_r,E_r),D_s^i P_r \rangle 
\exp \biggl( - \int_r^t \gamma \partial_e v(u,P_u,E_u) du \biggr) dr
\\
&= \int_s^t  \langle \alpha - \partial_p v(r,P_r,E_r),\sigma_{\cdot,i} \rangle 
\exp \biggl( - \int_r^t \gamma \partial_e v(u,P_u,E_u) du \biggr) dr,
\end{split}
\end{equation}
that is
\begin{equation*}
D_s E_t = \sigma^{\top} \int_s^t \bigl[\alpha - \gamma \partial_p v(r,P_r,E_r)\bigr] \exp \biggl( - \int_r^t \gamma \partial_e v(u,P_u,E_u) du \biggr) dr.
\end{equation*}
By Proposition \ref{prop:14:2:1}, there exists a bounded $\RR^d$-valued process 
$(\theta_t)_{t_0 \leq t < T}$, with positive coordinates, such that
\begin{equation*}
\bigl[ (\sigma^{\top})^{-1} D_s E_t \bigr]_i  =   \alpha_i \int_s^t (\theta_r)_i  \exp \biggl( - \int_r^t \gamma \partial_e v(u,P_u,E_u) du \biggr) dr. 
\end{equation*}
(By Bouleau and Hirsch criterion (see Theorem 2.1.3 in \cite{Nualart}), we recover that 
the distribution of $E_t$ is absolutely continuous with respect to the Lebesgue measure. Obviously, this is already known by H\"ormander Theorem.) Moreover, we can write
\begin{equation*}
{\mathbb E} \bigl[ \bigl( \sigma^{\top} \bigr)^{-1} D_s E_t |{\mathcal F}_s \bigr]
=   \int_s^t {\mathbb E} \biggl[ 
\alpha \star \theta_r  \exp \biggl( - \int_r^t \gamma \partial_e v(u,P_u,E_u) du \biggr) \bigl|{\mathcal F}_s 
\biggr] dr,
\end{equation*}
where $\alpha \star \theta_r$ is understood as $(\alpha_i (\theta_r)_i)_{1 \leq i \leq d}$. 
We then follow the proof of the lower bound in Proposition \ref{prop:14:2:1}.
Setting $\theta^{*} = \inf_{1 \leq i \leq d} \inf_{s \leq r \leq t}
[(\theta_r)_i]$ and using the lower bound explicitly, we deduce
\begin{equation}
\label{eq:15:2:10}
\begin{split}
\bigl| {\mathbb E} \bigl[ \bigl( \sigma^{\top} \bigr)^{-1} D_s E_t |{\mathcal F}_s \bigr] \bigr|
&\geq  |\alpha| {\mathbb E} \biggl[ 
\theta^*  \int_s^t 
\exp \biggl( - \int_r^t \gamma \partial_e v(u,P_u,E_u) du \biggr) dr\bigl|{\mathcal F}_s 
\biggr]
\\
&\geq \frac{|\alpha|}{c} \exp\bigl( - \frac{c}{T-t} \bigr)
{\mathbb E} \biggl[  \int_s^t 
\exp \biggl( - \int_r^t \gamma \partial_e v(u,P_u,E_u) du \biggr) dr \bigl|{\mathcal F}_s 
\biggr],
\end{split}
\end{equation}
Since $\partial_e v(u,\cdot,\cdot) \in [0,\gamma^{-1} (T-u)^{-1}]$, we can modify $c$ so that
\begin{equation}
\label{eq:3:8:4}
\bigl| {\mathbb E} \bigl[ \bigl( \sigma^{\top} \bigr)^{-1} D_s E_t |{\mathcal F}_s \bigr] \bigr|^2
\geq c^{-1} |\alpha|^2 \exp\bigl( - c/(T-t) \bigr) (t-s)^2,
\end{equation}
for $t_0 \leq s < t \leq (T+t_0)/2$.
Since $\sigma$ is invertible, we deduce that
\begin{equation*}
\bigl| {\mathbb E} \bigl[ D_s E_t |{\mathcal F}_s \bigr] \bigr|^2
\geq c^{-1} |\alpha|^2 \exp\bigl( - c/(T-t) \bigr) (t-s)^2.
\end{equation*}
Therefore,
\begin{equation}
\label{eq:14:2:15}
{\mathbb E} \int_{t_0}^t \bigl| {\mathbb E} \bigl[ D_s E_t |{\mathcal F}_s \bigr] \bigr|^2 ds
\geq c^{-1} |\alpha|^2 \exp\bigl( - c/(T-t) \bigr) (t-t_0)^3, \quad 0 \leq t-t_0 \leq \frac{T-t_0}{2}.
\end{equation}

We now seek an upper bound
for $| {\mathbb E} [ D_s E_t |{\mathcal F}_s ] |^2$ 
when $(t_0,\bar{e})$ belongs to the critical cone
${\mathcal C}''= \{ (t,e') \in [0,T) \times \RR : 
(e'-\Lambda)/(T-t) \in [4\gamma/9,5\gamma/9]\}$. Assume that $(s,\bar{E}_s) \in {\mathcal C}' =  \{ (t,e') \in [0,T) \times \RR : 
(e'-\Lambda)/(T-t) \in [3\gamma/8,5\gamma/8]\}$. Then, following the proof of \eqref{eq:08:03:4}--\eqref{eq:3:8:1}, we can find $c_2'\geq 1$ such that, for $T-s \leq 1/c_2'$,
\begin{equation}
\label{eq:14:2:20}
{\mathbb P}\bigl\{ \exists r \in [s,t] : (r,\bar{E}_r) \in {\mathcal C}^{\complement} |{\mathcal F}_s \bigr\} \leq c_2'  \exp \bigl( - \frac{1}{c_2'(T-s)} \bigr),
\end{equation}
with ${\mathcal C}$ as in \eqref{eq:7:1:10}. 
Using the upper bound in Proposition \ref{prop:14:2:1} and assuming
$T-s \leq 1/c_2'$ (modifying $c_2'$ if necessary), we deduce
(by reversing the inequality in \eqref{eq:15:2:10})
\begin{equation*}
\bigl| {\mathbb E} \bigl[ \bigl( \sigma^{*} \bigr)^{-1} D_s E_t |{\mathcal F}_s \bigr] \bigr|^2
\leq  c_2' \exp \bigl( - \frac{1}{c_2'(T-s)} \bigr) (t-s)^2 {\mathbf 1}_{{\mathcal C}'}(s,\bar{E}_s)
+ c_2' (t-s)^2{\mathbf 1}_{{{\mathcal C}'}^{\complement}}(s,\bar{E}_s).
\end{equation*}
Taking the expectation and applying a similar bound to \eqref{eq:14:2:20}, we deduce
\begin{equation*}
{\mathbb E} \bigl[ \bigl| {\mathbb E} \bigl[ \bigl( \sigma^{*} \bigr)^{-1} D_s E_t |{\mathcal F}_s \bigr] \bigr|^2 \bigr]
\leq  c_2' \exp \bigl( - \frac{1}{c_2'(T-s)} \bigr) (t-s)^2 
+ c_2' 
\exp \bigl( - \frac{1}{c_2'(T-t_0)} \bigr) (t-s)^2,
\end{equation*}
whenever $(t_0,\bar{e}) \in {\mathcal C}''= \{ (t,e) \in [0,T) \times \RR : 
(e-\Lambda)/(T-t) \in [4\gamma/9,5\gamma/9]\}$ and $T-t_0 \leq 1/c_2'$. 
In such a case,
\begin{equation}
\label{eq:14:2:16}
\int_{t_0}^t {\mathbb E} \bigl[ \bigl| {\mathbb E} \bigl[ \bigl( \sigma^{*} \bigr)^{-1} D_s E_t |{\mathcal F}_s \bigr] \bigr|^2 \bigr] ds
\leq   c_2' 
\exp \bigl( - \frac{1}{c_2'(T-t_0)} \bigr) (t-t_0)^3.
\end{equation}
By Clark-Ocone formula, we deduce from \eqref{eq:14:2:15} and \eqref{eq:14:2:16}
that
\begin{equation*}
(c_2')^{-1} \exp \bigl( - \frac{c_2'}{(T-t)} \bigr) (t-t_0)^3
\leq \text{ var} \bigl( E_t^{t_0,p,e} \bigr)
\leq c_2' \exp \bigl( - \frac{1}{c_2'(T-t_0)} \bigr) (t-t_0)^3,
\end{equation*}
when $t_0 \leq t \leq (T+t_0)/2$ and $(t_0,\bar{e}) \in {\mathcal C}''$. 
\qed

\subsection{Enlightening Example}
We provide an example showing that the structure of the diffusion process
$(P_t)_{t \geq 0}$ can affect the validity of condition \eqref{eq:6:1:11} in the sense that \eqref{eq:6:1:11} can fail and
\begin{equation}
\label{eq:08:03:11}
\partial_p \bigl[ f\bigl(p,v(t,p,e) \bigr) \bigr] = 0,
\end{equation}
at some point. This doesn't mean that the hypoelliptic property fails since the noise may be transmitted from the first to the second equation through higher order derivatives. However, this suggests that a transition exists in the regime
of the pair process $(P,E)$. Actually, we believe that some noise is indeed transmitted from the first to the second equation since the solution $v$ to the PDE is shown to be smooth inside $[0,T) \times \RR \times \RR$ despite the degeneracy property 
\eqref{eq:08:03:11}. For this reason, the example is quite striking since there is some degeneracy, but some smoothing as well.

The idea is to go back to the setting of Proposition \ref{prop:7:1:1} when $d=1$ and to assume therein that the drift has the form $b(p) = b + \lambda p$ for some
$\lambda \in \RR$. As noticed above, condition \eqref{eq:6:1:11} is fulfilled when 
$\lambda=0$, but the noise transmitted to the process $E$ is exponentially small with respect to the distance to the singularity.

\begin{proposition}
\label{prop:08:03:1}
Assume that $d=1$, $b(p) = b + \lambda p$, $p \in \RR$ for some $\lambda \in \RR$,  
$\sigma$ is equal to a strictly positive constant, and that 
$-f(p,y) = \alpha p - \gamma y$, $p,y \in \RR$, for some $\alpha,\gamma >0$.
Then, $v$ is infinitely differentiable on $[0,T) \times \RR \times \RR$.

Moreover, if $\lambda <0$, then for any starting point $(t_0,p,e) \in [0,T) \times \RR \times \RR$ the pair $(P_t^{t_0,p},E_t^{t_0,p,e})$ has an infinitely differentiable density at all time $t \in [t_0,T)$. And, there exists $c' \geq 1$, depending on known parameters only, such that, for $T-t_0 \leq 1/c'$ and
$0 \leq t-t_0 \leq (T-t_0)/2$, 
\begin{equation}
\label{eq:08:03:12}
(c')^{-1}(T-t_0)^2 (t-t_0)^3 \leq \text{ var }(E_t^{t_0,p,e}) \leq
c' (T-t_0)^2 (t-t_0)^3,
\end{equation}
when $(t_0,\bar{e}) \in {\mathcal C}''
=\{(t,e') \in [0,T) \times \RR : (e'-\Lambda)/(T-t) 
\in [4 \gamma/9,5 \gamma/9]\}$, $\bar{e} = e + w(t_0,p)$. 

Finally, if $\lambda >0$, then for any $(t_0,p) \in [0,T) \times \RR$ such that 
$T-t_0 \leq 1/c'$, for some constant $c'$ depending on known parameters only,
there exists $e \in \RR$ such that $\alpha - \partial_p v(t_0,p,e)=0$. 
\end{proposition}

The reader should compare \eqref{eq:08:03:12} with \eqref{eq:14:2:18}: clearly,
the value $\lambda=0$ appears as a critical threshold. Notice that 
$c'$ depends on $\lambda$ in \eqref{eq:08:03:12}.
We also notice that the coefficient $b$ doesn't satisfy Assumption 
(A.4) since it is not bounded. Anyhow, we know that the Dirac mass exists by Remark \ref{rem:4:9:1}. Actually, we know a little bit more: as emphasized in Remark 
\ref{rem:4:9:1}, the assumption $\|b\|_{\infty} < + \infty$  in 
Proposition \ref{prop:30:12:1} is required to bound the increments of the process $P$
in the proof of Proposition \ref{prop:26:2:1}; in the current framework, 
Proposition \ref{prop:26:2:1} is useless since Lemma \ref{lem:22:12:1} applies directly. (We let the reader check that only the boundedness of $\sigma$ is used in the proof of Lemma
\ref{lem:22:12:1}.) Therefore, the original version of Proposition \ref{prop:30:12:1} is still valid in the current framework.

\vspace{5pt}\noindent
{\bf Proof.}  We first prove that $v$ is infinitely differentiable 
on $[0,T) \times \RR \times \RR$. To do so, we compute
\begin{equation*}
{\mathbb E}[P_t^{t_0,p}] = p + b(t-t_0) + \lambda 
\int_{t_0}^t {\mathbb E}[P_s^{t_0,p}] ds,
\end{equation*}
that is
\begin{equation*}
{\mathbb E}[P_t^{t_0,p}] = \exp[\lambda(t-t_0)] \bigl[ p 
+ \frac{b}{\lambda} \bigr]  - \frac{b}{\lambda}.
\end{equation*}
Therefore,
\begin{equation*}
\begin{split}
w(t_0,p) 
= \alpha \int_{t_0}^T {\mathbb E}[P_s^{t_0,p}] ds
= \alpha \int_{t_0}^T \exp[\lambda(s-t_0)] \bigl[ p 
+ \frac{b}{\lambda} \bigr] ds - \frac{\alpha b}{\lambda}(T-t_0).
\end{split}
\end{equation*}
In particular,
\begin{equation*}
\partial_p w(t_0,p) = \alpha \int_{t_0}^T \exp[\lambda(s-t_0)] ds
= \frac{\alpha}{\lambda} \bigl[ \exp[\lambda(T-t_0)] - 1 \bigr].
\end{equation*}
This shows that $\bar{E}$ has autonomous dynamics, i.e.
\begin{equation*}
d \bar{E}_t = - \gamma Y_t dt + 
\frac{\alpha}{\lambda} \bigl[ \exp[\lambda(T-t)] - 1 \bigr] dW_t.
\end{equation*}
The infinite differentiability of $v$ is then proven as in Lemma \ref{lem:7:3:1}.
\vskip 2pt
Next, we use \eqref{eq:13:03:1} and we follow \eqref{eq:15:2:1} and \eqref{eq:15:2:2}. Again with 
$U_t = \partial_p v(t,P_t,E_t)$, we write
\begin{equation*}
dU_t = - \bigl[ \alpha - \gamma U_t \bigr] \partial_e v(t,P_t,E_t) dt
- \lambda U_t dt, \quad 0 \leq t < T.
\end{equation*}
The additional $\lambda$ here comes from $\partial_p b(p) = \lambda$ in
\eqref{eq:13:03:1}. The above expression also holds for a smooth terminal condition $\phi$. In such a case, $\partial_p v^{\phi}(T,\cdot,\cdot)=0$.
By variation of the constant, we obtain
\begin{equation*}
\begin{split}
\partial_p v^{\phi}(t_0,p,e) =  \alpha 
{\mathbb E} \biggl[ \int_{t_0}^T \partial_e 
v^{\phi}(t,P_t^{t_0,p},E_t^{\phi,t_0,p,e}) \exp \biggl( \int_{t_0}^t
\bigl[ \lambda - \gamma \partial_e v^{\phi}(s,P_s^{t_0,p},E_s^{\phi,t_0,p,e}) \bigr]
ds \biggr) dt \biggr],
\end{split}
\end{equation*}
when the terminal condition $\phi$ is smooth.
By integration by parts, we deduce
\begin{equation*}
\begin{split}
\partial_p v^{\phi}(t_0,p,e) &= \gamma^{-1} \alpha \biggl[1 - \exp \bigl( \lambda(T-t_0) \bigr) 
{\mathbb E} \biggl[ 
\exp \biggl( - \gamma \int_{t_0}^T \partial_e v^{\phi}(s,P_s^{t_0,p},E_s^{\phi,t_0,p,e}) ds \biggr)
\biggr] \biggr] 
\\
&\hspace{15pt}
+ \lambda \gamma^{-1} \alpha 
{\mathbb E} \biggl[ \int_{t_0}^T  
\exp \biggl( \int_{t_0}^t
\bigl[ \lambda - \gamma \partial_e v^{\phi}(s,P_s^{t_0,p},E_s^{\phi,t_0,p,e}) \bigr]
ds \biggr) dt \biggr].
\end{split}
\end{equation*}
Finally, for a smooth boundary condition $\phi$, we obtain the analogue of \eqref{eq:7:1:5}:
\begin{equation}
\label{eq:15:2:55}
\begin{split}
\alpha - \gamma \partial_p v^{\phi}(t_0,p,e) &= \alpha \exp \bigl( \lambda(T-t_0) \bigr) 
{\mathbb E} \biggl[ 
\exp \biggl( - \gamma \int_{t_0}^T \partial_e v^{\phi}(s,P_s^{t_0,p},E_s^{\phi,t_0,p,e}) ds \biggr)
\biggr] 
\\
&\hspace{15pt}
- \lambda  \alpha 
{\mathbb E} \biggl[ \int_{t_0}^T  
\exp \biggl( \int_{t_0}^t
\bigl[ \lambda - \gamma \partial_e v^{\phi}(s,P_s^{t_0,p},E_s^{\phi,t_0,p,e}) \bigr]
ds \biggr) dt \biggr].
\end{split}
\end{equation}
As in the proof of the lower bound in Lemma \ref{lem:2:8:1}, we aim
at passing to the limit in the above expression along a mollification of the Heaviside terminal condition. Applying Dominated Convergence Theorem to the second term in the right-hand side and handling the first term as in the proof of Lemma \ref{lem:2:8:1},
we can 
find a constant $c \geq 1$ depending on known parameters only (and possibly on $\lambda$) such that (below, the initial condition is $(P_{t_0},E_{t_0})=(p,e)$)
\begin{equation}
\label{eq:15:2:50}
\begin{split}
\alpha - \gamma \partial_p v(t_0,p,e) &\geq c^{-1} \alpha \exp \bigl( \lambda(T-t_0) \bigr) \exp \bigl( - c/(T-t_0) \bigr)
\\
&\hspace{15pt}
- \lambda  \alpha 
{\mathbb E} \biggl[ \int_{t_0}^T  
\exp \biggl( \int_{t_0}^t
\bigl[ \lambda - \gamma \partial_e v(s,P_s,E_s) \bigr]
ds \biggr) dt \biggr].
\end{split}
\end{equation}
Therefore, for $\lambda <0$, absolute continuity follows as in the proof of Proposition 
\ref{prop:7:1:1} (with $D_s P_t = \sigma \exp (\lambda(t-s))$ in \eqref{eq:18:2:10}). Moreover, 
we can specify the \textit{size} of the transmission coefficient:
\begin{equation}
\label{eq:15:2:51}
\begin{split}
\alpha - \gamma \partial_p v(t_0,p,e) 
&\geq c^{-1} \alpha \exp \bigl( \lambda(T-t_0) \bigr) \exp \bigl( - c/(T-t_0) \bigr)
\\
&\hspace{15pt}
- \lambda \alpha  \exp \bigl( \lambda(T-t_0) \bigr)   \int_{t_0}^T  
\exp \bigl( \ln[(T-t)/(T-t_0)] \bigr)  dt
\\
&= \frac{\alpha}{c} \exp \bigl( \lambda(T-t_0) \bigr) \exp \bigl( - \frac{c}{T-t_0} \bigr)
- \frac{1}{2} \lambda \alpha  \exp \bigl( \lambda(T-t_0) \bigr) (T-t_0).
\end{split}
\end{equation}
The drift $\lambda$ makes the coefficient much larger than in the critical case when $\lambda=0$. In particular, we can compute the Malliavin derivative as
in \eqref{eq:18:2:10}. With
$D_s P_t = \sigma \exp (\lambda(t-s))$, we follow \eqref{eq:18:2:10}--\eqref{eq:3:8:4}
and write, for $T-t \geq (T-t_0)/2$ and $t_0 \leq s <t$:
\begin{equation}
\label{eq:18:2:20}
\begin{split}
&\bigl| {\mathbb E} \bigl[ \bigl( \sigma^{\top} \bigr)^{-1} D_s E_t |{\mathcal F}_s \bigr] \bigr|
\\
&\hspace{15pt} \geq \frac{T-t_0}{c}
{\mathbb E} \biggl[  \int_s^t 
\exp \biggl( -  \int_r^t \gamma \partial_e v(u,P_u,E_u) du \biggr) dr \bigl|{\mathcal F}_s 
\biggr]
\geq \frac{(T-t) (t-s)}{c}.
\end{split}
\end{equation}
Therefore, there exists $c \geq 1$ such that, for $T-t \geq (T-t_0)/2$,
\begin{equation}
\label{eq:18:2:21}
\text{ var } \bigl( E_t^{t_0,p,e} \bigr) \geq
c^{-1}(T-t_0)^2(t-t_0)^3. 
\end{equation} 

The bound is shown to be sharp for $(t_0,\bar{e}) \in {\mathcal C}''
=\{(t,e') \in [0,T) \times \RR : (e'-\Lambda)/(T-t) 
\in [4 \gamma/9,5 \gamma/9]\}$ and $T-t_0$ small enough by the same
argument as in Proposition \ref{prop:7:1:1}. In the critical area, the first 
term in the right-hand side in \eqref{eq:15:2:55} is exponentially small; the second one is always less than $|\lambda| \alpha(T-t_0)$. The end of the proof is an in 
\eqref{eq:18:2:20}.

\vskip 2pt
Finally, we consider the case $\lambda >0$.
Again, we can repeat the proof of the upper bound in Lemma
\ref{lem:2:8:1}. Passing to the limit in \eqref{eq:15:2:55} and applying 
(two-sided) Fatou's Lemma, there exists a constant $c \geq 1$ such that
that, for $(t_0,\bar{e}) \in {\mathcal C}''$ and $T-t_0$ small enough,
\begin{equation*}
\begin{split}
\alpha - \gamma \partial_p v(t_0,p,e) 
\leq  c \alpha \exp \bigl( \lambda(T-t_0) \bigr) \exp \bigl( - c^{-1}/(T-t_0) \bigr)
- \frac{1}{2} \lambda \alpha  \exp \bigl( \lambda(T-t_0) \bigr) (T-t_0).
\end{split}
\end{equation*}
Clearly, this says that, for $\bar{e}$ in 
${\mathcal C}''$ and $T-t_0$ small enough, the transmission coefficient
is negative! 

Now, for $\bar{e}$ away from the critical cone, we know from Proposition
\ref{prop:20:07:5} below that $\partial_e v$ doesn't explode and tends to $0$ as $t$ tends to $T$. In particular,
we let the reader check from \eqref{eq:15:2:55} that, for $\bar{e}$ away from the critical cone,
\begin{equation*}
\alpha - \gamma \partial_p v(t_0,p,e) 
\geq  c^{-1} - c (T-t_0).
\end{equation*}
As above, the proof consists in localizing the trajectories of the process
$\bar{E}$ and in taking advantage of the specific shape of $\partial_p w$.
Obviously, $c$ is independent of $p$ and $e$ provided that 
$\bar{e}$ is far enough from the critical cone. In particular, the transmission coefficient is positive for $\bar{e}$ away from the critical cone and $T-t_0$ small enough: by continuity of $\partial_p v$, it has a zero!
\qed

\subsection{Sufficient Condition for Hypoellipticity in Dimension $1$}
We provide a one-dimensional nonlinear example for which the first-order hypoelliptic condition holds and the conditional variance of the process $E$ is bounded from below. 

\begin{proposition}
\label{prop:20:07:1}
Assume that $d=1$ and that $f$ has the form $f(p,y) = -f_0(\mu p-y)$, for some real $\mu >0$ and some continuously differentiable function $f_0 : \RR \rightarrow \RR$ satisfying $\ell_1 \leq f_0' \leq \ell_2$ with $\ell_1$, $\ell_2$ and $L$ as in (A.1-2). Assume also that $b$ and $\sigma$ 
satisfy (A.1-2) w.r.t. $L$ and 
are continuously differentiable with H\"older continuous derivatives, that $\partial_p b(p) \leq 0$ for $p \in \RR$, and that $\sigma$ is bounded by $L$
and satisfies $\inf_{p \in \RR} \sigma(p) \geq L^{-1} >0$. Then, for any initial condition $(t_0,p,e)$, the process $(E_t^{t_0,p,e})_{t_0 \leq t < T}$ has absolutely continuous marginal distributions. Moreover, if there exists a constant 
$\lambda \in (0,L]$ such that $\partial_p b(p) \leq - \lambda$ for any $p \in \RR$, then,
there exists a constant $c \geq 1$, depending on known parameters only (but not on 
$t_0$), such that, for any $t_0 \leq t \leq (T+t_0)/2$,
\begin{equation*}
\text{ var }(E_t^{t_0,p,e}) \geq c^{-1}(T-t)^2(t-t_0)^3. 
\end{equation*}
\end{proposition}

The role of the assumption $d=1$ is twofold. 
First, it permits to specify the form of $f$. Second, the proof relies
on a variation of the strong maximum principle for the PDE satisfied
by $\partial_p v$, and, in higher dimension, $\partial_p v$ satisfies
a system of partial differential equations for which a maximum principle 
of this type is more problematic. From an intuitive point of view, the restriction to the one-dimensional setting is not satisfactory: additional noise in the dynamics for $P$ should favor non-degeneracy of the dynamics of $E$, so that increasing the dimension should help and not be a hindrance. We leave this question to further investigations.
\vspace{5pt}

We also emphasize that we say nothing about the existence of a density to the pair process. The reason is rather technical: we know very little about the smoothness of the value function $v$ so that standard hypoellipticity arguments fail. In Delarue and Menozzi \cite{DelarueMenozzi2}, the existence of a density for the pair is investigated under ${\mathcal C}^1$ conditions (that is much less than standard hypoellipticity conditions), but here we are unable to check these conditions. When investigating the process $E$ separately, we avoid the smoothness conditions by using the Bouleau and Hirsch criterion, which is very simple to check. In order to apply the Bouleau and Hirsch criterion to the pair process, it would be necessary to check the invertibility of the corresponding Malliavin matrix: it seems that extra smoothness is needed in order to be able to do it. 

\vskip 2pt\noindent
{\bf Proof.}
\textit{First Step.} Since nothing is known about the smoothness of $v$, we start with the case when the terminal condition $v^{\phi}(T,p,e)=\phi(e)$ is a smooth
non-decreasing
 function with values in $[0,1]$. The point is then to prove an estimate for
 the transmission coefficient $\partial_p [f(p,v^{\phi}(t,p,e))]$, independently of the smoothness of $\phi$. It is shown in Proposition \ref{prop:08:03:5} below that $v^{\phi}$ is continuously differentiable when $\phi$ is smooth and that, with 
 $(P,E^{\phi},Y^{\phi}) = (P^{t_0,p},E^{\phi,t_0,p,e},Y^{\phi,t_0,p,e})$ and
$$\frac{d {\mathbb Q}}{d {\mathbb P}}
= \exp \biggl( \int_{t_0}^T \partial_p \sigma(P_s) dW_s - \frac{1}{2}
\int_{t_0}^T \bigl[ \partial_p \sigma(P_s) \bigr]^2 ds 
\biggr),$$
it holds 
\begin{equation*}
\begin{split}
\partial_p v^{\phi}(t_0,p,e)
&= - {\mathbb E}^{\mathbb Q}
\biggl[ \int_{t_0}^T \partial_e v^{\phi}(t,P_t,Y_t^{\phi}) \partial_p f(t,P_t,Y_t^{\phi}) 
\\
&\hspace{-10pt} \times
\exp \biggl( \int_{t_0}^t 
\bigl[
- \partial_y f(P_s,Y_s^{\phi}) \partial_e v^{\phi}(s,P_s,E_s^{\phi}) 
+ \partial_p b(P_s) \bigr] ds \biggr)
dt
\biggr],
\end{split}
\end{equation*}
that is
\begin{equation*}
\begin{split}
\partial_p v^{\phi}(t_0,p,e)
&=  \mu {\mathbb E}^{\mathbb Q}
\biggl[ \int_{t_0}^T 
\exp \biggl( \int_{t_0}^t 
 \partial_p b(P_s) ds \biggr)
\\
&\hspace{-10pt} \times
\partial_y f(P_t,Y_t^{\phi}) \partial_e v^{\phi}(t,P_t,Y_t^{\phi}) 
\exp \biggl( \int_{t_0}^t 
- \partial_y f(P_s,Y_s^{\phi}) \partial_e v^{\phi}(s,P_s,E_s^{\phi}) 
 ds \biggr)
dt
\biggr],
\end{split}
\end{equation*}
by the specific form of $f$: 
$[\partial_p f/\partial_y f](p,e) = - \mu$. As in the linear counter-example, 
we then make an integration by parts. We obtain
\begin{equation}
\label{eq:4:8:1}
\begin{split}
&1 - \mu^{-1} \partial_p v^{\phi}(t_0,p,e)
\\
&=  {\mathbb E}^{\mathbb Q}
\biggl[ 
\exp \biggl( \int_{t_0}^T
 \partial_p b(P_s) ds \biggr)
\exp \biggl( -\int_{t_0}^T
\partial_y f(P_s,Y_s^{\phi}) \partial_e v^{\phi}(s,P_s,E_s^{\phi}) 
 ds \biggr) \biggr]
 \\
 &\hspace{5pt}
 - {\mathbb E}^{\mathbb Q}
 \biggl[ \int_{t_0}^T \partial_p b(P_t) 
 \exp \biggl( 
 \int_{t_0}^t 
\bigl[ \partial_p b(P_s)- \partial_y f(P_s,Y_s^{\phi}) \partial_e v^{\phi}(s,P_s,E_s^{\phi}) 
\bigr] ds
  \biggr)
dt
\biggr].
\end{split}
\end{equation}

\textit{Second Step.}  When $\partial_p b(P_s) \leq - \lambda$, we can follow 
\eqref{eq:15:2:55}--\eqref{eq:15:2:51} and then obtain (still in the smooth setting):
\begin{equation}
\label{eq:18:2:11}
\begin{split}
1- \mu^{-1} \partial_p v^{\phi}(t_0,p,e)
&\geq \lambda \exp\bigl[-L(T-t_0) \bigr] \int_{t_0}^T 
\exp \bigl(  L^2 \ln\bigl[(T-t)/(T-t_0)\bigr] \bigr) dt
\\
&\geq \lambda (L^2+1)^{-1} (T-t_0) \exp \bigl[ - L(T-t_0) \bigr].
\end{split}
\end{equation}

We wish we could pass to the limit in \eqref{eq:18:2:11} along a mollification of the terminal condition. Obviously, we can't since 
$\partial_p v$ is not known to exist in the singular setting. Nevertheless, 
\eqref{eq:18:2:11} says that the function $\RR \ni p \hookrightarrow p- \mu^{-1} v(t,p,e)$ is increasing (in the singular setting), the monotonicity constant being greater than $\lambda (L^2+1)^{-1} (T-t_0) \exp [ - L(T-t_0)]$. Below, we will consider the singular setting only and \eqref{eq:18:2:11} will be understood as a lower bound on the Lipschitz constant of $v$.

By Theorem 2.2.1 in Nualart \cite{Nualart}, $P$ is known to be differentiable in the Malliavin sense and 
\begin{equation}
\label{eq:18:2:12}
D_s P_t = \sigma(P_s) \exp\biggl(\int_s^t \partial_p b(P_r) dr
+ \int_s^t \partial_p \sigma(P_r) dW_r - (1/2) \int_s^t [\partial_p \sigma(P_r)]^2 dr
\biggr).
\end{equation}
(Here and below, $(P_{t_0},E_{t_0})=(p,e)$.)
By Proposition 1.2.3 and Theorem 2.2.1 in \cite{Nualart}, we can also compute the Malliavin derivative of 
$E$ despite the lack of differentiability of $v$. Following \eqref{eq:18:2:10},
\begin{equation*}
D_s E_t
=  \int_s^t \bigl[ \mu - \partial_p v(r,P_r,E_r)
\bigr] \partial_y f(P_r,Y_r) 
D_s P_r \exp \biggl( - \int_r^t \partial_y f(P_u,Y_u)
\partial_e v(u,P_u,E_u) du \biggr) dr.
\end{equation*}
Here, $(\partial_p v(r,P_r,E_r))_{t_0 \leq r <T}$ and $(\partial_e v(r,P_r,E_r))_{t_0 \leq r <T}$ stand for progressively-measurable processes that coincide with the true derivatives whenever they exist and are continuous. Following the proof of Proposition 1.2.3 in \cite{Nualart}, processes $(\partial_p v(r,P_r,E_r))_{t_0 \leq r <T}$ and $(\partial_e v(r,P_r,E_r))_{t_0 \leq r <T}$ are constructed by approximating $v$ by a standard convolution argument: by the bounds we have on the Lipschitz and monotonicity constants of 
$v$ with respect to $p$ and $e$, we deduce that the process
$(\mu - \partial_p v(r,P_r,E_r))_{t_0 \leq r <T}$ 
is bounded and greater than $(\mu \lambda (L^2+1)^{-1} (T-r) \exp [ - L(T-r)])_{t_0 \leq r <T}$  and that the process $(\partial_e v(r,P_r,E_r))_{t_0 \leq r <T}$ is non-negative and less than $(L(T-r)^{-1})_{t_0 \leq r <T}$. Therefore, Bouleau and Hirsch criterion applies as in \eqref{eq:18:2:10}.
From \eqref{eq:18:2:11}, \eqref{eq:18:2:12} and the identity above, we deduce, as in 
\eqref{eq:18:2:20},
\begin{equation*}
\bigl| {\mathbb E} \bigl[D_s E_t | {\mathcal F}_s \bigr]
\bigr| \geq c^{-1} (T-t_0) \biggl| \int_s^t  \bigl( \frac{T-t}{T-r} \bigr)^{L^2} dr \biggr|, \quad t_0 \leq t \leq \frac{T+t_0}{2},
\end{equation*}
the constant $c \geq 1$ being here a positive constant depending on 
known parameters only, and not on $t_0$. We then complete the lower bound for the variance by Clark-Ocone formula.
\vspace{5pt}

\textit{Third Step.} We now consider the case when $\partial_p b(P_s) \leq 0$ only.
We then approximate the Heaviside terminal condition by a sequence $(\phi^n)_{n \geq 1}$ of smooth terminal conditions. Since
the lower bound in \eqref{eq:18:2:11} fails for any $n \geq 1$, 
the point is to consider the first term only in the right-hand side in \eqref{eq:4:8:1} and to bound it from below. Assuming that $\phi^n(e)=1$ for $e \geq \Lambda$, we recover \eqref{eq:4:8:10}, but under the probability ${\mathbb Q}$. Passing to the limit, we deduce the analogue of the left-hand side in \eqref{eq:1:8:1}:
\begin{equation}
\label{eq:4:8:11}
\liminf_{n \rightarrow + \infty} 
\bigl[ 1 - \mu^{-1} \partial_{p_i} v^n(t_0,p,e) \bigr]
\geq C^{-1} 
{\mathbb Q} \bigl\{ \inf_{(t_0+T)/2 \leq t \leq T} \bigl[\bar{E}_t - \Lambda \bigr] >  C(T-t_0)
\bigr\},
\end{equation}
with $v^n = v^{\phi^n}$. (Above, $\bar{E}_t = E_t + w(t,P_t)$, with 
$(P_{t_0},E_{t_0})=(p,e)$.) In the right-hand side above, we can switch back from ${\mathbb Q}$ to ${\mathbb P}$
since ${\mathbb P}(A) \leq C ({\mathbb Q}(A))^{1/2}$, $A \in {\mathcal F}$, for some $C>0$. It then remains to bound from below ${\mathbb P} \{ \inf_{(t_0+T)/2 \leq t \leq T} [\bar{E}_t - \Lambda ] >  C(T-t_0)\}$. We then follow \eqref{eq:14:2:5}, replacing $|\bar{E}-\Lambda|$ therein by $\bar{E}-\Lambda$. 

The estimate of $\pi_2$ in \eqref{eq:14:2:6} is then similar.
(The bound of $\partial_p w$ is kept preserved even if (A.4) may not be satisfied: go back to the statement of Lemma \ref{lem:22:12:3}.)
To estimate $\pi_1$ (without the absolute value), it is sufficient to bound from below
\begin{equation*}
\pi_1'=
{\mathbb P}\biggl\{ \int_{t_0}^{(T+t_0)/2}  \sigma(P_s) \partial_p w(s,P_s) dW_s  \geq \bigl(C+\frac{L}{2}+1 \bigr) (T-t_0) \biggr\}
\end{equation*}
when $\bar{e} \geq \Lambda$. As in the proof of Proposition \ref{prop:4:1:1}, we can prove that $\partial_p w$ is bounded from below. Indeed, by  \eqref{eq:09:03:5}, we know that 
$\partial_p w$ has the form
\begin{equation*}
\partial_p w(t,p) = \mu {\mathbb E} \biggl[ \int_t^T 
f_0'( \mu P_s^{t,p}) \partial_p P_s^{t,p} ds \biggr].
\end{equation*}
Since we are in dimension 1,
\begin{equation*}
\partial_p P_s^{t,p} = \exp \biggl( \int_t^s \partial_p b(P_r) dr
+ \int_t^s \partial_p \sigma(P_r) dW_r - \frac{1}{2}
\int_t^s [\partial_p \sigma(P_r)]^2 dr \biggr), 
\end{equation*}
so that
\begin{equation}
\label{eq:09:03:2}
\partial_p w(t,p) \geq \mu L^{-1} 
{\mathbb E}^{\mathbb Q} \biggl[ \int_t^T \exp \biggl( \int_t^s 
\partial_p b(P_r) dr \biggr) ds \biggr]
\geq \mu L^{-1} \exp ( -L T) (T-t). 
\end{equation}

To bound $\pi_1'$, we can proceed as follows. Setting $(\theta_t
= \sigma(P_t) \partial_p w(t,P_t))_{t_0 \leq t \leq T}$, we can find some constants $c,C'\geq 1$ such that, for any $a>0$, 
\begin{equation*}
\begin{split}
1 &= {\mathbb E} \biggl[ \exp \biggl( a\int_{t_0}^{(T+t_0)/2}
\theta_s dW_s - \frac{a^2}{2} \int_{t_0}^{(T+t_0)/2} \theta_s^2 ds \biggr) \biggr]
\\
&\leq \exp \bigl( - c^{-2} a^2 (T-t_0)^3 \bigr)
{\mathbb E} \biggl[ \exp \biggl( a\int_{t_0}^{(T+t_0)/2}
\theta_s dW_s  \biggr) \biggr]
\\
&\leq \exp \bigl( - c^{-2} a^2 (T-t_0)^3 \bigr) \exp \bigl( a [C+(L/2)+1] (T-t_0) \bigr)
\\
&\hspace{15pt} + \exp \bigl( -  c^{-2} a^2 (T-t_0)^3 \bigr)
{\mathbb E} \biggl[ \exp \biggl( 2 a\int_{t_0}^{(T+t_0)/2}
\theta_s dW_s  \biggr) \biggr]^{1/2} (\pi_1')^{1/2}
\\
&\leq \exp \bigl( - c^{-2} a^2 (T-t_0)^3  \bigr) \exp \bigl( a [ C+ (L/2)+1](T-t_0) \bigr)
 + \exp \bigl(C' a^2 (T-t_0)^3 \bigr) (\pi_1')^{1/2}.
\end{split}
\end{equation*}
Obviously, $c$ and $C'$ are independent of $t_0$ and $p$. Choosing 
$a(T-t_0)^2$ large enough, we can make the first term in the above right-hand side as small as desired. We deduce that $\pi_1' \geq \exp(-C''(T-t_0)^{-1})$, 
for some $C'' >0$. We recover \eqref{eq:09:03:1}. 
From \eqref{eq:4:8:11}, we deduce that $\liminf_{n \rightarrow + \infty} [1-\mu^{-1} \partial_p v^n(t_0,p,e)] \geq (C'')^{-1} \exp(-C''(T-t_0)^{-1})$ when $\bar{e} = e+ w(t_0,p) \geq \Lambda$ (modifying $C''$ if necessary). Since $w(t_0,\cdot)$ is continuous and increasing, the set $I_+ = \{p \in \RR : e+w(t,p) \geq \Lambda\}$ is an interval: passing to the limit, we deduce that the function $I_+ \ni p \hookrightarrow p - \mu^{-1} v(t_0,p,e)$ is increasing, the monotonicity constant being greater than 
$(C'')^{-1} \exp(-C''(T-t_0)^{-1})$. Choosing the approximating sequence $(\phi^n)_{n \geq 1}$ such that $\phi^n(e)=0$ for $e \leq \Lambda$ and repeating the argument, we obtain a similar bound on the set $I_- = \{p \in \RR : e+w(t,p) \leq \Lambda\}$.
Absolute continuity of the density of $E_t$, $t<T$, then follows by Bouleau and Hirsch criterion again.
\qed

\section{Appendix}
\label{se:appendix}
\subsection{Smoothness of $v$ for a Smooth Terminal Condition}

\begin{proposition}
\label{prop:08:03:5}
Assume that, in addition to (A.1), (A.2), the coefficients $b$ and $\sigma$ are continuously differentiable with H\"older continuous derivatives, that the function $f$ is continuously differentiable w.r.t. $p$ and $y$ and that $\sigma \sigma^{\top}$ is uniformly non-degenerate as in \eqref{eq:4:8:15}. Assume also that the terminal condition $\phi$ is smooth. Then, the function
$v^{\phi}$ is continuously differentiable with respect to 
$p$ and $e$ on $[0,T) \times \RR^d \times \RR$. 

Moreover, if $d=1$, then 
\begin{equation*}
\begin{split}
\partial_p v^{\phi}(t_0,p,e)
&= - {\mathbb E}^{{\mathbb Q}}
\biggl[ \int_{t_0}^T \partial_e v(t,P_t^{t_0,p},Y_t^{\phi,t_0,p,e}) \partial_p f(t,P_t^{t_0,p},Y_t^{\phi,t_0,p,e}) 
\\
&\hspace{10pt} \times
\exp \biggl( \int_{t_0}^t 
\bigl[
- \partial_y f(P_s^{t_0,p},Y_s^{\phi,t_0,p,e}) \partial_e v(s,P_s^{t_0,p},E_s^{\phi,t_0,p,e}) 
+ \partial_p b(P_s^{t_0,p}) \bigr] ds \biggr)
dt
\biggr].
\end{split}
\end{equation*}
with
\begin{equation*}
\frac{d {\mathbb Q}}{d {\mathbb P}}
= \exp \biggl( \int_{t_0}^T \partial_p \sigma(P_s^{t_0,p}) dW_s - \frac{1}{2}
\int_{t_0}^T \bigl[ \partial_p \sigma(P_s^{t_0,p}) \bigr]^2 ds 
\biggr).
\end{equation*}
\end{proposition}

{\bf Proof.} 
\textit{First Step: Continuous Differentiability of $v^{\phi}$.} 
When the coefficients $b$, $\sigma$, $f$ and $\phi$ are smooth with bounded derivatives of any order, we know from Lemma 3.5 in \cite{CrisanDelarue} (and by the bound \reff{eq:9:11:1}) that $v^{\varepsilon}$ in \eqref{eq:25:7:5}--\eqref{eq:7:11:10} is infinitely differentiable w.r.t. $(p,e)$ on $(0,T) \times \RR^d \times \RR$ with time-space continuous derivatives of any order. Differentiating 
Eq. \eqref{eq:7:11:10}, we obtain
\begin{equation}
\label{eq:6:3:1}
\begin{split}
&\partial_t \bigl[ \partial_e v^{\varepsilon}\bigr](t,p,e) + {\mathcal L}_p 
\bigl[ \partial_e v^{\varepsilon}\bigr](t,p,e)
+ \frac{\varepsilon^2}{2} \partial^2_{pp} \bigl[\partial_e v^{\varepsilon}\bigr](t,p,e)
+ \frac{\varepsilon^2}{2} \partial^2_{ee} \bigl[\partial_e v^{\varepsilon}\bigr](t,p,e) 
\\
&\hspace{15pt}
- f\bigl(p,v^{\varepsilon}(t,p,e)\bigr) \partial_e \bigl[\partial_e v^{\varepsilon}\bigr](t,p,e) 
- \partial_y f\bigl(p,v^{\varepsilon}(t,p,e) \bigr) 
\bigl[ \partial_e v^{\varepsilon} \bigr]^2 (t,p,e) = 0,
\end{split}
\end{equation} 
with $\partial_e v^{\varepsilon}(T,p,e) = \phi'(e)$ as terminal 
condition. By It\^o's formula,
\begin{equation*}
\partial_e v^{\varepsilon}(t_0,p,e)
= {\mathbb E} \biggl[ 
\phi'\bigl(E_T^{\varepsilon,t_0,p,e} \bigr)
- \int_{t_0}^T 
\partial_y f\bigl(P_s^{\varepsilon,t_0,p},Y_s^{\varepsilon,t_0,p,e} \bigr)
\bigl[  \partial_e v^{\varepsilon}\bigl(s,P_s^{\varepsilon,t_0,p},E_s^{\varepsilon,t_0,p,e} \bigr)
\bigr]^2 ds \biggr],
\end{equation*}
with the same notation as in \eqref{eq:4:11:5}. Following the approximation argument 
in Corollary \ref{cor:25:7:1}, the above expression holds true under the assumption of Proposition \ref{prop:08:03:5}. 
We deduce that the process $(\Psi_t^{\varepsilon}= \partial_e v^{\varepsilon}(t,P_t^{\varepsilon,t_0,p},E_t^{\varepsilon,t_0,p,e}))_{t_0 \leq t \leq T}$ satisfies the BSDE:
\begin{equation*}
d \Psi_t^{\varepsilon} = \partial_y f\bigl(P_t^{\varepsilon,t_0,p},Y_t^{\varepsilon,t_0,p,e} \bigr)
\bigl[ \Psi_t^{\varepsilon} \bigr]^2 dt + dm_t^{\varepsilon}, \quad t_0 \leq t \leq T,
\end{equation*}
with $\Psi_T^{\varepsilon} = \phi'(E_T^{\varepsilon,t_0,p,e})$ as terminal condition,
$(m_t^{\varepsilon})_{t_0 \leq t\leq T}$ standing for 
a martingale term. Therefore, for $\varepsilon,\varepsilon'>0$,
\begin{equation*}
\begin{split}
d \bigl[ (\Psi_t^{\varepsilon} - \Psi_t^{\varepsilon'} )^2 \bigr]
&= 2 \partial_y f(P_t^{\varepsilon,t_0,p},Y_t^{\varepsilon,t_0,p,e})
\bigl[ \Psi_t^{\varepsilon} + \Psi_t^{\varepsilon'} \bigr]
\bigl[ (\Psi_t^{\varepsilon} - \Psi_t^{\varepsilon'} )^2 \bigr] 
\\
&\hspace{15pt} + 
2 \bigl[ \partial_y f(P_t^{\varepsilon,t_0,p},Y_t^{\varepsilon,t_0,p,e})
- \partial_y f(P_t^{\varepsilon',t_0,p},Y_t^{\varepsilon',t_0,p,e}) \bigr]
\bigl[ \Psi_t^{\varepsilon'} \bigr]^2 
\bigl[ \Psi_t^{\varepsilon} - \Psi_t^{\varepsilon'}  \bigr] dt 
\\
&\hspace{15pt} 
+ d \bigl[ m^{\varepsilon} - m^{\varepsilon'} \bigr]_t + dn_t^{\varepsilon,\varepsilon'},
\end{split}
\end{equation*}
$(n_t^{\varepsilon,\varepsilon'})_{t_0 \leq t \leq T}$ standing for a new martingale term and 
$[m^{\varepsilon} - m^{\varepsilon'}]$ for the quadratic variation of 
$m^{\varepsilon} - m^{\varepsilon'}$. 
Since $\partial_y f \geq 0$ and $\Psi^{\varepsilon}$ is non-negative and bounded by a constant depending on 
$L$, $T$ and $\|\phi'\|_{\infty}$ only (see Proposition \ref{prop:4:11:1}), 
\begin{equation*}
\begin{split}
\bigl|\partial_e v^{\varepsilon}(t_0,p,e) - \partial_e v^{\varepsilon'}(t_0,p,e) \bigr|^2
&\leq {\mathbb E} \bigl[ \bigl( \phi'(E_T^{\varepsilon,t_0,p,e})
- \phi'(E_T^{\varepsilon',t_0,p,e}) \bigr)^2 \bigr]
\\
&\hspace{5pt} + C {\mathbb E} \int_{t_0}^T 
\bigl| \partial_y f(P_t^{\varepsilon,t_0,p},Y_t^{\varepsilon,t_0,p,e})
- \partial_y f(P_t^{\varepsilon',t_0,p},Y_t^{\varepsilon',t_0,p,e}) \bigr| dt.
\end{split}
\end{equation*}
This is sufficient to prove that the sequence $(\partial_e v^{\varepsilon})_{\varepsilon >0}$ is uniformly convergent on compact subsets 
as $\varepsilon$ tends to $0$: the limit is $\partial_e v^{\phi}$; it is continuous. 

We now investigate $\partial_p v^{\phi}$. Using the same notation as in 
\eqref{eq:4:11:5} and applying  It\^o's formula to the process $(v^{\varepsilon}(s,P_s^{\varepsilon,t_0,p},e+\varepsilon B_{s-t_0}))_{t_0 \leq s \leq T}$, 
we obtain
\begin{equation*}
\begin{split}
v^{\varepsilon}(t_0,p,e) &= {\mathbb E} \biggl[ 
\phi(e+\varepsilon B_{T-t_0})
\\
&\hspace{15pt}
- \int_{t_0}^T f\bigl(P_s^{\varepsilon,t_0,p},v^{\varepsilon}
(s,P_s^{\varepsilon,t_0,p},e+\varepsilon B_{s-t_0}) \bigr)
\partial_e v^{\varepsilon}\bigl(s,P_s^{\varepsilon,t_0,p},e+\varepsilon B_{s-t_0} \bigr) ds \biggr].
\end{split}
\end{equation*}
Since $\partial_e v^{\varepsilon}$ converges towards $\partial_e v^{\phi}$ uniformly on compact subsets, we deduce
\begin{equation*}
v^{\phi}(t_0,p,e) = 
\phi(e)
- {\mathbb E} \biggl[ 
 \int_{t_0}^T f\bigl(P_s^{t_0,p},v^{\phi}
(s,P_s^{t_0,p},e) \bigr)
\partial_e v^{\phi}\bigl(s,P_s^{t_0,p},e\bigr) ds \biggr].
\end{equation*}
By Bismut-Elworthy formula (see Nualart \cite{Nualart})
\begin{equation*}
\begin{split}
&\partial_p v^{\phi}(t_0,p,e) 
\\
&= 
-  \int_{t_0}^T \frac{1}{s-t_0} {\mathbb E}  \biggl[ f\bigl(P_s^{t_0,p},v^{\phi}
(s,P_s^{t_0,p},e) \bigr)
\partial_e v^{\phi}\bigl(s,P_s^{t_0,p},e\bigr)
\int_{t_0}^s 
 \bigl( \sigma^{-1}(P_r^{t_0,p})  \partial_p P_r^{t_0,p} \bigr)^{\top}
dW_r  \biggr] ds.
\end{split}
\end{equation*}
The term inside the integral from $t_0$ to $T$ is bounded by $(s-t_0)^{-1/2}$, uniformly in $t_0$, $s$, $p$ and $e$ in compact sets. By uniform integrability, continuity
of $\partial_p v^{\phi}$ easily follows.  
\vspace{5pt}

\textit{Second Step: Representation of $\partial_p v^{\phi}$.}
We now assume that $d=1$. We then make use of \eqref{eq:13:03:10}. With $U_t^{\varepsilon}= \partial_p v^{\varepsilon}(t,P_t^{\varepsilon,t_0,p},E_t^{\varepsilon,t_0,p,e})$
and $V_t^{\varepsilon} = \partial^2_{pp} v^{\varepsilon}(t,P_t^{\varepsilon,t_0,p},E_t^{\varepsilon,t_0,p,e})$
for some initial condition $(t_0,p,e)$, we obtain as dynamics for $U^{\varepsilon}$ (below, we do not specify $(t_0,p,e)$ in $P^{\varepsilon}$ and $E^{\varepsilon}$), 
\begin{equation}
\label{eq:16:12:1}
\begin{split}
dU_t^{\varepsilon} &= - \partial_p b(P_t^{\varepsilon}) U_t^{\varepsilon} dt - 
\sigma(P_t^{\varepsilon}) \partial_p \sigma(P_t^{\varepsilon}) V_t^{\varepsilon}  dt 
\\
&\hspace{15pt} + \partial_e v^{\varepsilon}(t,P_t^{\varepsilon},E_t^{\varepsilon})
\bigl[ \partial_p f(P_t^{\varepsilon},Y_t^{\varepsilon}) + \partial_y f(P_t^{\varepsilon},Y_t^{\varepsilon}) U_t^{\varepsilon} \bigr] dt  
\\
&\hspace{15pt} + V_t^{\varepsilon} \sigma(P_t^{\varepsilon}) dW_t
+ \varepsilon V_t^{\varepsilon} dW_t' + \varepsilon 
\partial_{pe}^2 v^{\varepsilon}(t,P_t^{\varepsilon},E_t^{\varepsilon}) dB_t, \quad 0 \leq t \leq T,
\end{split}
\end{equation}
with $U_T^{\varepsilon}=0$ as terminal condition. We then introduce the exponential weight:
\begin{equation*}
\begin{split}
&{\mathcal E}_t^{\varepsilon} 
\\
&=
\exp \biggl( \int_{t_0}^t 
\bigl[-
\partial_y f(P_s^{\varepsilon},Y_s^{\varepsilon}) \partial_e v^{\varepsilon}(s,P_s^{\varepsilon},E_s^{\varepsilon}) 
+ \partial_p b(P_s^{\varepsilon}) 
- \frac{1}{2}
\bigl( \partial_p \sigma(P_s^{\varepsilon}) \bigr)^2\bigr] ds+
\int_{t_0}^t \partial_p \sigma(P_s^{\varepsilon}) dW_s 
\biggr).
\end{split}
\end{equation*}
By It\^o's formula, we obtain the following Feynman-Kac formula:
\begin{equation*}
U_{t_0}^{\varepsilon} = - {\mathbb E} \biggl[ \int_{t_0}^T
\partial_e v^{\varepsilon}(t,P_t^{\varepsilon},Y_t^{\varepsilon}) \partial_p f(P_t^{\varepsilon},Y_t^{\varepsilon}) {\mathcal E}_t^{\varepsilon} dt \biggr].
\end{equation*}
Defining the density
\begin{equation*}
\frac{d {\mathbb Q}^{\varepsilon}}{d {\mathbb P}}
= \exp \biggl( \int_{t_0}^T \partial_p \sigma(P_s^{\varepsilon}) dW_s - \frac{1}{2}
\int_{t_0}^T \bigl[ \partial_p \sigma(P_s^{\varepsilon}) \bigr]^2 ds 
\biggr),
\end{equation*}
we get:
\begin{equation*}
\begin{split}
U_{t_0}^{\varepsilon} 
&= - {\mathbb E}^{{\mathbb Q}^{\varepsilon}}
\biggl[ \int_{t_0}^T \partial_e v^{\varepsilon}(t,P_t^{\varepsilon},Y_t^{\varepsilon}) \partial_p f(t,P_t^{\varepsilon},Y_t^{\varepsilon}) 
\\
&\hspace{30pt} \times
\exp \biggl( \int_{t_0}^t 
\bigl[
- \partial_y f(P_s^{\varepsilon},Y_s^{\varepsilon}) \partial_e v^{\varepsilon}(s,P_s^{\varepsilon},E_s^{\varepsilon}) 
+ \partial_p b(P_s^{\varepsilon}) \bigr] ds \biggr)
dt
\biggr],
\end{split}
\end{equation*}
and letting $\varepsilon$ tend to $0$, we complete the proof. \qed

\begin{proposition}
\label{prop:20:07:5}
Assume that, in addition to (A.1), (A.2), the coefficients $b$ and $\sigma$ are continuously differentiable with H\"older continuous derivatives, that the function $f$ is continuously differentiable w.r.t. $p$ and $y$ and that 
$\sigma$ is bounded by $L$ and $\sigma \sigma^{\top}$ uniformly non-degenerate as in 
(A.4) and \eqref{eq:4:8:15}. (Pay attention that $b$ may not be bounded.) Assume also that the terminal condition $\phi$ is smooth. Then, there exists a constant $C\geq 1$, depending on $L$ and $T$ only, such that,
if the terminal condition $\phi$ satisfies $\phi(e)=1$ for $e \geq \Lambda$, then, for any $(t_0,p,e) \in [0,T) \times \RR^d \times \RR$,
\begin{equation}
\label{eq:1:8:2}
\bar{e} = e + w(t_0,p) > \Lambda + C (T-t_0) \Rightarrow \partial_e v^{\phi}(t,p,e) \leq C (T-t_0)^2.
\end{equation}

Similarly, if $\phi$ satisfies $\phi(e)=0$ for $e \leq \Lambda$, then, for any $(t_0,p,e) \in [0,T) \times \RR^d \times \RR$,
\begin{equation}
\label{eq:1:8:3}
\bar{e} = e + w(t_0,p) < \Lambda - C (T-t_0) \Rightarrow \partial_e v^{\phi}(t,p,e) \leq C (T-t_0)^2.
\end{equation}
\end{proposition}

{\bf Proof.} We prove (\ref{eq:1:8:2}) only, i.e.  we assume $\phi(e)=1$ for $e \geq \Lambda$. Given an initial condition $(t_0,p,e)$, we consider
$\partial_e E^{\phi}$ and for convenience we omit to specify the initial condition in 
$P$ and $E^{\phi}$. Differentiating $(Y^{\phi}_t=v^{\phi}(t,P_t,E_t^{\phi}))_{t_0 \leq t <T}$ w.r.t. $e$, we deduce that 
$(M_t = \partial_e v^{\phi}(t,P_t,E_t^{\phi}) \partial_e E_t^{\phi})_{t_0 \leq t \leq T}$ is a bounded martingale.
(By Proposition \ref{prop:08:03:5}, $\partial_e v^{\phi}$ exists as a true function.
By Proposition \ref{prop:4:11:1}, it is bounded up to the boundary.)

When $E_T^{\phi} > \Lambda$, $\partial_e v^{\phi}(T,P_T,E_T)$ is zero, so that $M_T=0$ as well.
Since $\bar{E}_T^{\phi}= E_T^{\phi}$, we also have that $M_T=0$ when $\bar{E}_T^{\phi} > \Lambda$. (Here $\bar{E}_t^{\phi} = E_t^{\phi} + 
w(t,P_t)$ as in \eqref{eq:7:1:3}.)
We now use the hitting time $\tau$ defined in Lemma \ref{lem:18:12:1} below. We also
assume that $\bar{e} = e + w(t_0,p) \geq \Lambda + (L+1) (T-t_0)$. On the event
$\{\tau =T\}$, it is proven below that $\bar{E}_T^{\phi} \geq \Lambda + (T-t_0)$, so that $M_T=0$. By Doob's Theorem, we get that $\partial_e v^{\phi}(t_0,p,e)={\mathbb E}[M_{\tau}]$,
so that $\partial_e v^{\phi}(t_0,p,e) = {\mathbb E}[M_{\tau}; \tau < T ]$.

Since 
$\partial_e v^{\phi}(\tau,P_{\tau},E_{\tau}) \leq L (T-\tau)^{-1}$
on the event $\{\tau < T\}$
(see Proposition \ref{prop:4:11:1}), we deduce 
\begin{equation}
\label{eq:5:8:1}
\partial_e v^{\phi}(t_0,p,e) \leq L {\mathbb E} \bigl[(T-\tau)^{-1}; \tau < T \bigr]\leq L \sum_{n \geq 0} (T-t_{n+1})^{-1} {\mathbb P}\bigl\{\tau \in (t_n,t_{n+1}] \bigr\},
\end{equation}
for the net $(t_n = t_0 + (T-t_0) (1- 2^{-4n}))_{n \geq 0}$ given in Lemma \ref{lem:18:12:1}.

Therefore, $(T-t_{n+1})^{-1} = 2^{4(n+1)}(T-t_0)^{-1}$ in \eqref{eq:5:8:1}. By 
\eqref{eq:5:8:2} in Lemma \ref{lem:18:12:1}, we deduce that $\partial_e v^{\phi}(t_0,p,e) \leq  C  (T-t_0)^2$. The proof is similar when $\bar{e} \leq \Lambda - (L+1) (T-t_0)$.
\qed

\begin{lemma}
\label{lem:18:12:1}
For a given $t_0 \in [0,T)$, define $\displaystyle \delta_n = 15(T-t_0) 2^{-4n}$, $\displaystyle
h_n = 16 L(T-t_0) 2^{-n}$, $n \in {\mathbb N}^*$. Given the time and space nets
\begin{equation*}
\begin{split}
&\bigl[t_n = t_0 + \sum_{k=1}^n \delta_k = t_0 + (T-t_0) \frac{15}{16} \sum_{k=0}^{n-1} 16^{-k}
= t_0 + (T-t_0) \bigl(1 - 2^{-4n} \bigr) \bigr]_{n \geq 1},
\\
&\bigl[H_n = \sum_{k=1}^n h_k = 16 L \frac{T-t_0}{2} \sum_{k=0}^{n-1} 2^{-k}
= 16  L (T-t_0) \bigl( 1 - 2^{-n} \bigr) \bigr]_{n \geq 1},
\end{split}
\end{equation*}
we define the envelope function $\psi_t = H_1 {\mathbf 1}_{[t_0,t_1]}(t) +
\sum_{n \geq 2} H_n {\mathbf 1}_{(t_{n-1},t_n]}(t)$, $ t \in [t_0,T]$,
together with the tube
${\mathcal T} = \{ (t,\bar{e}') \in [t_0,T] \times \RR : |\bar{e}'-\bar{e}| < \psi_t\}$.

Then, there exists a constant $C$, depending on known parameters only, 
such that the distribution of the hitting time $\tau = \inf\{t \geq t_0 : (t,\bar{E}_t^{t_0,p,\bar{e}}) \not \in {\mathcal T}\}
\wedge T$ satisfies
\begin{equation}
\label{eq:5:8:2}
\forall n \geq 1, \quad {\mathbb P} \bigl\{ \tau \in (t_{n-1},t_n] \bigr\}
\leq C  2^{-6n} (T-t_0)^{3}.
\end{equation}
\end{lemma}

{\bf Proof.} 
From \eqref{eq:26:2:5}, we know that the drift of $(\bar{E}_t)_{t_0 \leq t \leq T}$
is in $[-L,0]$ and its diffusion coefficient less than $C(T-t_0)$. 
Given $n \geq 1$, 
\begin{equation*}
\begin{split}
&{\mathbb P} \bigl\{ \tau \in (t_{n-1},t_n] \bigr\}
\\
&\leq \sup_{p \in \RR^d,|\bar{e}'-\bar{e}| \leq H_{n-1}} {\mathbb P} \bigl\{
\sup_{t_{n-1} \leq t \leq t_n} |\bar{E}_t^{t_{n-1},p,\bar{e}'} - \bar{e}| \geq H_n \bigr\}.
\\
&\leq \sup_{p \in \RR^d,|\bar{e}'-\bar{e}| \leq H_{n-1}} {\mathbb P} \bigl\{
\sup_{t_{n-1} \leq t \leq t_n} |\bar{E}_t^{t_{n-1},p,\bar{e}'} - \bar{e}'| + H_{n-1} \geq H_n \bigr\}
\\
&= \sup_{p \in \RR^d, |\bar{e}'-\bar{e}| \leq H_{n-1}} {\mathbb P} \bigl\{
\sup_{t_{n-1} \leq t \leq t_n} \biggl|\int_{t_{n-1}}^t \langle \partial_p w(s,P_s),\sigma(P_s) dB_s \rangle \biggr| \geq h_{n} - L \delta_n
\geq \frac{L(T-t_0)}{2^n} \bigr\}.
\end{split}
\end{equation*}
By Markov and Doob inequalities (with 6 as exponent), we conclude that there exists a constant 
$C$, depending on $L$ and $T$ only, such that
\begin{equation*}
{\mathbb P} \bigl\{ \tau \in (t_{n-1},t_n] \bigr\}
\leq  C  2^{6n} \delta_n^3
= C (T-t_0)^3 2^{-6n}.\qed
\end{equation*}

\bibliographystyle{plain}

\end{document}